\documentclass[aps,prb]{revtex4}

\usepackage{amsmath,amssymb,amsthm}
\usepackage{mathtools}
\usepackage{xspace}
\usepackage[usenames,dvipsnames]{xcolor}
\usepackage[colorlinks,hyperindex,breaklinks]{hyperref}
\usepackage{braket}
\usepackage{hyphenat}
\usepackage{bbold}
\usepackage{enumerate}
\usepackage{tabularx}
\usepackage[hang,flushmargin]{footmisc}

\usepackage[british]{babel}

\usepackage{srcltx}
\usepackage{color}

\usepackage{graphicx}
\usepackage{amsmath}
\usepackage{amssymb}
\usepackage{bm}

\begin{document}

\title{On an infinite number of nonlinear Euler sums}
\author{J. Braun}
\affiliation{Ludwig Maximilians-Universit{\"a}t, M{\"u}nchen, Germany}
\author{D. Romberger}
\affiliation{Fakult\"at IV, Abt. BWL, Hochschule Hannover, Germany}
\author{H. J. Bentz}
\affiliation{Institut f\"{u}r Mathematik und Informatik, Universit\"{a}t Hildesheim, Germany}
\date{\today}

\begin{abstract}
Linear harmonic number sums had been studied by a variety of authors during the last centuries, but only few results are known about nonlinear Euler sums of quadratic
or even higher degree. The first systematic study on nonlinear Euler sums consisting of products of hyperharmonic sums had been published by Flajolet and Salvy in
1997 followed by similar studies presented during the last years by different authors. Although these studies had been restricted to sums where the nominator
consists of a product of even or odd hyperharmonic sums, where the denominator is of the type $1/k^n$. We have generalized these results to nonlinear Euler sums with
different denominators and nominators which consist in addition of mixed products between even and odd hyperharmonic numbers. In detail we present eight families of
quadratic Euler sums which are expressible by zeta values and special types of linear Euler sums only where the order of the nonlinear Euler sums is always an even
number. The resulting eight different families of nonlinear Euler sums which we discovered consist of various products between even and odd hyperharmonic numbers,
divided by three different types of denominators $1/k^n$, $1/((2k-1)^n)$ and $1/(k(2k-1))$. The calculational scheme is based on proper two-valued integer functions, which
allow us to compute these sequences explicitly in terms of zeta values and pairs of odd-type linear harmonic numbers and even hyperharmonic numbers of second order.
\end{abstract}
\maketitle
%Keywords: Harmonic Numbers, Nonlinear Euler sums, Riemann Zeta function

\section{Introduction}
Harmonic numbers and their generalizations called hyper-harmonic numbers are defined by
\begin{eqnarray}
H_k = H_k^{(1)} = \sum^{k}_{i=1} \frac{1}{i}~, ~~~~~~~~~~~~  H_k^{(n)} = \sum^{k}_{i=1} \frac{1}{i^n}~.
\end{eqnarray}
It has been discovered by many authors in the past starting with Goldbach 1742 and later on with Euler, Ramanujan, Borwein, Salvi and others that linear sums of the type
\begin{eqnarray}
t(n,m) = \sum^{\infty}_{k=1} \frac{H_k^{(n)}}{k^m} 
\end{eqnarray}
are expressible in terms of zeta values only for an odd weight $p=n+m$, where $p=n+m$ defines the order of the corresponding sum. An excellent overview about these works is
found in the paper of Flajolet and Salvy \cite{fla97} and in citations therein. Further work on similar types of Euler sums can be found, for example, in \cite{ce16,ce17,ce18,ce20}.
Although, only few sums of the type
\begin{eqnarray}
s(n,m) = \sum^{\infty}_{k=1} \frac{h_k^{(n)}}{k^m} 
\end{eqnarray}
had been discovered in the past \cite{sit85,zeh07,ade16,ade16a,ce20,bra21}. In analogy these sums are expressible in terms of zeta values only for an odd weight $p=n+m$, where again
$p=n+m$ defines the order of the corresponding sum. Here the symbol $h_k$ denote odd-type (hyper-)harmonic sums which are defined as follows:
\begin{eqnarray}
h_k = h_k^{(1)} = \sum^{k}_{i=1} \frac{1}{2i-1}~,~~~~~~~~~~~~  h_k^{(n)} = \sum^{k}_{i=1} \frac{1}{(2i-1)^n}~.
\end{eqnarray} 
This way the subject of the present paper is to present a formalism that allows for an explicit calculation of quadratic Euler sums in terms zeta values and linear Euler
sums. In detail we introduce a calculational scheme which is based on proper two-valued integer functions in correspondence to our last works on Euler sums that were devoted
to odd harmonic numbers and central binomial coefficients as well as to special types of nonlinear Euler sums \cite{bra20,bra21}. We present in the first section a special
type of Euler sums, which we define as the 'zero family'. We demonstrate then in the following sections that a variety of quadratic Euler sums can be calculated explicitly
in terms of these linear sums and zeta values only. In this sense our work can be seen as a generalization of the paper presented by Flajolet and Salvy \cite{fla97}, where
sums of products of even hyperharmonic numbers divided by $r^n$-type denominators had been discovered. 

\section{Zero family}

This way, We start with the so called zero family of linear Euler sums which are defined as:
\begin{eqnarray}
s(2m,2n) = \sum^{\infty}_{k=1} \frac{h^{(2m)}_k}{k^{2n}}~.
\end{eqnarray}~,
where $m, n \in \mathbb{N}$. As an example, the corresponding order 4 sum with $m=2$ and $ n=2$ is defined as:
\begin{eqnarray}
s(2,2) = \sum^{\infty}_{k=1} \frac{h^{(2)}_k}{k^2} = \frac{45}{16} \zeta(4) - \sum^{\infty}_{k=1} \frac{h_k}{k^3}~.
\end{eqnarray}
Related to these sums are the following ones:
\begin{eqnarray}
u(2,2) = \sum^{\infty}_{k=1} \frac{H^{(2)}_k}{(2k-1)^2} = -\frac{15}{16} \zeta(4) + 4\zeta(2) - 8ln(2) ) + \sum^{\infty}_{k=1} \frac{h_k}{k^3}
\end{eqnarray}
and
\begin{eqnarray}
v(2,2) = \sum^{\infty}_{k=1} \frac{h^{(2)}_k}{(2k-1)^2} = \frac{75}{64} \zeta(4)~.
\end{eqnarray}

For the following order 6 sum it holds:
\begin{eqnarray}
s(2,4) = \sum^{\infty}_{k=1} \frac{h^{(2)}_k}{k^4} = \frac{135}{24} \zeta(6) - \frac{7}{4}\zeta(3)^2 - 
2\sum^{\infty}_{k=1} \frac{h_k}{k^5}
\end{eqnarray}

The corresponding calculational procedure for these Euler sums is based on the following identity \cite{ade16a}: 
\begin{eqnarray}
H^{(n)}_{2k} - H^{(n)}_{k} =  \sum^{k}_{i=1} \frac{1}{(i+k)^n}~.
\end{eqnarray}
From this we get immediately
\begin{eqnarray}
h^{(n)}_{k} =  \sum^{k}_{i=1} \frac{1}{(i+k)^n} + \frac{2^n - 1}{2^n} H^{(n)}_{k}~.
\end{eqnarray}
Choosing n=2, dividing both sides by $k^4$ and summing up over k we arrive at:
\begin{eqnarray}
\sum^{\infty}_{k=1} \frac{h^{(2)}_k}{k^4} = \frac{3}{4} \frac{H^{(2)}_k}{k^4} + \sum^{\infty}_{k=1} \frac{1}{k^4} \left( \sum^{k}_{i=1} \frac{1}{(i+k)^2} \right)~.
\end{eqnarray}
Rearranging the double summation on the right side it follows:
\begin{eqnarray}
\sum^{\infty}_{k=1} \frac{h^{(2)}_k}{k^4} = \frac{1}{4}\zeta(6) + \frac{3}{4} \frac{H^{(2)}_k}{k^4} + \sum^{\infty}_{i=1}
\left( \sum^{\infty}_{k=1} \frac{1}{(k+i)^4(k+2i)^2} \right)~.
\end{eqnarray}
By partial fraction decomposition on gets:
\begin{eqnarray}
\sum^{\infty}_{k=1} \frac{1}{(k+i)^4(k+2i)^2} &=& \frac{1}{i^2}\sum^{\infty}_{k=1} \frac{1}{(k+i)^4} - \frac{2}{i^2}\sum^{\infty}_{k=1} \frac{1}{(k+i)^3} +
\frac{3}{i^4}\sum^{\infty}_{k=1} \frac{1}{(k+i)^2} - \frac{4}{i^2}\sum^{\infty}_{k=1} \frac{1}{(k+i)(k+2i)} \nonumber \\ &+&
\frac{1}{i^2}\sum^{\infty}_{k=1} \frac{1}{(k+2i)^2}~.
\end{eqnarray}
Using the three identities \cite{ade16a}:
\begin{eqnarray}
\sum^{\infty}_{k=1} \frac{1}{(k+i)^n} = \frac{2^n - 1}{2^n}\zeta(n) - H^{(n)}_{i}~,
\end{eqnarray}
\begin{eqnarray}
\sum^{\infty}_{k=1} \frac{1}{(k+i)(k+2i)} = \frac{H_{2i}-H_i}{i}
\end{eqnarray}
and
\begin{eqnarray}
\sum^{\infty}_{k=1} \frac{1}{(k+2i)^2} = \zeta(2) - H^{(2)}_{2i} 
\end{eqnarray}
it follows: 
\begin{eqnarray}
\sum^{\infty}_{k=1} \frac{h^{(2)}_k}{k^4} = \frac{9}{2}\zeta(6) - \zeta(3)^2 - \frac{7}{4}\sum^{\infty}_{k=1} \frac{H^{(2)}_k}{k^4} + \sum^{\infty}_{k=1}\frac{H^{(3)}_k}{k^3} +
\sum^{\infty}_{k=1} \frac{H_k}{k^5} -2\sum^{\infty}_{k=1}\frac{h_k}{k^5} = s(2,4)~.
\end{eqnarray}

For the remaining order 6 sums it follows:
\begin{eqnarray}
s(3,3) = \sum^{\infty}_{k=1} \frac{h^{(3)}_k}{k^3} = -\frac{135}{16} \zeta(6) + \frac{91}{16}\zeta(3)^2 + 
\frac{3}{2}\sum^{\infty}_{k=1} \frac{h_k}{k^5}
\end{eqnarray}
and
\begin{eqnarray}
s(4,2) = \sum^{\infty}_{k=1} \frac{h^{(4)}_k}{k^2} = \frac{495}{64} \zeta(6) - \frac{63}{16}\zeta(3)^2 - 
\frac{1}{2}\sum^{\infty}_{k=1} \frac{h_k}{k^5}~.
\end{eqnarray}

Related to these sums are the following ones:
\begin{eqnarray}
u(2,4) = \sum^{\infty}_{k=1} \frac{H^{(2)}_k}{(2k-1)^4} = -\frac{195}{32} \zeta(6) + \frac{63}{16}\zeta(3)^2 + \frac{15}{4}\zeta(4) -7\zeta(3) + 10\zeta(2) 
-16 ln(2)+ \frac{1}{2}\sum^{\infty}_{k=1} \frac{h_k}{k^5}~,
\end{eqnarray}
\begin{eqnarray}
u(4,2) = \sum^{\infty}_{k=1} \frac{H^{(4)}_k}{(2k-1)^2} = -\frac{207}{48} \zeta(6) + \frac{7}{4}\zeta(3)^2 + \zeta(4) +4\zeta(3) + 24\zeta(2) -64ln(2)
+ 2\sum^{\infty}_{k=1} \frac{h_k}{k^5}~, 
\end{eqnarray}
\begin{eqnarray}
v(2,4) = \sum^{\infty}_{k=1} \frac{h^{(2)}_k}{(2k-1)^4} = \frac{27}{32} \zeta(6) + \frac{7}{64}\zeta(3)^2~,
\end{eqnarray}
and
\begin{eqnarray}
v(4,2) = \sum^{\infty}_{k=1} \frac{h^{(4)}_k}{(2k-1)^2} = \frac{351}{256} \zeta(6) - \frac{7}{64}\zeta(3)^2.
\end{eqnarray}
As one can see all of these s-, t- and v-type Euler sums are expressible in terms of odd-type linear Euler sums and zeta values. \\

Again we apply this calculational scheme for the order 8 case. For example, choosing n=2, dividing both sides by $k^6$ and summing up over k one arrives, after
rearranging the double summation, at the following expression for $s(2,6)$:
\begin{eqnarray}
\sum^{\infty}_{k=1} \frac{h^{(2)}_k}{k^6} = \frac{1}{4}\zeta(8) + \frac{3}{4} \frac{H^{(2)}_k}{k^6} + \sum^{\infty}_{i=1}
\left( \sum^{\infty}_{k=1} \frac{1}{(k+i)^6(k+2i)^2} \right)~.
\end{eqnarray}
From this $s(2,6)$ explicitly results as:
\begin{eqnarray}
s(2,6) = \sum^{\infty}_{k=1} \frac{h^{(2)}_k}{k^6} = \frac{13}{2} \zeta(8) + \frac{3}{2}\zeta(3)\zeta(5) - 
3\sum^{\infty}_{k=1} \frac{h_k}{k^7} - \frac{17}{4} \sum^{\infty}_{k=1} \frac{H^{(2)}_k}{k^6}~,
\end{eqnarray}
although an additional linear Euler sum appears on the right side which is not expressible in terms of zeta values only. For the remaining Euler sums of order 8 it
follows:
\begin{eqnarray}
s(3,5) = \sum^{\infty}_{k=1} \frac{h^{(3)}_k}{k^5} = \frac{1073}{64} \zeta(8) - \frac{143}{4}\zeta(3)\zeta(5) + 
\frac{15}{4}\sum^{\infty}_{k=1} \frac{h_k}{k^7} + \frac{391}{16} \sum^{\infty}_{k=1} \frac{H^{(2)}_k}{k^6}~,
\end{eqnarray}
\begin{eqnarray}
s(4,4) = \sum^{\infty}_{k=1} \frac{h^{(4)}_k}{k^4} = -\frac{2661}{64} \zeta(8) + \frac{271}{4}\zeta(3)\zeta(5) - 
\frac{5}{2}\sum^{\infty}_{k=1} \frac{h_k}{k^7} - \frac{153}{4} \sum^{\infty}_{k=1} \frac{H^{(2)}_k}{k^6}~,
\end{eqnarray}
\begin{eqnarray}
s(5,3) = \sum^{\infty}_{k=1} \frac{h^{(5)}_k}{k^3} = \frac{5335}{256} \zeta(8) - \frac{1123}{32}\zeta(3)\zeta(5) + 
\frac{15}{16}\sum^{\infty}_{k=1} \frac{h_k}{k^7} + \frac{1445}{64} \sum^{\infty}_{k=1} \frac{H^{(2)}_k}{k^6}
\end{eqnarray}
and
\begin{eqnarray}
s(6,2) = \sum^{\infty}_{k=1} \frac{h^{(6)}_k}{k^2} = \frac{569}{128} \zeta(8) + \frac{51}{32}\zeta(3)\zeta(5) - 
\frac{3}{16}\sum^{\infty}_{k=1} \frac{h_k}{k^7} - \frac{289}{64} \sum^{\infty}_{k=1} \frac{H^{(2)}_k}{k^6}~.
\end{eqnarray}

Related to these sums are the following ones:
\begin{eqnarray}
u(6,2) = \sum^{\infty}_{k=1} \frac{H^{(6)}_k}{(2k-1)^2} &=& -\frac{21}{4}\zeta(8) - \frac{3}{2}\zeta(3)\zeta(5) + \zeta(6) + 4\zeta(5) + 12\zeta(4)
+ 32\zeta(3) + 128\zeta(2) - 384ln(2) \nonumber \\ &+& 3\sum^{\infty}_{k=1} \frac{h_k}{k^7} + \frac{17}{4}\sum^{\infty}_{k=1} \frac{H^{(2)}_k}{k^6}~,
\end{eqnarray}

\begin{eqnarray}
u(2,6) = \sum^{\infty}_{k=1} \frac{H^{(2)}_k}{(2k-1)^6} &=& -\frac{1077}{384}\zeta(8) + \frac{63}{16}\zeta(6) - \frac{51}{32}\zeta(3)\zeta(5) - \frac{31}{4}\zeta(5)
+ \frac{45}{4}\zeta(4) - 14\zeta(3) + 16\zeta(2) - 24ln(2) \nonumber \\ &+& \frac{3}{16}\sum^{\infty}_{k=1} \frac{h_k}{k^7}
+ \frac{289}{64}\sum^{\infty}_{k=1} \frac{H^{(2)}_k}{k^6}~,
\end{eqnarray}

\begin{eqnarray}
v(2,6) = \sum^{\infty}_{k=1} \frac{h^{(2)}_k}{(2k-1)^6} = \frac{307}{512}\zeta(8) + \frac{3}{8}\zeta(3)\zeta(5)- \frac{17}{256}\sum^{\infty}_{k=1} \frac{H^{(2)}_k}{k^6}~,
\end{eqnarray}

\begin{eqnarray}
v(3,5) = \sum^{\infty}_{k=1} \frac{h^{(3)}_k}{(2k-1)^5} = -\frac{1093}{2048}\zeta(8) + \frac{407}{256}\zeta(3)\zeta(5) - \frac{221}{512}\sum^{\infty}_{k=1} \frac{H^{(2)}_k}{k^6}~,
\end{eqnarray}

\begin{eqnarray}
v(4,4) = \sum^{\infty}_{k=1} \frac{h^{(4)}_k}{(2k-1)^4} = \frac{1035}{1024}\zeta(8)~,
\end{eqnarray}

\begin{eqnarray}
v(5,3) = \sum^{\infty}_{k=1} \frac{h^{(5)}_k}{(2k-1)^3} = \frac{3133}{2048}\zeta(8) - \frac{95}{128}\zeta(3)\zeta(5) + \frac{221}{512}\sum^{\infty}_{k=1} \frac{H^{(2)}_k}{k^6}~,
\end{eqnarray}

\begin{eqnarray}
v(6,2) = \sum^{\infty}_{k=1} \frac{h^{(6)}_k}{(2k-1)^2} = \frac{833}{512}\zeta(8) - \frac{3}{8}\zeta(3)\zeta(5)+ \frac{17}{256}\sum^{\infty}_{k=1} \frac{H^{(2)}_k}{k^6}~.
\end{eqnarray}

The calculational scheme introduced above can be extended to higher order Euler sums, where for each order two new linear Euler sums appear, which are not expressible
in terms of zeta values only. Exemplarily, this will be shown for the order 10 case, where
\begin{eqnarray}
\sum^{\infty}_{k=1} \frac{h_k}{k^9}~,~~~\sum^{\infty}_{k=1} \frac{H^{(2)}_k}{k^8} 
\end{eqnarray}
show up in addition.
\begin{eqnarray}
\sum^{\infty}_{k=1} \frac{h^{(2)}_k}{k^8} = \frac{1}{4}\zeta(10) + \frac{3}{4}\frac{H^{(2)}_k}{k^8} + \sum^{\infty}_{i=1}
\left( \sum^{\infty}_{k=1} \frac{1}{(k+i)^8(k+2i)^2} \right)~.
\end{eqnarray}
By partial fraction decomposition on gets:
\begin{eqnarray}
\sum^{\infty}_{k=1}\frac{h^{(2)}_k}{k^8} &=& \frac{239}{20}\zeta(10) - 4\zeta(3)\zeta(7) -2\zeta(5)^2 + 2\sum^{\infty}_{k=1}\frac{H_k}{k^9} - 4\sum^{\infty}_{k=1}\frac{h_k}{k^9}
-\frac{13}{4}\sum^{\infty}_{k=1}\frac{H^{(2)}_k}{k^8} + 3\sum^{\infty}_{k=1}\frac{H^{(3)}_k}{k^7} \nonumber \\ &-& \frac{5}{2}\sum^{\infty}_{k=1}\frac{H^{(4)}_k}{k^6} +
2\sum^{\infty}_{k=1}\frac{H^{(5)}_k}{k^5} - \frac{3}{2}\sum^{\infty}_{k=1}\frac{H^{(6)}_k}{k^4} + \sum^{\infty}_{k=1}\frac{H^{(7)}_k}{k^3} -
\frac{1}{2}\sum^{\infty}_{k=1}\frac{H^{(8)}_k}{k^2}~.
\end{eqnarray}
With
\begin{eqnarray}
\sum^{\infty}_{k=1}\frac{H_k}{k^9} = \frac{11}{4}\zeta(10) - \zeta(3)\zeta(7) - \frac{1}{2}\zeta(5)^2~,
\end{eqnarray}
\begin{eqnarray}
\sum^{\infty}_{k=1}\frac{H^{(5)}_k}{k^5} = \frac{1}{2}\zeta(10) + \frac{1}{2}\zeta(5)^2~,
\end{eqnarray}
\begin{eqnarray}
\sum^{\infty}_{k=1}\frac{H^{(6)}_k}{k^4} = \frac{21}{10}\zeta(10) - \sum^{\infty}_{k=1}\frac{H^{(4)}_k}{k^6}~,
\end{eqnarray}
\begin{eqnarray}
\sum^{\infty}_{k=1}\frac{H^{(7)}_k}{k^3} = \zeta(10) + \zeta(3)\zeta(7) - \sum^{\infty}_{k=1}\frac{H^{(3)}_k}{k^7}~,
\end{eqnarray}
\begin{eqnarray}
\sum^{\infty}_{k=1}\frac{H^{(8)}_k}{k^2} = \frac{53}{20}\zeta(10) - \sum^{\infty}_{k=1}\frac{H^{(2)}_k}{k^8}~,
\end{eqnarray}
\begin{eqnarray}
\sum^{\infty}_{k=1}\frac{H^{(3)}_k}{k^7} = -\frac{33}{4}\zeta(10) + 7\zeta(3)\zeta(7) + 4\zeta(5)^2 -\frac{7}{2}\sum^{\infty}_{k=1}\frac{H^{(2)}_k}{k^8}
\end{eqnarray}
and
\begin{eqnarray}
\sum^{\infty}_{k=1}\frac{H^{(4)}_k}{k^6} = \frac{227}{20}\zeta(10) - 7\zeta(3)\zeta(7) - 5\zeta(5)^2 +\frac{7}{2}\sum^{\infty}_{k=1}\frac{H^{(2)}_k}{k^8}
\end{eqnarray}
it follows:
\begin{eqnarray}
s(2,8) = \sum^{\infty}_{k=1}\frac{h^{(2)}_k}{k^8} = -\frac{103}{8}\zeta(10) + 16\zeta(3)\zeta(7) + 11\zeta(5)^2 -4\sum^{\infty}_{k=1}\frac{h_k}{k^9} -
\frac{53}{4}\sum^{\infty}_{k=1}\frac{H^{(2)}_k}{k^8}~.
\end{eqnarray}
All other s-type sums of order 10 can be obtained successively by applying the following scheme. In a second step we calculate $s(8,2)$ as a function of $s(2,8)$.
First, it follows:
\begin{eqnarray}
\sum^{\infty}_{k=1}\frac{H^{(2)}_k}{k^8} = \frac{1}{256}\sum^{\infty}_{k=1}\frac{H^{(2)}_{2k}}{k^8} + \sum^{\infty}_{k=1}\frac{H^{(2)}_{2k-1}}{(2k-1)^8}~.
\end{eqnarray}
From this we get:
\begin{eqnarray}
\frac{1023}{1024}\sum^{\infty}_{k=1}\frac{H^{(2)}_k}{k^8} &=& \frac{1}{256}\sum^{\infty}_{k=1}\frac{h^{(2)}_{k}}{k^8} + \frac{1}{4} 
\sum^{\infty}_{k=1}\frac{H^{(2)}_k}{(2k-1)^8} - \frac{1}{4}\sum^{\infty}_{k=1}\frac{1}{k^2(2k-1)^8} + \sum^{\infty}_{k=1}\frac{h^{(2)}_k}{(2k-1)^8}~.
\end{eqnarray}
Furthermore, it follows from a rearrangement of the double summation: 
\begin{eqnarray}
\sum^{\infty}_{k=1}\frac{H^{(2)}_k}{(2k-1)^8} = \sum^{\infty}_{k=1}\frac{1}{(2k-1)^8} \left( \sum^{k}_{i=1}\frac{1}{i^2} \right) = \frac{255}{256}\zeta(2)\zeta(8)
+ \sum^{\infty}_{k=1}\frac{1}{k^2(2k-1)^8} - \sum^{\infty}_{k=1}\frac{h^{(8)}_{k}}{k^2}~. 
\end{eqnarray}
Using this identity we get:
\begin{eqnarray}
\frac{1}{4}\sum^{\infty}_{k=1}\frac{h^{(8)}_{k}}{k^2} = \frac{255}{1024}\zeta(2)\zeta(8) + \frac{1}{256}\sum^{\infty}_{k=1}\frac{h^{(2)}_{k}}{k^8} - 
\frac{1023}{1024}\sum^{\infty}_{k=1}\frac{H^{(2)}_k}{k^8} + \sum^{\infty}_{k=1}\frac{h^{(2)}_k}{(2k-1)^8}~.
\end{eqnarray}
The last sum on the right side can be calculated by use of the following identity:
\begin{eqnarray}
\sum^{\infty}_{k=1}\frac{1}{(2k-1)^8} \left( \sum^{\infty}_{i=1} \frac{1}{(2i+2k-1)^2} \right) = \frac{765}{1024}\zeta(2)\zeta(8) - 
\sum^{\infty}_{k=1}\frac{h^{(2)}_k}{(2k-1)^8}~.
\end{eqnarray}
Again, by rearranging the double summation on the left side and by corresponding partial fraction decomposition we get:
\begin{eqnarray}
\sum^{\infty}_{k=1}\frac{h^{(2)}_k}{(2k-1)^8} = \frac{486}{1024}\zeta(2)\zeta(8) + \frac{148}{512}\zeta(3)\zeta(7) - \frac{33}{128}\zeta(4)\zeta(6) +
\frac{31}{256}\zeta(5)^2 + \frac{1}{64}\sum^{\infty}_{k=1}\frac{h_k}{k^9} + \frac{1}{256}\sum^{\infty}_{k=1}\frac{h^{(2)}_{k}}{k^8}~.
\end{eqnarray}
This way it follows finally:
\begin{eqnarray}
\sum^{\infty}_{k=1}\frac{h^{(8)}_{k}}{k^2} = \frac{3317}{1024}\zeta(10) + \frac{53}{32}\zeta(3)\zeta(7) + \frac{53}{64}\zeta(5)^2 -
\frac{1}{16}\sum^{\infty}_{k=1}\frac{h_k}{k^9} - \frac{1129}{256}\sum^{\infty}_{k=1}\frac{H^{(2)}_k}{k^8}~.
\end{eqnarray}
Analogously we get:
\begin{eqnarray}
\sum^{\infty}_{k=1}\frac{h^{(7)}_{k}}{k^3} = \frac{127}{128}\zeta(3)\zeta(7) + \frac{1}{16}\sum^{\infty}_{k=1}\frac{h^{(3)}_k}{k^7} - 
\frac{1023}{128}\sum^{\infty}_{k=1}\frac{H^{(3)}_k}{k^7} + 8\sum^{\infty}_{k=1}\frac{h^{(3)}_k}{(2k-1)^7}~.
\end{eqnarray}
Using the identity
\begin{eqnarray}
\sum^{\infty}_{k=1}\frac{h^{(3)}_k}{(2k-1)^7} = \frac{889}{1024}\zeta(3)\zeta(7) - \sum^{\infty}_{k=1}\frac{1}{(2k-1)^7} \left( 
\sum^{\infty}_{i=1} \frac{1}{(2i+2k-1)^3} \right)
\end{eqnarray}
it follows again by rearranging the double summation on the left side and by corresponding partial fraction decomposition:
\begin{eqnarray}
\sum^{\infty}_{k=1}\frac{h^{(3)}_k}{(2k-1)^7} = \frac{1023}{2048}\zeta(10) + \frac{83}{128}\zeta(3)\zeta(7) - \frac{93}{512}\zeta(5)^2 -
\frac{7}{128}\sum^{\infty}_{k=1}\frac{h_k}{k^9} - \frac{7}{256}\sum^{\infty}_{k=1}\frac{h^{(2)}_{k}}{k^8} - 
\frac{1}{128}\sum^{\infty}_{k=1}\frac{h^{(3)}_{k}}{k^7}~.
\end{eqnarray}
From this it follows:
\begin{eqnarray}
\sum^{\infty}_{k=1}\frac{h^{(7)}_{k}}{k^3} = \frac{37247}{512}\zeta(10) - \frac{3409}{64}\zeta(3)\zeta(7) - \frac{2293}{64}\zeta(5)^2 +
\frac{7}{16}\sum^{\infty}_{k=1}\frac{h_k}{k^9} + \frac{7903}{256}\sum^{\infty}_{k=1}\frac{H^{(2)}_{k}}{k^8}~.
\end{eqnarray}
Furthermore, using the identity:
\begin{eqnarray}
\sum^{\infty}_{k=1}\frac{h^{(4)}_{k}}{k^6} = \frac{1}{16}\zeta(10) + \frac{15}{16}\sum^{\infty}_{k=1}\frac{H^{(4)}_{k}}{k^6} +
\sum^{\infty}_{i=1} \left( \sum^{\infty}_{k=1} \frac{1}{(k+i)^6(2i+2k)^4} \right)
\end{eqnarray}
we get
\begin{eqnarray}
\sum^{\infty}_{k=1}\frac{h^{(4)}_{k}}{k^6} &=& \frac{6461}{64}\zeta(10) - \frac{1491}{16}\zeta(3)\zeta(7) - \frac{1035}{16}\zeta(5)^2 +
14\sum^{\infty}_{k=1}\frac{h_k}{k^9} + \frac{2275}{32}\sum^{\infty}_{k=1}\frac{H^{(2)}_{k}}{k^8} -3\sum^{\infty}_{k=1}\frac{h^{(3)}_{k}}{k^7}
\nonumber \\ &=&
-\frac{28051}{64}\zeta(10) + \frac{4365}{16}\zeta(3)\zeta(7) + \frac{3951}{16}\zeta(5)^2 - 7\sum^{\infty}_{k=1}\frac{h_k}{k^9} - \frac{4757}{32}
\sum^{\infty}_{k=1}\frac{H^{(2)}_k}{k^8}~.
\end{eqnarray}
Analogously it follows:
\begin{eqnarray}
\sum^{\infty}_{k=1}\frac{h^{(6)}_{k}}{k^4} &=& -\frac{35525}{256}\zeta(10) + \frac{5621}{64}\zeta(3)\zeta(7) + \frac{4093}{64}\zeta(5)^2 +
\frac{7}{2}\sum^{\infty}_{k=1}\frac{h_k}{k^9} - \frac{4837}{128}\sum^{\infty}_{k=1}\frac{H^{(2)}_{k}}{k^8} - \frac{3}{4}\sum^{\infty}_{k=1}\frac{h^{(3)}_{k}}{k^7}
\nonumber \\ &=&
-\frac{70037}{256}\zeta(10) + \frac{11477}{64}\zeta(3)\zeta(7) + \frac{9079}{64}\zeta(5)^2 - \frac{7}{4}\sum^{\infty}_{k=1}\frac{h_k}{k^9} - \frac{11869}{128}
\sum^{\infty}_{k=1}\frac{H^{(2)}_k}{k^8}~.
\end{eqnarray}
and
\begin{eqnarray}
\sum^{\infty}_{k=1}\frac{h^{(5)}_{k}}{k^5} &=& -\frac{25095}{128}\zeta(10) + \frac{5215}{32}\zeta(3)\zeta(7) + \frac{8231}{64}\zeta(5)^2 -
\frac{175}{8}\sum^{\infty}_{k=1}\frac{h_k}{k^9} - \frac{7665}{64}\sum^{\infty}_{k=1}\frac{H^{(2)}_{k}}{k^8} + \frac{15}{4}\sum^{\infty}_{k=1}\frac{h^{(3)}_{k}}{k^7}~,
\nonumber \\ &=&
\frac{61185}{128}\zeta(10) - \frac{9425}{32}\zeta(3)\zeta(7) - \frac{16699}{64}\zeta(5)^2 + \frac{35}{8}\sum^{\infty}_{k=1}\frac{h_k}{k^9} + \frac{9915}{64}
\sum^{\infty}_{k=1}\frac{H^{(2)}_k}{k^8}~,
\end{eqnarray}
where the last expression follows by use of the identity:
\begin{eqnarray}
\sum^{\infty}_{k=1}\frac{h^{(5)}_{k}}{k^5} = \frac{1}{32}\zeta(10) + \frac{31}{32}\sum^{\infty}_{k=1}\frac{H^{(5)}_{k}}{k^5} +
\sum^{\infty}_{i=1} \left( \sum^{\infty}_{k=1} \frac{1}{(k+i)^5(2i+2k)^5} \right)~.
\end{eqnarray}
Finally, using the identity:
\begin{eqnarray}
\sum^{\infty}_{k=1}\frac{h^{(7)}_k}{(2k-1)^3} = \frac{889}{1024}\zeta(3)\zeta(7) - \sum^{\infty}_{k=1}\frac{1}{(2k-1)^3} 
\left( \sum^{\infty}_{i=1} \frac{1}{(2i+2k-1)^7} \right)
\end{eqnarray}
we get:
\begin{eqnarray}
\frac{889}{1024}\zeta(3)\zeta(7) = \frac{7}{64}\sum^{\infty}_{k=1}\frac{h_k}{k^9} + \frac{7}{64}\sum^{\infty}_{k=1}\frac{h^{(2)}_{k}}{k^8} +
\frac{1}{8}\sum^{\infty}_{k=1}\frac{h^{(3)}_{k}}{k^7} + \frac{5}{32}\sum^{\infty}_{k=1}\frac{h^{(4)}_{k}}{k^6} + \frac{3}{16}\sum^{\infty}_{k=1}\frac{h^{(5)}_{k}}{k^5}
+ \frac{3}{16}\sum^{\infty}_{k=1}\frac{h^{(6)}_{k}}{k^4} + \frac{1}{8}\sum^{\infty}_{k=1}\frac{h^{(7)}_{k}}{k^3}~. 
\end{eqnarray}
Using the corresponding expressions for all the sums on the right side of the last identity one arrives at:
\begin{eqnarray}
\sum^{\infty}_{k=1}\frac{h^{(3)}_{k}}{k^7} = \frac{719}{4}\zeta(10) - 122\zeta(3)\zeta(7) - \frac{831}{8}\zeta(5)^2 + 7\sum^{\infty}_{k=1}\frac{h_k}{k^9} +
\frac{293}{4}\sum^{\infty}_{k=1}\frac{H^{(2)}_{k}}{k^8}~.
\end{eqnarray}
Therefore, as in the order 6 and 8 case the s-type linear Euler sums of order 10 can be written in terms of zeta values and the two sums s(1,9) and H(2,8). 
The order 10 series can be found in the appendix. 

All these Euler sums have been numerically calculated and compared to the corresponding analytical expressions to an accuracy of $10^{-16}$. 

\section{New Quadratic Euler sums}
In this section we will show that the following families of quadratic Euler sums can be expressed in terms of the linear Euler sums introduced in the last chapter,
where both the odd and even order sums are needed. It has to be mention here that only the linear Euler sums of odd oder are expressible in terms of zeta values only.
This means that the even order linear Euler sums $t(m,n)$ and $s(m,n)$ with $n+m=$ even can be seen as basic numbers like $ln(2), \zeta(3), Li_4(1/2)$ and so on
\cite{rem99}. As a remark we want to mention here that the order $p$ of the corresponding linear Euler sums has to be clearly distinguished from the degree
of the corresponding nonlinear Euler sums which is always $2$, where the attributed order is always even.

\section{First family}

The following type of nonlinear Euler sums can be expressed in terms of linear Euler sums and zeta values:
\begin{eqnarray}
\sum^{\infty}_{k=1} \frac{H^{(a)}_k h^{(b)}_k}{k^c}
\end{eqnarray}
with $a,b,c \in \mathbb{N}$, $a+b+c = $even. The corresponding calculational scheme is recursive as the different Euler sums depend on each other.
This will be shown in the following. We start with the order 4 case which is the most simple one. From Eq.~(42) defined in \cite{bra20} it follows:
\begin{eqnarray}
\sum^{\infty}_{k=1,i\ne k} \frac{1}{k(k-i)} = \frac{2}{i^2} - \frac{H_i}{i}~.
\end{eqnarray}
Multiplying both sides by $h_i$, dividing by i und summing up over i we get:
\begin{eqnarray}
\sum^{\infty}_{i=1} \frac{h_i}{i}\left( \sum^{\infty}_{k=1,i\ne k} \frac{1}{k(k-i)} \right) = 2\sum^{\infty}_{i=1} \frac{h_i}{i^3} - \sum^{\infty}_{i=1} \frac{H_i h_i}{i^2}~.
\end{eqnarray}
Rearranging the double summation it follows:
\begin{eqnarray}
\sum^{\infty}_{k=1} \frac{1}{k}\left( \sum^{\infty}_{i=1,k\ne i} \frac{h_i}{i(k-i)} \right) = 2\sum^{\infty}_{i=1} \frac{h_i}{i^3} - \sum^{\infty}_{i=1} \frac{H_i h_i}{i^2}~.
\end{eqnarray}
By use of identity (88) from \cite{bra20} one arrives at: 
\begin{eqnarray}
2\sum^{\infty}_{k=1} \frac{h_k}{k^3} - \sum^{\infty}_{k=1} \frac{H_k h_k}{k^2} = \sum^{\infty}_{k=1} \frac{1}{k} \left( \sum^{\infty}_{i=1} \frac{h_i}{i(k+i)}
- \frac{h_k}{k^2} - 2\frac{h^{(2)}_k}{k} \right)~.
\end{eqnarray}
From this we get:
\begin{eqnarray}
2\sum^{\infty}_{k=1} \frac{h_k}{k^3} - \sum^{\infty}_{k=1} \frac{H_k h_k}{k^2} = \sum^{\infty}_{i=1} \frac{h_i}{i} \left( \sum^{\infty}_{k=1} \frac{1}{k(i+k)} \right)
- \sum^{\infty}_{k=1} \frac{h_k}{k^3} - 2\sum^{\infty}_{k=1} \frac{h^{(2)}_k}{k^2}
\end{eqnarray}
and
\begin{eqnarray}
2\sum^{\infty}_{k=1} \frac{h_k}{k^3} - \sum^{\infty}_{k=1} \frac{H_k h_k}{k^2} = \sum^{\infty}_{k=1} \frac{H_k h_k}{k^2}
-\sum^{\infty}_{k=1}  \frac{h_k}{k^3} - 2\sum^{\infty}_{k=1} \frac{h^{(2)}_k}{k^2}
\end{eqnarray}
or
\begin{eqnarray}
\sum^{\infty}_{k=1} \frac{H_k h_k}{k^2} = \frac{1}{2}\sum^{\infty}_{k=1} \frac{h_k}{k^3} -\sum^{\infty}_{k=1} \frac{h^{(2)}_k}{k^2}~.
\end{eqnarray}
Using Eq.~(6) we finally arrive at Eq.~(75). This calculational scheme can be applied to all higher order Euler sums of this special type with $a=b=1$. 
For an explicit calculation the corresponding Euler sums with $a=1$ and $b>1$ one can use the identities (88) and (89) from \cite{bra20}. A calculational example
is given below for the order 8 case. In order to calculate explicitly the Euler sums with $a>1$ and $b=1$ for a given order we need a new class of two-valued integer 
functions:

The following identity holds:
\subsection{Lemma 1}
\begin{eqnarray}
\sum^{\infty}_{k=1,k\ne i} \frac{H_k}{k(i-k)} = \sum^{\infty}_{k=1} \frac{H_k}{k(i+k)} -\zeta(2)\frac{1}{k} - 2\frac{H^{(2)}_k}{k}~.
\end{eqnarray}

\subsection{proof}
For $i=1$ it follows:
\begin{eqnarray}
\sum^{\infty}_{k=2}\frac{H_k}{k(1-k)} = \sum^{\infty}_{k=2}\frac{H_k}{k(1-k)} = -\sum^{\infty}_{k=2}\frac{H_k}{k(k-1)} = - \sum^{\infty}_{k=1}\frac{H_{k+1}}{k(k+1)} = -2~.
\end{eqnarray}
In analogy we get for $i=2$ and $i=3$:
\begin{eqnarray}
\sum^{\infty}_{k=1,k\ne 2}\frac{H_k}{k(2-k)} = H_1 + \sum^{\infty}_{k=3}\frac{H_k}{k(2-k)} = H_1 - \sum^{\infty}_{k=3}\frac{H_k}{k(k-2)} = 
H_1 - \sum^{\infty}_{k=1}\frac{H_{k+2}}{k(k+2)}
\end{eqnarray}
and
\begin{eqnarray}
\sum^{\infty}_{k=1,k\ne 3}\frac{H_k}{k(3-k)} = \sum^{2}_{k=1} \frac{H_k}{k(3-k)} - \sum^{\infty}_{k=4}\frac{H_k}{k(3-k)} = \sum^{2}_{k=1} \frac{H_k}{k(3-k)}
- \sum^{\infty}_{k=1}\frac{H_{k+3}}{k(k+3)}~.
\end{eqnarray}
This way it follows for all $i \in \mathbb{N}$:
\begin{eqnarray}
\sum^{\infty}_{k=1,k\ne i} \frac{H_k}{k(i-k)} = \sum^{i-1}_{k=1} \frac{H_k}{k(i-k)} - \sum^{\infty}_{k=1}\frac{H_{i+k}}{k(i+k)}~.
\end{eqnarray}
Using the identities \cite{ce16}:
\begin{eqnarray}
\sum^{i-1}_{k=1} \frac{H_k}{(i-k)} = H^2_i - H^{(2)}_i
\end{eqnarray}
and
\begin{eqnarray}
\sum^{\infty}_{k=1}\frac{H_{i+k}}{k(i+k)} = \frac{1}{i} \left( H^{(2)}_i + H^2_i \right)~, 
\end{eqnarray}
one arrives at 
\begin{eqnarray}
\sum^{\infty}_{k=1,k\ne i} \frac{H_k}{k(i-k)} &=& \frac{1}{i} \sum^{i-1}_{k=1} \frac{H_k}{k} + \frac{1}{i} \sum^{i-1}_{k=1} \frac{H_k}{i-k} - \sum^{\infty}_{k=1}\frac{H_{i+k}}{k(i+k)}
\nonumber \\ &=& \frac{H^{(2)}_i}{2i} - \frac{1}{2i^3} - 2\frac{H^{(2)}_i}{i} + \frac{1}{2i} \left( H_i - \frac{1}{i} \right)^2
\nonumber \\ &=& \frac{H^2_i}{2i} - \frac{3}{2}\frac{H^{(2)}_i}{i} - \frac{H_i}{i^2}~.
\end{eqnarray}
Using the identity \cite{ce16}:
\begin{eqnarray}
\sum^{\infty}_{k=1} \frac{H_k}{k(i+k)} = \frac{1}{2i}H^{(2)}_i + \frac{1}{2i}H^2_i - \frac{H_i}{i^2} + \zeta(2)\frac{1}{i}
\end{eqnarray}
and calculating the corresponding difference between Eq.~(83) and Eq.~(84) one arrives at Eq.~(76). Thus the lemma is proved. 

These calculational scheme can be easily applied to $n>1$ This way, the following identities hold:
\begin{eqnarray}
\sum^{\infty}_{k=1,k\ne i}\frac{H^{(2n-1)}_k}{k(i-k)}&=&\sum^{\infty}_{k=1}\frac{H^{(2n-1)}_k}{k(i+k)}+\sum^{2n-1}_{m=1}\frac{(-)^m}{i^m}\zeta(2n-m+1)+
\frac{H_i}{i^{2n}}+2\sum^{n-1}_{m=1}\zeta(2n-2m)\frac{H^{(2m)}_i}{i} \nonumber \\ &-& \frac{H^{(2n-1)}_i}{i^2}-\frac{2nH^{(2n)}_i}{i}
\end{eqnarray}
for $n>1$ and 
\begin{eqnarray}
\sum^{\infty}_{k=1,k\ne i}\frac{H^{(2n)}_k}{k(i-k)}&=& -\sum^{\infty}_{k=1}\frac{H^{(2n)}_k}{k(i+k)}-\sum^{2n}_{m=1}\frac{(-)^m}{i^m}\zeta(2n-2m+2)+
\frac{H_i}{i^{2n+1}}+2\sum^{n}_{m=1}\zeta(2n-2m+2)\frac{H^{(2m-1)}_i}{i} \nonumber \\ &-& \frac{H^{(2n)}_i}{i^2}-\frac{(2n+1)H^{(2n+1)}_i}{i}
\end{eqnarray}
for $n\ge 1$.

\subsection{examples}
The order 4 case:
\begin{eqnarray}
\sum^{\infty}_{k=1} \frac{H_k h_k}{k^2} = \frac{45}{16} \zeta(4) + \frac{1}{2}\sum^{\infty}_{k=1} \frac{h_k}{k^3}~. 
\end{eqnarray}

The order 6 case:
\begin{eqnarray}
\sum^{\infty}_{k=1} \frac{H_k h_k}{k^4} = \frac{135}{24} \zeta(6) - \frac{7}{2}\zeta(3)^2
+\frac{1}{2}\sum^{\infty}_{k=1} \frac{h_k}{k^5}~,
\end{eqnarray}

\begin{eqnarray}
\sum^{\infty}_{k=1} \frac{H_k h^{(2)}_k}{k^3} = -\frac{135}{24}\zeta(6) + \frac{63}{8}\zeta(3)^2 - \frac{3}{2}\zeta(2)
\sum^{\infty}_{k=1} \frac{h_k}{k^3} - \sum^{\infty}_{k=1} \frac{h_k}{k^5}~,
\end{eqnarray}

\begin{eqnarray}
\sum^{\infty}_{k=1} \frac{H^{(2)}_k h_k}{k^3} = -\frac{7}{8}\zeta(3)^2 + \zeta(2)\sum^{\infty}_{k=1} \frac{h_k}{k^3} + \frac{1}{2}\sum^{\infty}_{k=1} \frac{h_k}{k^5}~,
\end{eqnarray}

\begin{eqnarray}
\sum^{\infty}_{k=1} \frac{H_k h^{(3)}_k}{k^2} = \frac{405}{128}\zeta(6) - \frac{105}{32}\zeta(3)^2 + \frac{3}{2}\zeta(2)
\sum^{\infty}_{k=1} \frac{h_k}{k^3} + \frac{3}{4} \sum^{\infty}_{k=1} \frac{h_k}{k^5}~,
\end{eqnarray}

\begin{eqnarray}
\sum^{\infty}_{k=1} \frac{H^{(2)}_k h^{(2)}_k}{k^2} = \frac{45}{32}\zeta(6) + \frac{21}{8}\zeta(3)^2 - \zeta(2)
\sum^{\infty}_{k=1} \frac{h_k}{k^3} - \sum^{\infty}_{k=1} \frac{h_k}{k^5}~
\end{eqnarray}
and
\begin{eqnarray}
\sum^{\infty}_{k=1} \frac{H^{(3)}_k h_k}{k^2} = \frac{135}{32}\zeta(6) - \frac{7}{4}\zeta(3)^2 + \frac{1}{2}\sum^{\infty}_{k=1} \frac{h_k}{k^5}~.
\end{eqnarray}

The order 8 case: \\

As an example we show the calculational scheme for the following order 8 Euler sum: 
\begin{eqnarray}
\sum^{\infty}_{k=1} \frac{H_k h^{(3)}_k}{k^4} &=& -\frac{225}{16}\zeta(8) - \frac{93}{8}\zeta(3)\zeta(5) + \frac{63}{8}\zeta(2)\zeta(3)^2 +
3\zeta(2)\sum^{\infty}_{k=1} \frac{h_k}{k^5} + \frac{5}{2}\sum^{\infty}_{k=1} \frac{h^{(3)}_k}{k^5} + 3\sum^{\infty}_{k=1} \frac{h^{(4)}_k}{k^4}~.
\end{eqnarray}
With the help of identity (88) from \cite{bra20} we have: 
\begin{eqnarray}
\sum^{\infty}_{i=1} \frac{1}{i^3}\left( \sum^{\infty}_{k=1,k\ne i} \frac{h^{(3)}_k}{k(i-k)} \right) = 
\sum^{\infty}_{i=1} \frac{1}{i^3}\left( \sum^{\infty}_{k=1} \frac{h^{(3)}_k}{k(i+k)}  \right) + 3\zeta(2)\sum^{\infty}_{i=1} \frac{h^{(2)}_i}{i^4} -
\sum^{\infty}_{i=1} \frac{h^{(3)}_i}{i^5} - 6\sum^{\infty}_{i=1} \frac{h^{(4)}_i}{i^4}~.
\end{eqnarray}
Rearranging the double summation we get:
\begin{eqnarray}
\sum^{\infty}_{k=1} \frac{h^{(3)}_k}{k} \left( \sum^{\infty}_{i=1,i\ne k} \frac{1}{i^3(i-k)} \right) = 
\sum^{\infty}_{k=1} \frac{h^{(3)}_k}{k} \left( \sum^{\infty}_{i=1} \frac{1}{i^3(i+k)}  \right) + 3\zeta(2)\sum^{\infty}_{k=1} \frac{h^{(2)}_k}{k^4} -
\sum^{\infty}_{i=1} \frac{h^{(3)}_k}{k^5} - 6\sum^{\infty}_{k=1} \frac{h^{(4)}_k}{k^4}~.
\end{eqnarray}
Using identity (42) from \cite{bra20} it follows:
\begin{eqnarray}
\sum^{\infty}_{k=1} \frac{h^{(3)}_k}{k} \left( \frac{4}{k^4} -\zeta(3)\frac{1}{k} -\zeta(2)\frac{1}{k^2} - \frac{H_k}{k^3} \right) &=& \zeta(3)
\sum^{\infty}_{k=1} \frac{h^{(3)}_k}{k^2} - \zeta(2)\sum^{\infty}_{k=1} \frac{h^{(3)}_k}{k^3} + \sum^{\infty}_{k=1}\frac{H_k h^{(3)}_k}{k^4} +
3\zeta(2)\sum^{\infty}_{k=1} \frac{h^{(2)}_k}{k^4} \nonumber \\ &-& \sum^{\infty}_{k=1} \frac{h^{(3)}_k}{k^5} - 6\sum^{\infty}_{k=1} \frac{h^{(4)}_k}{k^4}~.
\end{eqnarray}
By some simple algebraic manipulations it finally follows:
\begin{eqnarray}
\sum^{\infty}_{k=1}\frac{H_k h^{(3)}_k}{k^4} &=& \frac{5}{2}\sum^{\infty}_{k=1} \frac{h^{(3)}_k}{k^5} - \frac{3}{2}\zeta(2)\sum^{\infty}_{k=1} \frac{h^{(2)}_k}{k^4}
-\zeta(3)\sum^{\infty}_{k=1} \frac{h^{(3)}_k}{k^2} + 3\sum^{\infty}_{k=1} \frac{h^{(4)}_k}{k^4} \nonumber \\ &=&
-\frac{12401}{128}\zeta(8) + \frac{409}{4}\zeta(3)\zeta(5) + \frac{63}{8}\zeta(2)\zeta(3)^2 + 3\zeta(2)\sum^{\infty}_{k=1} \frac{h_k}{k^5} + 
\frac{15}{8}\sum^{\infty}_{k=1} \frac{h_k}{k^7} - \frac{1717}{32}\sum^{\infty}_{k=1} \frac{H^{(2)}_k}{k^6}~.
\end{eqnarray}
This expression is identical to Eq.~(94).

More example are given below.
\begin{eqnarray}
\sum^{\infty}_{k=1} \frac{H_k h_k}{k^6} = - \frac{19}{2}\zeta(3)\zeta(5) + \frac{7}{2}\zeta(2)\zeta(3)^2 +\frac{7}{2}\sum^{\infty}_{k=1} \frac{h_k}{k^7} 
+\sum^{\infty}_{k=1} \frac{h^{(2)}_k}{k^6}~, 
\end{eqnarray}

\begin{eqnarray}
\sum^{\infty}_{k=1} \frac{H_k h^{(2)}_k}{k^5} &=& 3\sum^{\infty}_{k=1} \frac{h^{(2)}_k}{k^6} - \zeta(3)\sum^{\infty}_{k=1} \frac{h^{(2)}_k}{k^3}
-\frac{3}{2}\zeta(2)\sum^{\infty}_{k=1} \frac{h_k}{k^5} + 2\sum^{\infty}_{k=1} \frac{h^{(3)}_k}{k^5} \nonumber \\ &=&
\frac{1697}{32}\zeta(8) - \frac{103}{2}\zeta(3)\zeta(5) - \frac{35}{4}\zeta(2)\zeta(3)^2 - \frac{3}{2}\zeta(2)\sum^{\infty}_{k=1} \frac{h_k}{k^5} - 
\frac{3}{2}\sum^{\infty}_{k=1} \frac{h_k}{k^7} + \frac{289}{8}\sum^{\infty}_{k=1} \frac{H^{(2)}_k}{k^6}~.
\end{eqnarray}

\begin{eqnarray}
\sum^{\infty}_{k=1} \frac{H^{(2)}_k h_k}{k^5} &=& \frac{97}{4}\zeta(3)\zeta(5) - 7\zeta(2)\zeta(3)^2 + \zeta(2)\sum^{\infty}_{k=1} \frac{h_k}{k^5}
- 7\sum^{\infty}_{k=1} \frac{h_k}{k^7} - 5\sum^{\infty}_{k=1} \frac{h^{(2)}_k}{k^6} - 2\sum^{\infty}_{k=1} \frac{h^{(3)}_k}{k^5} \nonumber \\ &=&
-\frac{2113}{32}\zeta(8) + \frac{353}{4}\zeta(3)\zeta(5) - 7\zeta(2)\zeta(3)^2 + \zeta(2)\sum^{\infty}_{k=1} \frac{h_k}{k^5} + 
\frac{1}{2}\sum^{\infty}_{k=1} \frac{h_k}{k^7} - \frac{221}{8}\sum^{\infty}_{k=1} \frac{H^{(2)}_k}{k^6}~.
\end{eqnarray}

\begin{eqnarray}
\sum^{\infty}_{k=1} \frac{H^{(2)}_k h^{(2)}_k}{k^4} &=& \frac{75}{2}\zeta(8) - 31\zeta(3)\zeta(5) + \frac{21}{2}\zeta(2)\zeta(3)^2 -
2\zeta(2)\sum^{\infty}_{k=1} \frac{h_k}{k^5} - 8\sum^{\infty}_{k=1} \frac{h^{(3)}_k}{k^5} \nonumber \\ &-& 6\sum^{\infty}_{k=1} \frac{h^{(4)}_k}{k^4}
- \frac{9}{2}\sum^{\infty}_{k=1} \frac{h^{(2)}_k}{k^6} \nonumber \\ &=&
\frac{3955}{32}\zeta(8) - \frac{633}{4}\zeta(3)\zeta(5) + \frac{21}{2}\zeta(2)\zeta(3)^2 - 2\zeta(2)\sum^{\infty}_{k=1} \frac{h_k}{k^5} - 
\frac{3}{2}\sum^{\infty}_{k=1} \frac{h_k}{k^7} + \frac{425}{8}\sum^{\infty}_{k=1} \frac{H^{(2)}_k}{k^6}
\end{eqnarray}
and
\begin{eqnarray}
\sum^{\infty}_{k=1} \frac{H^{(3)}_k h_k}{k^4} &=& -\frac{75}{4}\zeta(8) - \frac{35}{2}\zeta(3)\zeta(5) + \frac{7}{2}\zeta(2)\zeta(3)^2 +
\frac{21}{2}\sum^{\infty}_{k=1} \frac{h_k}{k^7} + 8\sum^{\infty}_{k=1} \frac{h^{(3)}_k}{k^5} \nonumber \\ &+& 4\sum^{\infty}_{k=1} \frac{h^{(4)}_k}{k^4}
+ 10\sum^{\infty}_{k=1} \frac{h^{(2)}_k}{k^6}
\nonumber \\ &=& \frac{225}{16}\zeta(8) - \frac{35}{2}\zeta(3)\zeta(5) + \frac{7}{2}\zeta(2)\zeta(3)^2 + \frac{1}{2}\sum^{\infty}_{k=1} \frac{h_k}{k^7}~. 
\end{eqnarray}

As a second example the calculational scheme for the last Euler sum will be presented in the following.

Using identity (10) from \cite{bra20} it follows:
\begin{eqnarray}
\sum^{\infty}_{k=1}\frac{H^{(2)}_k}{k^3} \left(\sum^{\infty}_{i=1,i\ne k} \frac{1}{i(i-k)} \right) = \sum^{\infty}_{k=1}\frac{H^{(2)}_k}{k^3} \left(\sum^{\infty}_{i=1}
\frac{1}{i(i+k)} \right) - \sum^{\infty}_{k=1} \frac{H^{(2)}_k h_k}{k^5} -2\sum^{\infty}_{k=1} \frac{H^{(2)}_k h^{(2)}_k}{k^4}~.
\end{eqnarray}
By rearrangement of the double summation and by partial fraction decomposition we get:
\begin{eqnarray}
&-&\sum^{\infty}_{i=1}\frac{h_i}{i^2} \left(\sum^{\infty}_{k=1,k\ne i} \frac{H^{(2)}_k}{k^3}\right) - \sum^{\infty}_{i=1}\frac{h_i}{i^3}\left(
\sum^{\infty}_{k=1,k\ne i} \frac{H^{(2)}_k}{k^2}\right) - \sum^{\infty}_{i=1}\frac{h_i}{i^3} \left(\sum^{\infty}_{k=1,k\ne i} \frac{H^{(2)}_k}{k(i-k)}\right) =
\sum^{\infty}_{i=1}\frac{h_i}{i^2} \left(\sum^{\infty}_{k=1} \frac{H^{(2)}_k}{k^3}\right) \nonumber \\ &-& \sum^{\infty}_{i=1}\frac{h_i}{i^3}\left(
\sum^{\infty}_{k=1} \frac{H^{(2)}_k}{k^2}\right) + \sum^{\infty}_{i=1}\frac{h_i}{i^3} \left(\sum^{\infty}_{k=1} \frac{H^{(2)}_k}{k(i+k)}\right) -
\sum^{\infty}_{k=1} \frac{H^{(2)}_k h_k}{k^5} -2\sum^{\infty}_{k=1} \frac{H^{(2)}_k h^{(2)}_k}{k^4}~.
\end{eqnarray}
Using now identity (86) for $n=1$:
\begin{eqnarray}
\sum^{\infty}_{k=1,k\ne i} \frac{H^{(2)}_k}{k(i-k)} = - \sum^{\infty}_{k=1} \frac{H^{(2)}_k}{k(i+k)} + \zeta(3)\frac{1}{i} - \zeta(2)\frac{1}{i^2} + 2\zeta(2)\frac{H_i}{i} +
\frac{H_i}{i^3} - \frac{H^{(2)}_i}{i^2} - 3\frac{H^{(3)}_i}{i}
\end{eqnarray}
one arrives after some standard algebraic manipulations at:
\begin{eqnarray}
\sum^{\infty}_{k=1} \frac{H^{(3)}_k h_k}{k^4} &=& \frac{2}{3}\sum^{\infty}_{i=1}\frac{h_i}{i^2} \left(\sum^{\infty}_{k=1} \frac{H^{(2)}_k}{k^3}\right) +
\frac{1}{3}\zeta(3)\sum^{\infty}_{k=1}\frac{h_k}{k^4} - \frac{1}{3}\zeta(2)\sum^{\infty}_{k=1}\frac{h_k}{k^5} + \frac{2}{3}\zeta(2)\sum^{\infty}_{k=1}\frac{H_k h_k}{k^4} +
\frac{1}{3}\zeta(2)\sum^{\infty}_{k=1}\frac{H_k h_k}{k^6} \nonumber \\ &-& \frac{4}{3} \sum^{\infty}_{k=1} \frac{H^{(2)}_k h_k}{k^5} -
\frac{2}{3} \sum^{\infty}_{k=1} \frac{H^{(2)}_k h^{(2)}_k}{k^4}~. 
\end{eqnarray}

As one can see, in the last expression Euler sums of lower order as well as two other sums of order 8 appear. This clearly shows that the calculational scheme is highly
recursive. The next step is then to calculate these two sums of order 8. We start by using identity (76):
\begin{eqnarray}
\sum^{\infty}_{k=1}\frac{h_k}{k^4} \left(\sum^{\infty}_{i=1,i\ne k} \frac{H_i}{i(k-i)} \right) = \sum^{\infty}_{k=1}\frac{h_k}{k^4} \left(\sum^{\infty}_{i=1}
\frac{H_i}{i(i+k)} \right) - \zeta(2)\sum^{\infty}_{k=1} \frac{h_k}{k^5} -2\sum^{\infty}_{k=1} \frac{H^{(2)}_k h_k}{k^5}~.
\end{eqnarray}
From this we get by rearrangement of the double summation, use of identity (88) from \cite{bra20} and by further algebraic manipulations:
\begin{eqnarray}
\sum^{\infty}_{k=1} \frac{H^{(2)}_k h_k}{k^5} = \sum^{\infty}_{k=1}\frac{H_k}{k^2} \left(\sum^{\infty}_{i=1} \frac{h_i}{i^4}\right) +
\sum^{\infty}_{k=1}\frac{H_k}{k^4} \left(\sum^{\infty}_{i=1} \frac{h_i}{i^2}\right) - 
\frac{1}{2}\zeta(2)\sum^{\infty}_{k=1}\frac{h_k}{k^5} - \frac{1}{2}\sum^{\infty}_{k=1}\frac{H_k h_k}{k^6} -
\sum^{\infty}_{k=1} \frac{H_k h^{(2)}_k}{k^5}~.
\end{eqnarray}
Here, only the last to terms are in principle unknown, but is has been shown before how to calculate these sums explicitly. Therefore, using Eqs.~(99) and (100) which
represent the results of the corresponding calculations one arrives at Eq.~(101)~. This way, it remains to calculate the last sum in Eq.~(107).

We start with the following expression:
\begin{eqnarray}
\sum^{\infty}_{k=1}\frac{h^{(2)}_k}{k^3} \left(\sum^{\infty}_{i=1,i\ne k} \frac{H_i}{i(k-i)} \right) = \sum^{\infty}_{k=1}\frac{h^{(2)}_k}{k^3} \left(\sum^{\infty}_{i=1}
\frac{H_i}{i(i+k)} \right) - \zeta(2)\sum^{\infty}_{k=1} \frac{h^{(2)}_k}{k^4} -2\sum^{\infty}_{k=1} \frac{H^{(2)}_k h^{(2)}_k}{k^4}~.
\end{eqnarray}
Analogously, by rearrangement of the double summation, use of identity (88) from \cite{bra20} and by further algebraic manipulations it follows:
\begin{eqnarray}
\sum^{\infty}_{k=1} \frac{H^{(2)}_k h^{(2)}_k}{k^4} = \sum^{\infty}_{k=1}\frac{H_k}{k^2} \left(\sum^{\infty}_{i=1} \frac{h^{(2)}_i}{i^3}\right) - \frac{1}{2}
\zeta(2)\sum^{\infty}_{k=1}\frac{h^{(2)}_k}{k^4} + \frac{3}{2}\zeta(2)\sum^{\infty}_{k=1}\frac{H_k h_k}{k^4} - \frac{3}{2}\sum^{\infty}_{k=1} \frac{H_k h^{(2)}_k}{k^5}
-\sum^{\infty}_{k=1}\frac{H_k h^{(3)}_k}{k^4}~.
\end{eqnarray}
As all Euler sums on the right side of Eq.~(111) are known, Eq.~(103) has been proven. We stop here with the calculations of further Euler sums of this type for the
order 8 case. The remaining Euler sums of order 8 can be found in the appendix. \\

Last but not least we present three examples of this type for the order 10 case. It follows: 
\begin{eqnarray}
\sum^{\infty}_{k=1} \frac{H_k h_k}{k^8} &=& - \frac{67}{2}\zeta(3)\zeta(7) - \frac{31}{4}\zeta(5)^2 + 19\zeta(2)\zeta(3)\zeta(5) + 
\frac{7}{2}\zeta(3)^2\zeta(4)  + \frac{9}{2}\sum^{\infty}_{k=1} \frac{h_k}{k^9} + \sum^{\infty}_{k=1} \frac{h^{(2)}_k}{k^8}
\nonumber \\ &=&
-\frac{103}{8}\zeta(10) - \frac{35}{2}\zeta(3)\zeta(7) + \frac{13}{4}\zeta(5)^2 + 19\zeta(2)\zeta(3)\zeta(5) + \frac{7}{2}\zeta(3)^2\zeta(4) + 
\frac{1}{2}\sum^{\infty}_{k=1} \frac{h_k}{k^9} \nonumber \\ &-& \frac{53}{4}\sum^{\infty}_{k=1} \frac{H^{(2)}_k}{k^8}~.
\end{eqnarray}
\begin{eqnarray}
\sum^{\infty}_{k=1} \frac{H_k h^{(2)}_k}{k^7} &=& \frac{381}{4}\zeta(3)\zeta(7) + \frac{31}{2}\zeta(5)^2 - 63\zeta(2)\zeta(3)\zeta(5) - 
\frac{7}{2}\zeta(3)^2\zeta(4) - \frac{3}{2}\zeta(2)\sum^{\infty}_{k=1} \frac{h_k}{k^7} + 4\sum^{\infty}_{k=1} \frac{h^{(2)}_k}{k^8} \nonumber \\ &+&
2\sum^{\infty}_{k=1} \frac{h^{(3)}_k}{k^7}
\nonumber \\ &=&
308\zeta(10) - \frac{339}{4}\zeta(3)\zeta(7) - \frac{593}{4}\zeta(5)^2 - 63\zeta(2)\zeta(3)\zeta(5) - \frac{7}{2}\zeta(3)^2\zeta(4) - 
\frac{3}{2}\zeta(2)\sum^{\infty}_{k=1} \frac{h_k}{k^7} \nonumber \\ &-& 2\sum^{\infty}_{k=1}\frac{h_k}{k^9} + \frac{187}{2}\sum^{\infty}_{k=1} \frac{H^{(2)}_k}{k^8}~.
\end{eqnarray}
\begin{eqnarray}
\sum^{\infty}_{k=1} \frac{H^{(2)}_k h_k}{k^7} &=& \frac{303}{4}\zeta(3)\zeta(7) + 31\zeta(5)^2 - 45\zeta(2)\zeta(3)\zeta(5) - 
7\zeta(3)^2\zeta(4) + \zeta(2)\sum^{\infty}_{k=1} \frac{h_k}{k^7} - \frac{27}{2}\sum^{\infty}_{k=1}\frac{h_k}{k^9} \nonumber \\ &-& 
7\sum^{\infty}_{k=1} \frac{h^{(2)}_k}{k^8} - 2\sum^{\infty}_{k=1} \frac{h^{(3)}_k}{k^7}
\nonumber \\ &=&
-\frac{2155}{8}\zeta(10) + \frac{831}{4}\zeta(3)\zeta(7) + \frac{647}{4}\zeta(5)^2 - 45\zeta(2)\zeta(3)\zeta(5) - 7\zeta(3)^2\zeta(4) + 
\zeta(2)\sum^{\infty}_{k=1}\frac{h_k}{k^7} \nonumber \\ &+& \frac{1}{2}\sum^{\infty}_{k=1}\frac{h_k}{k^9} - \frac{215}{4}\sum^{\infty}_{k=1} \frac{H^{(2)}_k}{k^8}~.
\end{eqnarray}
The remaining Euler sums of order 10 belonging to this family can be calculated recursively starting with the sums in exactly this order: 
\begin{eqnarray}
\sum^{\infty}_{k=1} \frac{H_k h_k^{(3)}}{k^6},~~~~\sum^{\infty}_{k=1} \frac{H^{(2)}_k h_k^{(2)}}{k^6},~~~~\sum^{\infty}_{k=1} \frac{H^{(3)}_k h_k}{k^6}, \\
\sum^{\infty}_{k=1} \frac{H_k h_k^{(4)}}{k^5},~~~~\sum^{\infty}_{k=1} \frac{H^{(2)}_k h_k^{(3)}}{k^5},~~~~\sum^{\infty}_{k=1} \frac{H^{(3)}_k h_k^{(2)}}{k^5},
~~~~\sum^{\infty}_{k=1} \frac{H^{(4)}_k h_k^{(1)}}{k^5}~...
\end{eqnarray}
and so on for denominators $n\le4$.

\section{Second family}

The following type of nonlinear Euler sums can be expressed in terms of linear Euler sums and zeta values:
\begin{eqnarray}
\sum^{\infty}_{k=1} \frac{H^{(a)}_k H^{(b)}_k}{(2k-1)^c}
\end{eqnarray}

with $a,b,c \in \mathbb{N}$, $a+b+c = $even.

The following identities hold:
\subsection{Lemma 2a}
\begin{eqnarray}
\sum^{\infty}_{k=1} \frac{H_k}{(2k+1)(i+k)} = \left(\zeta(2) - 2 \left(ln(2)\right)^2 \right)\frac{1}{2i-1} -\frac{H_i}{i} +\frac{1}{2} \frac{H^{(2)}_i}{2i-1} +
\frac{1}{2} \frac{H^2_i}{2i-1}~.
\end{eqnarray}
\subsection{proof}
We define:
\begin{eqnarray}
F(i) = \sum^{\infty}_{k=1} \frac{H_k}{(2k+1)(i+k)}~.
\end{eqnarray}
From this we get:
\begin{eqnarray}
F(i+1) &=& \sum^{\infty}_{k=1} \frac{H_k}{(2k+1)(i+k+1)} = \sum^{\infty}_{k=1} \frac{H_k}{(2k-1)(i+k)} - \sum^{\infty}_{k=1} \frac{1}{k(2k-1)(i+k)}
\nonumber \\ &=& \frac{1}{2i+1}\sum^{\infty}_{k=1} \frac{2H_k}{2k-1} - \frac{1}{2i+1}\sum^{\infty}_{k=1} \frac{H_k}{i+k} - \sum^{\infty}_{k=1} \frac{1}{k(2k-1)(i+k)}
\nonumber \\ &=& \frac{2}{2i+1} + \frac{1}{2i+1}\sum^{\infty}_{k=1}\frac{2H_{k+1}}{2k+1} - \frac{1}{2i+1}\sum^{\infty}_{k=1} \frac{H_k}{i+k} -
\sum^{\infty}_{k=1} \frac{1}{k(2k-1)(i+k)} \nonumber \\ &=& \frac{2}{2i+1} + \frac{2}{2i+1} \sum^{\infty}_{k=1}\frac{1}{(k+1)(2k+1)} +
\frac{1}{2i+1}\sum^{\infty}_{k=1} \left(\frac{2H_k}{2k+1}-\frac{H_k}{i+k} \right) - \sum^{\infty}_{k=1} \frac{1}{k(2k-1)(i+k)}
\nonumber \\ &=& \frac{2i-1}{2i+1}F(i) + \frac{2}{2i+1} + \frac{4ln(2)-2}{2i+1} - \frac{1}{2i} \sum^{\infty}_{k=1}\frac{1}{k(2k-1)} +
\frac{1}{2i} \sum^{\infty}_{k=1}\frac{1}{(2k-1)(i+k)} \nonumber \\ &=& \frac{2i-1}{2i+1}F(i) + \frac{4ln(2)}{2i+1} -\frac{2ln(2)}{2i} + \frac{H_i}{2i(2i+1)} +
\frac{2ln(2)}{2i(2i-1)}~.
\end{eqnarray}
From this we get:
\begin{eqnarray}
F(i+1) - \frac{2i-1}{2i+1}F(i) = \frac{2ln(2)}{2i+1} + \frac{H_i}{2i(2i+1)}~.
\end{eqnarray}
The solution of this first order difference equation is given by Eq.~(118). Thus the lemma is proved.

\subsection{Lemma 2b}
\begin{eqnarray}
\sum^{\infty}_{k=1} \frac{H_k}{(i+k)^n} = \frac{H_k}{(k+1)^n} - \sum^{i-1}_{k=1} \left( \sum^{\infty}_{m=1}\frac{1}{m(k+m)^n} \right)~.
\end{eqnarray}

\subsection{proof}
We define:
\begin{eqnarray}
F(i) = \sum^{\infty}_{k=1} \frac{H_k}{(i+k)^n}~.
\end{eqnarray}
From this we get:
\begin{eqnarray}
F(i+1) = \sum^{\infty}_{k=1} \frac{H_k}{(i+k+1)^n} = \sum^{\infty}_{k=1} \frac{H_{k-1}}{(i+k)^n} = F(i) - \sum^{\infty}_{k=1}\frac{1}{k(i+k)^n}~. 
\end{eqnarray}
The solution of the corresponding difference equation is simply given by Eq.~(122). Thus the lemma is proved. To give an example, we explicitly calculate the identity
for the case $n=2$ which is needed, for example, if $a=1$ and $b=2$. Using Eq.~(122) it follows by partial fraction decomposition:
\begin{eqnarray}
\sum^{\infty}_{k=1} \frac{H_k}{(i+k)^2} &=& \frac{H_k}{(k+1)^2} - \sum^{i-1}_{k=1} \frac{1}{k} \left( \sum^{\infty}_{m=1}\frac{1}{m(k+m)} \right) +
\sum^{i-1}_{k=1} \frac{1}{k} \left( \sum^{\infty}_{m=1}\frac{1}{(k+m)^2} \right) \nonumber \\ &=&
\frac{H_k}{(k+1)^2} - \sum^{i-1}_{k=1}\frac{H_k}{k^2} + \sum^{i-1}_{k=1}\frac{1}{k}\left(\zeta(2) - H^{(2)}_k \right)  \nonumber \\ &=&
\frac{H_k}{(k+1)^2} - \sum^{i-1}_{k=1}\frac{H_k}{k^2} + \zeta(2)H_{i-1} - \sum^{i-1}_{k=1}\frac{H^{(2)}_k}{k}~.
\end{eqnarray}
After some standard algebraic manipulations one gets the final result:
\begin{eqnarray}
\sum^{\infty}_{k=1} \frac{H_k}{(i+k)^2} = \zeta(3) - \zeta(2)\frac{1}{i} + \zeta(2)H_{i} - H_{i}H^{(2)}_{i} - H^{(3)}_{i} + \frac{H_i}{i^2} + \frac{H^{(2)}_i}{i}~. 
\end{eqnarray}

\subsection{Lemma 2c}
\begin{eqnarray}
\sum^{\infty}_{k=1} \frac{H^{(a)}_k H^{(b)}_k}{(2k-1)^c} &=& \frac{2^c-1}{2^c} \left( \zeta(a)\zeta(b)\zeta(c) + \zeta(a+b)\zeta(c) \right) +
\sum^{\infty}_{k=1} \frac{H^{(b)}_k}{k^a(2k-1)^c} + \sum^{\infty}_{k=1} \frac{H^{(a)}_k}{k^b(2k-1)^c} \nonumber \\ &-&
\sum^{\infty}_{k=1} \frac{H^{(a+b)}_k}{(2k-1)^c} - \sum^{\infty}_{k=1} \frac{H^{(a)}_k h^{(c)}_k}{k^b} - \sum^{\infty}_{k=1} \frac{H^{(b)}_k h^{(c)}_k}{k^a}~.
\end{eqnarray}

\subsection{proof}
We start with the identity \cite{ade16}:
\begin{eqnarray}
\sum^{k}_{i=1} \frac{H^{(a)}_i}{i^b} = H^{(a)}_k H^{(b)}_b + H^{(a+b)}_k - \sum^{k}_{i=1} \frac{H^{(b)}_i}{i^a}~.
\end{eqnarray}
From this it follows:
\begin{eqnarray}
\sum^{\infty}_{k=1} \frac{1}{(2k-1)^c}\left( \sum^{k}_{i=1} \frac{H^{(a)}_i}{i^b} \right) = \sum^{\infty}_{k=1} \frac{H^{(a)}_k H^{(b)}_b}{(2k-1)^c} + 
\sum^{\infty}_{k=1}\frac{H^{(a+b)}_k }{(2k-1)^c} - \sum^{\infty}_{k=1} \frac{1}{(2k-1)^c} \left( \sum^{k}_{i=1} \frac{H^{(b)}_i}{i^a} \right) ~.
\end{eqnarray}
Rearranging the corresponding double summations by use the following identity:
\begin{eqnarray}
\sum^{\infty}_{k=1} \frac{H^{(a)}_k}{k^b} = \sum^{\infty}_{k=1} \frac{1}{k^b} \left( \sum^{k}_{i=1} \frac{1}{i^a} \right) = \zeta(a)\zeta(b) + \zeta(a+b) -
\sum^{\infty}_{k=1} \frac{H^{(b)}_k}{k^a}
\end{eqnarray}
and performing some standard algebraic manipulations Eq.~(127) follows. Thus lemma 2c is proved.

This way, Lemma 2a has to be applied in the case where $a=b=1$. Lemma 2b has to be applied in the case where $a=1$ and $b>1$ or vice versa. Lemma 3c is devoted to the
case where both $a>1$ and $b>1$.
\subsection{examples}

The order 4 case:
\begin{eqnarray}
\sum^{\infty}_{k=1} \frac{H_k^2}{(2k-1)^2} &=& \frac{45}{16}\zeta(4) -7ln(2)\zeta(3) + 7\zeta(3) + 3\left(ln(2)\right)^2 \zeta(2) - 6ln(2)\zeta(2) + 2\zeta(2) -
8ln(2) \nonumber \\ &+& 4\left(ln(2)\right)^2 + \sum^{\infty}_{k=1} \frac{h_k}{k^3}~.
\end{eqnarray}
The order 5 case with $a=b=1$:
\begin{eqnarray}
\sum^{\infty}_{k=1} \frac{H_k^2}{(2k-1)^3} &=& \frac{31}{8}\zeta(5) - \frac{7}{8} \zeta(2)\zeta(3) - \frac{45}{8}ln(2)\zeta(4) + \frac{45}{8}\zeta(4) + 
\frac{7}{2}\left(ln(2)\right)^2\zeta(3) - 7ln(2) \zeta(3) - \frac{7}{2}\zeta(3) + 6ln(2)\zeta(2) \nonumber \\ &-& 5\zeta(2) + 12 ln(2) -4\left(ln(2)\right)^2~.
\end{eqnarray}
The order 6 case with $a=b=1$:
\begin{eqnarray}
\sum^{\infty}_{k=1} \frac{H^2_k}{(2k-1)^4} &=& \frac{585}{128}\zeta(6) - \frac{31}{2}ln(2)\zeta(5) + \frac{31}{2}\zeta(5) + \frac{15}{4}\left(ln(2)\right)^2
\zeta(4) - \frac{15}{2}ln(2)\zeta(4) - \frac{15}{8}\zeta(4) + \frac{21}{4}ln(2)\zeta(2)\zeta(3)  \nonumber \\  &-& \frac{35}{16}\zeta(3)^2 -\frac{21}{4}\zeta(2)\zeta(3) +
7ln(2)\zeta(3) - 6ln(2)\zeta(2) + 8\zeta(2) - 16ln(2) + 4\left(ln(2)\right)^2 + \frac{1}{2}\sum^{\infty}_{k=1} \frac{h_k}{k^5}~.
\end{eqnarray}
The order 6 case:
\begin{eqnarray}
\sum^{\infty}_{k=1} \frac{H^{(2)}_k H^{(2)}_k}{(2k-1)^2} &=& \frac{195}{32}\zeta(6) - 5\zeta(4) -7\zeta(3)^2 + 16\zeta(3) - 16ln(2)\zeta(2) +
24\zeta(2) - 64ln(2) \nonumber \\  &+& 8 \sum^{\infty}_{k=1} \frac{h_k}{k^3} + 2\zeta(2)\sum^{\infty}_{k=1} \frac{h_k}{k^3}~.
\end{eqnarray}
As a second example we present:
\begin{eqnarray}
\sum^{\infty}_{k=1} \frac{H_k H^{(2)}_k}{(2k-1)^3} &=& -\frac{1035}{128}\zeta(6) + \frac{93}{4}ln(2)\zeta(5) - \frac{93}{4}\zeta(5) +\frac{15}{2}\zeta(4)
+ \frac{203}{32}\zeta(3)^2 - \frac{49}{4}ln(2)\zeta(2)\zeta(3) + \frac{49}{4}\zeta(2)\zeta(3) \nonumber \\  &-& 7ln(2)\zeta(3) -  \frac{21}{2}\zeta(3)
+ 14ln(2)\zeta(2) -18\zeta(2) + 48ln(2) - 12\left(ln(2)\right)^2 - 2\sum^{\infty}_{k=1} \frac{h_k}{k^3}
\nonumber \\  &-& \frac{3}{4}\sum^{\infty}_{k=1} \frac{h_k}{k^5}~.
\end{eqnarray}

The order 7 case with $a=b=1$:
\begin{eqnarray}
\sum^{\infty}_{k=1} \frac{H_k^2}{(2k-1)^5} &=& \frac{635}{64}\zeta(7) - \frac{315}{32}ln(2)\zeta(6) + \frac{315}{32}\zeta(6) + \frac{31}{8}\left(ln(2)\right)^2\zeta(5)
- \frac{31}{32}\zeta(2)\zeta(5) - \frac{31}{4}ln(2)\zeta(5)  \nonumber \\ &-& \frac{93}{8}\zeta(5) - \frac{315}{64}\zeta(3)\zeta(4) + \frac{15}{2}ln(2)\zeta(4) -
\frac{15}{8}\zeta(4) + \frac{49}{16}ln(2)\zeta(3)^2 - \frac{49}{16}\zeta(3)^2 + \frac{21}{4}\zeta(2)\zeta(3) \nonumber \\ &-& 7ln(2)\zeta(3) + \frac{7}{2}\zeta(3) +
6ln(2)\zeta(2) - 11\zeta(2) + 20ln(2) - 4\left(ln(2)\right)^2~.
\end{eqnarray}
The order 8 case:
\begin{eqnarray}
\sum^{\infty}_{k=1} \frac{H^{(2)}_k H^{(3)}_k}{(2k-1)^3} &=& \frac{3375}{16}\zeta(8) + \frac{135}{4}\zeta(6) - \frac{279}{24}\zeta(3)\zeta(5) - 217\zeta(5) +
\frac{39}{2}\zeta(4) - \frac{2331}{48}\zeta(2)\zeta(3)^2 - \frac{77}{4}\zeta(3)^2 + 112\zeta(2)\zeta(3) \nonumber \\ &+& 12ln(2)\zeta(3) - 32\zeta(3) +
48ln(2)\zeta(2) - 200\zeta(2) + 480ln(2) -36\sum^{\infty}_{k=1} \frac{h_k}{k^3} - \frac{3}{2}\zeta(2)\sum^{\infty}_{k=1} \frac{h_k}{k^5} - 6\sum^{\infty}_{k=1} \frac{h_k}{k^7}
\nonumber \\ &+& \sum^{\infty}_{k=1} \frac{h^{(3)}_k}{k^5} -12 \sum^{\infty}_{k=1} \frac{h^{(5)}_k}{k^3} -40\sum^{\infty}_{k=1} \frac{h^{(6)}_k}{k^2}~.
\end{eqnarray}
One should notice here that, although lemma 2b and lemma 2c are applied to the case $a+b+c =$even only, lemma 2a can be used to calculate explicitly the corresponding Euler
sums for $a=b=1$ and c of arbitrary order. Even more remarkable is the fact that for $c =$ odd the corresponding Euler sums can be expressed in terms of zeta values only.

\section{Third family}

The following type of nonlinear Euler sums can be expressed in terms of linear Euler sums and zeta values:
\begin{eqnarray}
\sum^{\infty}_{k=1} \frac{H^{(a)}_k H^{(b)}_k}{k(2k-1)}
\end{eqnarray}

with $a,b,c \in \mathbb{N}$, $a+b+1 = $even.

The following identity holds:

\subsection{Lemma 3a}
\begin{eqnarray}
\sum^{\infty}_{k=1} \frac{H^{(n)}_k}{k(i+k)} &=& \zeta(n+1)\frac{1}{i} + \sum^{n-1}_{m=1}(-)^{m+1}\zeta(n+1-m)\frac{H^{(m)}_{i-1}}{i} + \frac{(-)^{n+1}}{i}
\sum^{i-1}_{k=1}\frac{H_k}{k^n}~.
\end{eqnarray}

\subsection{proof}
We define:
\begin{eqnarray}
F(i) = \sum^{\infty}_{k=1} \frac{H^{(n)}_k}{k(i+k)}~.
\end{eqnarray}
From this we get:
\begin{eqnarray}
F(i+1) = \sum^{\infty}_{k=1} \frac{H^{(n)}_k}{k(i+k+1)} &=& \frac{1}{i+1}\sum^{\infty}_{k=1}\frac{H^{(n)}_k}{k} - \frac{1}{i+1}\sum^{\infty}_{k=1}\frac{H^{(n)}_k}{i+k+1}
\nonumber \\ &=&
\frac{1}{i+1}\sum^{\infty}_{k=1}\frac{H^{(n)}_k}{k} - \frac{1}{i+1}\sum^{\infty}_{k=1}\frac{H^{(n)}_{k-1}}{i+k}
\nonumber \\ &=&
\frac{1}{i+1}\sum^{\infty}_{k=1}\frac{H^{(n)}_k}{k} - \frac{1}{i+1}\sum^{\infty}_{k=1}\frac{H^{(n)}_k}{i+k} + \frac{1}{i+1}\sum^{\infty}_{k=1}\frac{1}{k^n(i+k)}
\nonumber \\ &=&
\frac{i}{i+1}F(i) + \frac{1}{i+1}\sum^{\infty}_{k=1}\frac{1}{k^n(i+k)}~.
\end{eqnarray}
This way we have:
\begin{eqnarray}
F(i+1) - \frac{i}{i+1}F(i) = \frac{1}{i+1}\sum^{\infty}_{k=1}\frac{1}{k^n(i+k)}~. 
\end{eqnarray}
The solution of the corresponding difference equation is simply given by Eq.~(139). Thus the lemma is proved. To give an example, we explicitly calculate the identity
for $n=1$. It follows:
\begin{eqnarray}
F(i+1) - \frac{i}{i+1}F(i) = \frac{1}{i+1}\sum^{\infty}_{k=1}\frac{1}{k(i+k)}~,
\end{eqnarray}
with the solution:
\begin{eqnarray}
\sum^{\infty}_{k=1} \frac{H_k}{k(i+k)} = \zeta(2)\frac{1}{i} + \frac{1}{i}\sum^{i-1}_{k=1}\frac{H_k}{k}~.
\end{eqnarray}
Dividing now  Eq.~(139) on both sides by the factor $(2i+1)$ and summing up over i from $1$ to $\infty$ it follows:  
\begin{eqnarray}
\sum^{\infty}_{i=1} \frac{1}{2i+1} \left( \sum^{\infty}_{k=1} \frac{H^{(n)}_k}{k(i+k)} \right) &=& \zeta(n+1) \left( 2-2ln(2) \right)
+ \sum^{n-1}_{m=1}(-)^{m+1}\zeta(n+1-m) \sum^{\infty}_{i=1} \frac{H^{(m)}_{i-1}}{i(2i+1)} \nonumber \\ &+&
\sum^{\infty}_{i=1} \frac{(-)^{n+1}}{i(2i+1)} \left( \sum^{i-1}_{k=1}\frac{H_k}{k^n} \right)~.
\end{eqnarray}
As $a+b+1 = $even it follows in the case $b=1$ that $n=2a$ is demanded. This way we get: 
\begin{eqnarray}
\sum^{\infty}_{i=1} \frac{1}{2i+1} \left( \sum^{\infty}_{k=1} \frac{H^{(2a)}_k}{k(i+k)} \right) &=& \zeta(2a+1) \left( 2-2ln(2) \right)
+ \sum^{2a-1}_{m=1}(-)^{m+1}\zeta(2a+1-m) \sum^{\infty}_{i=1} \frac{H^{(m)}_{i-1}}{i(2i+1)} \nonumber \\ &-&
\sum^{\infty}_{i=1} \frac{1}{i(2i+1)} \left( \sum^{i-1}_{k=1}\frac{H_k}{k^{2a}} \right) ~.
\end{eqnarray}
It follows then:
\begin{eqnarray}
\sum^{\infty}_{k=1} \frac{H^{(2a)}_k}{k} \left( \sum^{\infty}_{i=1} \frac{1}{(2i+1)(i+k)} \right) &=&
\sum^{\infty}_{k=1} \frac{H^{(2a)}_k}{k} \left( \frac{H_k}{2k-1} - \sum^{\infty}_{k=1} \frac{2-2ln(2)}{2k-1} \right) =
\sum^{\infty}_{k=1} \frac{H_k H^{(2a)}_k}{k(2k-1)} \nonumber \\ &-& (2-2ln(2)) \sum^{\infty}_{k=1} \frac{H^{(2a)}_k}{k(2k-1)}~.
\end{eqnarray}
After some standard algebraic manipulations we end up with the following expression:
\begin{eqnarray}
\sum^{\infty}_{k=1} \frac{H_k H^{(2a)}_k}{k(2k-1)} &=& (2-2ln(2)) \left( \zeta(2a+1) + \sum^{\infty}_{k=1} \frac{H^{(2a)}_k}{k(2k-1)} \right) +
\sum^{2a-1}_{m=1}(-)^{m+1}\zeta(2a+1-m) \sum^{\infty}_{i=1} \frac{H^{(m)}_{i-1}}{i(2i+1)} \nonumber \\ &-&
\sum^{\infty}_{i=1} \frac{1}{i(2i+1)} \left( \sum^{i-1}_{k=1}\frac{H_k}{k^{2a}} \right)~.
\end{eqnarray}
As an example we calculate explicitly the corresponding Euler sum for $a=1$. It follows:
\begin{eqnarray}
\sum^{\infty}_{k=1} \frac{H_k H^{(2)}_k}{k(2k-1)} &=& (2-2ln(2)) \left( \zeta(3) + \sum^{\infty}_{k=1} \frac{H^{(2)}_k}{k(2k-1)} \right) +
\zeta(2)\sum^{\infty}_{k=1} \frac{H_k}{k(2k+1)} - \zeta(2)\sum^{\infty}_{k=1} \frac{1}{k^2(2k+1)} \nonumber \\ &-&
(2-2ln(2)) \sum^{\infty}_{k=1} \frac{H_k}{k^2} + 2\sum^{\infty}_{k=1}\frac{H_k (h_{k+1}-1)}{k^2}~.
\end{eqnarray}
As all terms on the right side are known, we get after some standard algebraic manipulations at Eq.~(154).

In order to calculate explicitly members of this family for $b>1$ an additional identity is needed. It follows:
\subsection{Lemma 3b} 
\begin{eqnarray}
\sum^{\infty}_{k=1} \frac{H^{(a)}_k H^{(b)}_k}{k(2k-1)} &=& 2ln(2) \left( \zeta(a)\zeta(b) + \zeta(a+b) \right) + \sum^{\infty}_{k=1} \frac{H^{(a)}_k}{k^{b+1}(2k-1)}
+ \sum^{\infty}_{k=1} \frac{H^{(b)}_k}{k^{a+1}(2k-1)} + \sum^{\infty}_{k=1} \frac{H_k H^{(a)}_k}{k^b} \nonumber \\ &+& \sum^{\infty}_{k=1} \frac{H_k H^{(b)}_k}{k^a} - 
2\sum^{\infty}_{k=1} \frac{H^{(a)}_k h_k}{k^b} - 2 \sum^{\infty}_{k=1} \frac{H^{(b)}_k h_k}{k^a} - \sum^{\infty}_{k=1} \frac{H^{(a+b)}_k}{k(2k-1)}~.
\end{eqnarray}

\subsection{proof}
Again, we start with the identity \cite{ade16}:
\begin{eqnarray}
\sum^{k}_{i=1} \frac{H^{(a)}_i}{i^b} = H^{(a)}_k H^{(b)}_b + H^{(a+b)}_k - \sum^{k}_{i=1} \frac{H^{(b)}_i}{i^a}~.
\end{eqnarray}
From this it follows:
\begin{eqnarray}
\sum^{\infty}_{k=1} \frac{1}{k(2k-1)}\left( \sum^{k}_{i=1} \frac{H^{(a)}_i}{i^b} \right) = \sum^{\infty}_{k=1} \frac{H^{(a)}_k H^{(b)}_b}{k(2k-1)} + 
\sum^{\infty}_{k=1}\frac{H^{(a+b)}_k }{k(2k-1)} - \sum^{\infty}_{k=1} \frac{1}{k(2k-1)} \left( \sum^{k}_{i=1} \frac{H^{(b)}_i}{i^a} \right) ~.
\end{eqnarray}
Rearranging the corresponding double summations by use the following identity:
\begin{eqnarray}
\sum^{\infty}_{k=1} \frac{H^{(a)}_k}{k^b} = \sum^{\infty}_{k=1} \frac{1}{k^b} \left( \sum^{k}_{i=1} \frac{1}{i^a} \right) = \zeta(a)\zeta(b) + \zeta(a+b) -
\sum^{\infty}_{k=1} \frac{H^{(b)}_k}{k^a}
\end{eqnarray}
and performing some standard algebraic manipulations Eq.~(150) follows. Thus lemma 3b is proved.

\subsection{examples}
The order 4 case:
\begin{eqnarray}
\sum^{\infty}_{k=1} \frac{H_k H^{(2)}_k}{k(2k-1)} &=& -\frac{31}{8}\zeta(4) + 7ln(2)\zeta(3) - 7\zeta(3) - 2\left(ln(2)\right)^2 \zeta(2) + 4ln(2)\zeta(2) +  
16ln(2) - 8\left(ln(2)\right)^2 \nonumber \\ &-& \sum^{\infty}_{k=1} \frac{h_k}{k^3}~.
\end{eqnarray}

The order 6 case:
\begin{eqnarray}
\sum^{\infty}_{k=1} \frac{H_k H^{(4)}_k}{k(2k-1)} &=& -\frac{535}{48}\zeta(6) + 31ln(2)\zeta(5) - 31\zeta(5) - 2\left(ln(2)\right)^2 \zeta(4) + 4ln(2)\zeta(4) -  
5\zeta(4) + \frac{7}{2}\zeta(3)^2 - 14ln(2)\zeta(2)\zeta(3) \nonumber \\ &+& 14\zeta(2)\zeta(3) - 16\zeta(3) + 64ln(2) -32\left(ln(2)\right)^2 + 2\zeta(2)\sum^{\infty}_{k=1}
\frac{h_k}{k^3} - \sum^{\infty}_{k=1} \frac{h_k}{k^5}~.
\end{eqnarray}
One should notice here that lemma 3b has to be applied to the case $a,b >1$, ($a+b+1$) $=$ even, whereas lemma 3a must be used to calculate explicitly the corresponding Euler
sums for $a=1, b>1$ or $a>1, b=1$, ($a+b+1$) $=$ even.

The order 8 case:
\begin{eqnarray}
\sum^{\infty}_{k=1} \frac{H^{(3)}_k H^{(4)}_k}{k(2k-1)} &=& -\frac{22207}{144}\zeta(8) - \frac{31}{3}\zeta(6) + \frac{427}{2}\zeta(3)\zeta(5) - 144\zeta(5) +
2ln(2)\zeta(3)\zeta(4) - 2\zeta(3)\zeta(4)  \nonumber \\ &+& 16ln(2)\zeta(4) + 4\zeta(4) -24\zeta(2)\zeta(3)^2 + 2\zeta(3)^2 + 64\zeta(2)\zeta(3) + 32ln(2)\zeta(3) -
32\zeta(3) \nonumber \\ &-& 64\zeta(2) + 256ln(2) - 32\sum^{\infty}_{k=1} \frac{h_k}{k^3} -2\zeta(4)\sum^{\infty}_{k=1} \frac{h_k}{k^3} - \frac{211}{4}
\sum^{\infty}_{k=1} \frac{H^{(2)}_k}{k^6}~.
\end{eqnarray}

\begin{eqnarray}
\sum^{\infty}_{k=1} \frac{H_k H^{(6)}_k}{k(2k-1)} &=& -\frac{275}{24}\zeta(8) + 127ln(2)\zeta(7) - 127\zeta(7) -2\left(ln(2)\right)^2\zeta(6) + 4ln(2)\zeta(6) -7\zeta(6) -
6\zeta(3)\zeta(5) \nonumber \\ &-& 62ln(2)\zeta(2)\zeta(5) + 62\zeta(2)\zeta(5) - 24\zeta(5) - 14ln(2)\zeta(3)\zeta(4) + 14\zeta(3)\zeta(4) - 20\zeta(4) \nonumber \\ &+&
\frac{1}{2}\zeta(2)\zeta(3)^2 + 2\zeta(3)^2 + 8\zeta(2)\zeta(3) - 64\zeta(3) + 256ln(2) -128\left(ln(2)\right)^2 + 2\zeta(4)\sum^{\infty}_{k=1} \frac{h_k}{k^3} 
\nonumber \\ &+& 2\zeta(2)\sum^{\infty}_{k=1} \frac{h_k}{k^5} - \sum^{\infty}_{k=1} \frac{h_k}{k^7} + \frac{19}{2}\sum^{\infty}_{k=1} \frac{H^{(2)}_k}{k^6}~.
\end{eqnarray}

\section{Fourth family}

The following type of nonlinear Euler sums can be expressed in terms of linear Euler sums and zeta values:
\begin{eqnarray}
\sum^{\infty}_{k=1} \frac{h^{(a)}_k h^{(b)}_k}{(2k-1)^{c+1}}
\end{eqnarray}

with $a,b,c \in \mathbb{N}$, $a+b+c+1 = $even.

The following identity holds:
\subsection{Lemma 4a}
\begin{eqnarray}
\sum^{\infty}_{k=1} \frac{h^{(a)}_k}{(2k-1)(2i+2k-1)} &=& \sum^{a}_{n=1}\frac{(-)^{n+1}}{i^n} \frac{2^{a+2-n}-1}{2^{a+2}}\zeta(a+2-n) \nonumber \\ &+&
\sum^{a-1}_{n=1} \left( (-)^{n+1}\zeta(a+1-n)\frac{2^{a+1-n}-1}{2^{a+2}}\frac{H^{(n)}_i}{i} \right) + \frac{(-)^{a+1}}{2^{a+1}} \frac{1}{i} \sum^{i-1}_{k=1}\frac{h_k}{k^a}~. 
\end{eqnarray}

\subsection{proof}
We start with:
\begin{eqnarray}
F(i) = \sum^{\infty}_{k=1} \frac{h^{(a)}_k }{(2k-1)(2i+2k-1)}~.
\end{eqnarray}
From this we get
\begin{eqnarray}
F(i+1) &=& \sum^{\infty}_{k=1} \frac{h^{(a)}_k }{(2k-1)(2i+2k+1)} =  \frac{1}{2i+2}\sum^{\infty}_{k=1} \frac{h^{(a)}_k }{2k-1} -
\frac{1}{2i+2}\sum^{\infty}_{k=1} \frac{h^{(a)}_k }{2i+2k+1} \nonumber \\ &=&
\frac{1}{2i+2}\sum^{\infty}_{k=1} \frac{h^{(a)}_k }{2k-1} - \frac{1}{2i+2}\sum^{\infty}_{k=1} \frac{h^{(a)}_{k-1}}{2i+2k-1} \nonumber \\ &=&
\frac{1}{2i+2}\sum^{\infty}_{k=1} \frac{h^{(a)}_k }{2k-1} - \frac{1}{2i+2}\sum^{\infty}_{k=1} \frac{h^{(a)}_k}{2i+2k-1} +
\frac{1}{2i+2}\sum^{\infty}_{k=1} \frac{1}{(2k-1)^a (2i+2k-1)} \nonumber \\ &=& \frac{i}{i+1}F(i) +  \frac{1}{2i+2}\sum^{\infty}_{k=1}\frac{1}{(2k-1)^a (2i+2k-1)}~.
\end{eqnarray}
This way we result in the following inhomogeneous difference equation of first order:
\begin{eqnarray}
F(i+1) - \frac{i}{i+1}F(i) = \frac{1}{2i+2}\sum^{\infty}_{k=1}\frac{1}{(2k-1)^a (2i+2k-1)}~.
\end{eqnarray}
The solution is given by Eq.~(159). 

As an example we present the explicit result for a = 1. In this case it follows:
\begin{eqnarray}
F(i+1) - \frac{i}{i+1}F(i) = \frac{h_i}{2i(2i+2)}~,
\end{eqnarray}
with the solution:
\begin{eqnarray}
\sum^{\infty}_{k=1} \frac{h_k}{(2k-1)(2i+2k-1)} &=& 
F(1)\frac{1}{i} + \frac{1}{4i}\sum^{i-1}_{k=1} \frac{h_k}{k}
\nonumber \\ &=& 
\frac{1}{i} \sum^{\infty}_{k=1}\frac{h_k}{(2k-1)(2k+1)} + \frac{1}{4i}\sum^{i-1}_{k=1} \frac{h_k}{k} 
\nonumber \\ &=& 
\frac{3}{8}\zeta(2)\frac{1}{i} + \frac{1}{4i}\sum^{i-1}_{k=1} \frac{h_k}{k}~. 
\end{eqnarray}

Dividing now both sides of Eq.~(159) by $i^c$ and summing up over i from 1 to $\infty$ we get:
\begin{eqnarray}
\sum^{\infty}_{k=1} \frac{h^{(a)}_k h_k}{(2k-1)^{c+1}} &=& ln(2) \sum^{\infty}_{k=1} \frac{h^{(a)}_k}{(2k-1)^{c+1}} + \frac{1}{2^{c+1}}
\sum^{c-1}_{n=1}(-)^{n+1}2^n \zeta(c-n+1) \sum^{\infty}_{k=1} \frac{h^{(a)}_k }{(2k-1)^{c-n+1}} \nonumber \\ &-&
\frac{(-)^c}{2^{a+c+2}} \Big( \sum^{a}_{n=1} (-)^{n+1} \left(2^{a+2-n} -1 \right) \zeta(a+2-n)\zeta(c+n) \nonumber \\ &+& 
\sum^{a-1}_{n=1} (-)^{n+1} \left( 2^{a+1-n}- 1\right) \zeta(a+1-n)\sum^{\infty}_{k=1}
\frac{H^{(n)}_k}{k^{c+1}} \nonumber \\ &+& 2(-)^{a+1} \zeta(c+1) \sum^{\infty}_{k=1} \frac{h_k}{k^a} + 2(-)^a \sum^{\infty}_{k=1}
\frac{H^{(c+1)}_k h_k}{k^a} \Big)~.
\end{eqnarray}

In order to calculate corresponding Euler sums where b$>1$ an additional expression is needed. It follows:
\begin{eqnarray}
\sum^{\infty}_{k=1} \frac{1}{(2k-1)^c} \left( \sum^{\infty}_{i=1} \frac{H^{(a)}_i}{(2i+2k-1)^b} \right) = \sum^{\infty}_{i=1}
H^{(a)}_i \left( \sum^{\infty}_{i=1} \frac{1}{(2k-1)^c (2i+2k-1)^b} \right)~.
\end{eqnarray}
with $a,b,c \in \mathbb{N}$, $a+b+c=$even. Partial fraction decomposition of the denominator result in new double valued help functions. We define: 
\begin{eqnarray}
F(i) = \sum^{\infty}_{k=1} \frac{H^{(a)}_k}{(2i+2k-1)^b}~.
\end{eqnarray}
From this we get
\begin{eqnarray}
F(i+1) = \sum^{\infty}_{k=1} \frac{H^{(a)}_k}{(2i+2k+1)^a} = \sum^{\infty}_{k=1} \frac{H^{(a)}_{k-1}}{(2i+2k-1)^b} = \sum^{\infty}_{k=1} \frac{H^{(a)}_k}{(2i+2k-1)^b} -
\sum^{\infty}_{k=1} \frac{1}{k^a(2i+2k-1)^b}~, 
\end{eqnarray}
or
\begin{eqnarray}
F(i+1) - F(i) = - \sum^{\infty}_{k=1} \frac{1}{k^a(2i+2k-1)^b}~. 
\end{eqnarray}
This way we result in an inhomogeneous difference equation of first order with a constant coefficient. The solution is simply obtained by a summation over the
inhomogeneity. This way, we obtain:  
\begin{eqnarray}
F(i) = \sum^{\infty}_{k=1} \frac{H^{(a)}_k}{(2k+1)^b} - \sum^{i-1}_{k=1} \left( \sum^{\infty}_{m=1}\frac{1}{m^a(2m+2k-1)^b} \right)~.
\end{eqnarray}
As an example it holds for $a=2, b=2$ and $a=3, b=2$ which are needed two calculated all Euler sums of order 8 belonging to the fourth and fifth family:

\subsection{Lemma 4b}
\begin{eqnarray}
\sum^{\infty}_{k=1} \frac{H^{(2)}_k}{(2i+2k-1)^2} = \sum^{\infty}_{k=1}\frac{H^{(2)}_k}{(2k+1)^2} -8ln(2)h^{(3)}_{i-1} -4\zeta(2)h^{(2)}_{i-1} +
8\sum^{i-1}_{k=1}\frac{h_k}{(2k-1)^3} + 4\sum^{i-1}_{k=1}\frac{h^{(2)}_k}{(2k-1)^2}~.
\end{eqnarray}

\subsection{proof}
We define:
\begin{eqnarray}
F(i) = \sum^{\infty}_{k=1} \frac{H^{(2)}_k}{(2i+2k-1)^2}~.
\end{eqnarray}
From this we get:
\begin{eqnarray}
F(i+1) = \sum^{\infty}_{k=1} \frac{H^{(2)}_k}{(2i+2k+1)^2} = \sum^{\infty}_{k=1} \frac{H^{(2)}_{k-1}}{(2i+2k-1)^2} = F(i) -
\sum^{\infty}_{k=1} \frac{1}{i^2(2i+2k-1)^2}~.
\end{eqnarray}
From partial fraction decomposition it follows:
\begin{eqnarray}
F(i+1) - F(i) &=& -\frac{1}{(2k-1)^2}\left( \zeta(2) - 4\sum^{\infty}_{i=1}\frac{1}{i(2i+2k+1)} + 4\sum^{\infty}_{i=1}\frac{1}{(2i+2k+1)^2} \right)
\nonumber \\ &=& - \frac{\zeta(2)}{(2k-1)^2} +\frac{4}{(2k-1)^2} \left( \frac{2h_k}{(2k-1)} - \frac{2ln(2)}{(2k-1)} \right) - \frac{3\zeta(2)}{(2k-1)^2} +
\frac{4h^{(2)}_k}{(2k-1)^2} \nonumber \\ &=&
- \frac{8ln(2)}{(2k-1)^3} - \frac{4\zeta(2)}{(2k-1)^2} + \frac{8h_k}{(2k-1)^3} + \frac{4h^{(2)}_k}{(2k-1)^2}~.
\end{eqnarray}
The solution of this difference equation is given by Eq.~(171). Thus the lemma is proved.

\subsection{Lemma 4c}
\begin{eqnarray}
\sum^{\infty}_{k=1} \frac{H^{(3)}_k}{(2i+2k-1)^2} &=& \sum^{\infty}_{k=1}\frac{H^{(3)}_k}{(2k+1)^2} + 24ln(2)h^{(4)}_{i-1} + 10\zeta(2)h^{(3)}_{i-1} -
\zeta(3) h^{(2)}_{i-1} + 8\sum^{i-1}_{k=1}\frac{h_k}{(2k-1)^3} \nonumber \\ &-& 24\sum^{i-1}_{k=1}\frac{h_k}{(2k-1)^4} -8\sum^{i-1}_{k=1}\frac{h^{(2)}_k}{(2k-1)^3}~.
\end{eqnarray}

\subsection{proof}
We define:
\begin{eqnarray}
F(i) = \sum^{\infty}_{k=1} \frac{H^{(3)}_k}{(2i+2k-1)^2}~.
\end{eqnarray}
From this we get:
\begin{eqnarray}
F(i+1) = \sum^{\infty}_{k=1} \frac{H^{(3)}_k}{(2i+2k+1)^2} = \sum^{\infty}_{k=1} \frac{H^{(3)}_{k-1}}{(2i+2k-1)^2} = F(i) -
\sum^{\infty}_{k=1} \frac{1}{i^3(2i+2k-1)^2}~.
\end{eqnarray}
From partial fraction decomposition it follows:
\begin{eqnarray}
F(i+1) - F(i) &=& -\frac{\zeta(3)}{(2k-1)^2} + \frac{4\zeta(2)}{(2k-1)^3} - \frac{12}{(2k-1)^3} \left( \frac{2h_k}{(2k-1)} - \frac{2ln(2)}{(2k-1)} \right) +
\frac{8}{(2k-1)^3} \left( \frac{3}{4}\zeta(2) - h^{(2)}_k  \right) \nonumber \\ &=&
-\frac{\zeta(3)}{(2k-1)^2} + \frac{10\zeta(2)}{(2k-1)^3} + \frac{24ln(2)}{(2k-1)^4} - \frac{24h_k}{(2k-1)^4} -\frac{8h^{(2)}_k}{(2k-1)^3}~.
\end{eqnarray}
The solution of this difference equation is given by Eq.~(175). Thus the lemma is proved. This procedure can be applied for all higher order Euler sums belonging
to this family.

\subsection{examples}
The order 4 case:
\begin{eqnarray}
\sum^{\infty}_{k=1} \frac{h^2_k}{(2k-1)^2} &=& \frac{15}{32}\zeta(4) + \frac{7}{8}ln(2)\zeta(3) + \frac{3}{4}\left(ln(2)\right)^2 \zeta(2) - \frac{1}{8}\sum^{\infty}_{k=1}
\frac{h_k}{k^3}~.
\end{eqnarray}

The order 6 case with $a=b=1$:
\begin{eqnarray}
\sum^{\infty}_{k=1} \frac{h^2_k}{(2k-1)^4} &=& \frac{147}{512}\zeta(6) + \frac{31}{32}ln(2)\zeta(5) + \frac{15}{16}\left(ln(2)\right)^2 \zeta(4) - \frac{7}{64} \zeta(3)^2 -
\frac{3}{16}ln(2)\zeta(2)\zeta(3) - \frac{1}{32}\sum^{\infty}_{k=1} \frac{h_k}{k^5}~.
\end{eqnarray}

The order 6 case:
\begin{eqnarray}
\sum^{\infty}_{k=1} \frac{h_k h^{(2)}_k}{(2k-1)^3} = \frac{369}{1024}\zeta(6) + \frac{31}{64}ln(2)\zeta(5) + \frac{9}{32}ln(2)\zeta(2)\zeta(3) -
\frac{1}{64}\sum^{\infty}_{k=1} \frac{h_k}{k^5}~,
\end{eqnarray}

\begin{eqnarray}
\sum^{\infty}_{k=1} \frac{h^{(2)}_k h^{(2)}_k}{(2k-1)^2} = \frac{1737}{1536}\zeta(6) + \frac{7}{64}\zeta(3)^2~,
\end{eqnarray}

\begin{eqnarray}
\sum^{\infty}_{k=1} \frac{h_k h^{(3)}_k}{(2k-1)^2} = \frac{357}{1024}\zeta(6) + \frac{31}{64}ln(2)\zeta(5) + \frac{35}{256}\zeta(3)^2
+ \frac{3}{8}ln(2)\zeta(2)\zeta(3) - \frac{1}{64}\sum^{\infty}_{k=1} \frac{h_k}{k^5}~.
\end{eqnarray}

The order 8 case with $a=b=1$:
\begin{eqnarray}
\sum^{\infty}_{k=1} \frac{h^2_k}{(2k-1)^6} &=& \frac{265}{1024}\zeta(8) + \frac{127}{128}ln(2)\zeta(7) + \frac{63}{64}\left(ln(2)\right)^2 \zeta(6) - 
\frac{19}{128} \zeta(3)\zeta(5) - \frac{3}{64}ln(2)\zeta(2)\zeta(5) - \frac{15}{64}ln(2)\zeta(3)\zeta(4)  \nonumber \\ &+& \frac{3}{256}\zeta(2)\zeta(3)^2 -
\frac{1}{128}\sum^{\infty}_{k=1} \frac{h_k}{k^7}~.
\end{eqnarray}

The order 8 case:
\begin{eqnarray}
\sum^{\infty}_{k=1} \frac{h_k h^{(5)}_k}{(2k-1)^2} &=& \frac{10249}{12288}\zeta(8) + \frac{127}{256}ln(2)\zeta(7) - \frac{57}{128}\zeta(3)\zeta(5)
+\frac{39}{64}ln(2)\zeta(2)\zeta(5) -\frac{15}{64}ln(2)\zeta(3)\zeta(4) \nonumber \\ &+& \frac{3}{256}\zeta(2)\zeta(3)^2 -
\frac{1}{256}\sum^{\infty}_{k=1} \frac{h_k}{k^7} + \frac{221}{1024} \sum^{\infty}_{k=1} \frac{H^{(2)}_k}{k^6}~.
\end{eqnarray}

\begin{eqnarray}
\sum^{\infty}_{k=1} \frac{h_k h^{(4)}_k}{(2k-1)^3} &=& \frac{465}{1024}\zeta(8) + \frac{127}{256}ln(2)\zeta(7) - \frac{19}{256}\zeta(3)\zeta(5) +
\frac{15}{64}ln(2)\zeta(2)\zeta(5) + \frac{15}{128}ln(2)\zeta(3)\zeta(4) \\ &-& \frac{3}{512}\zeta(2)\zeta(3)^2 -
\frac{1}{256}\sum^{\infty}_{k=1} \frac{h_k}{k^7}~.
\end{eqnarray}

With
\begin{eqnarray}
&=&\sum^{\infty}_{k=1} \frac{1}{(2k-1)^4} \left( \sum^{\infty}_{i=1} \frac{H^{(2)}_i}{(2i+2k-1)^2} \right) = \sum^{\infty}_{i=1} H^{(2)}_i \left(
\sum^{\infty}_{k=1} \frac{1}{(2k-1)^4(2i+2k-1)^2} \right) \nonumber \\ &=& \sum^{\infty}_{k=1} \frac{1}{(2k-1)^4} 
\left( \sum^{\infty}_{i=1}\frac{H^{(2)}_i}{(2i+1)^2} -8ln(2)h^{(3)}_{k-1} -4\zeta(2)h^{(2)}_{k-1} +
8\sum^{k-1}_{i=1}\frac{h_i}{(2i-1)^3} + 4 \sum^{k-1}_{i=1}\frac{h^{(2)}_i}{(2i-1)^2} \right)~.
\end{eqnarray}
it follows by partial fraction decomposition
\begin{eqnarray}
\sum^{\infty}_{k=1} \frac{1}{(2k-1)^4(2i+2k-1)^2} = \frac{15}{64}\zeta(4)\frac{1}{i^2} - \frac{7}{32}\zeta(3)\frac{1}{i^3} + \frac{3}{16}\zeta(2)\frac{1}{i^4} -
\frac{1}{8}\zeta(2)\frac{h_i}{i^5} - \frac{1}{16}\zeta(2)\frac{h^{(2)}_i}{i^4}~.
\end{eqnarray}
Using this expression we get:
\begin{eqnarray}
\sum^{\infty}_{k=1}\frac{h^{(2)}_k h^{(4)}_k}{(2k-1)^2} &=& \frac{13465}{6144}\zeta(8) + \frac{127}{64}ln(2)\zeta(7) - \frac{63}{256}\zeta(3)\zeta(5) +
\frac{1}{128}\zeta(2)\zeta(3)^2 + \frac{105}{64}ln(2)\zeta(3)\zeta(4) \nonumber \\ &-& 2ln(2) \sum^{\infty}_{k=1}\frac{h^{(3)}_k}{(2k-1)^4} -
2\sum^{\infty}_{k=1}\frac{h_k h^{(4)}_k}{(2k-1)^3} - \frac{1}{32}\sum^{\infty}_{k=1}\frac{H^{(2)}_k h_k}{k^5} +
\frac{1}{64}\sum^{\infty}_{k=1}\frac{H^{(2)}_k h^{(2)}_k}{k^4}~.
\end{eqnarray}
As all Euler sums on the right side are known it follows:
\begin{eqnarray}
\sum^{\infty}_{k=1}\frac{h^{(2)}_k h^{(4)}_k}{(2k-1)^2} &=& \frac{221}{192}\zeta(8) + \frac{3}{16}\zeta(3)\zeta(5) - \frac{9}{256}\zeta(2)\zeta(3)^2 -
\frac{17}{512} \sum^{\infty}_{k=1} \frac{H^{(2)}_k}{k^6}~. 
\end{eqnarray}
Furthermore it follows by a simple rearrangement of the double summation:
\begin{eqnarray}
\sum^{\infty}_{k=1}\frac{h^{(2)}_k h^{(2)}_k}{(2k-1)^4} &=& \frac{15}{8}\zeta(4)\sum^{\infty}_{k=1}\frac{h^{(2)}_k}{(2k-1)^2} + 
2\sum^{\infty}_{k=1}\frac{h^{(2)}_k}{(2k-1)^6} - \sum^{\infty}_{k=1}\frac{h^{(4)}_k}{(2k-1)^4} - 2\sum^{\infty}_{k=1}\frac{h^{(2)}_k h^{(4)}_k}{(2k-1)^2}~.
\end{eqnarray}
From this we get:
\begin{eqnarray}
\sum^{\infty}_{k=1}\frac{h^{(2)}_k h^{(2)}_k}{(2k-1)^4} &=& \frac{691}{1536}\zeta(8) + \frac{3}{8}\zeta(3)\zeta(5) + \frac{9}{128}\zeta(2)\zeta(3)^2 
- \frac{17}{256} \sum^{\infty}_{k=1}\frac{H^{(2)}_k}{k^6}~.
\end{eqnarray}

The remaining Euler sums of order 8 of this type can be obtained by explicit calculation of the following expression:
\begin{eqnarray}
&=&\sum^{\infty}_{k=1} \frac{1}{(2k-1)^3} \left( \sum^{\infty}_{i=1} \frac{H^{(3)}_i}{(2i+2k-1)^2} \right) = \sum^{\infty}_{i=1} H^{(3)}_i \left(
\sum^{\infty}_{k=1} \frac{1}{(2k-1)^3(2i+2k-1)^2} \right) \nonumber \\ &=& \sum^{\infty}_{k=1} \frac{1}{(2k-1)^3}
\Big( \sum^{\infty}_{i=1}\frac{H^{(3)}_i}{(2i+1)^2} + 24ln(2)h^{(4)}_{k-1} + 10\zeta(2)h^{(3)}_{k-1} - \zeta(3) h^{(2)}_{k-1}
\nonumber \\ &-& 24\sum^{k-1}_{i=1}\frac{h_i}{(2i-1)^4} - 8\sum^{k-1}_{i=1}\frac{h^{(2)}_i}{(2i-1)^3} \Big) \sim 
8\sum^{\infty}_{k=1}\frac{h^{(2)}_k h^{(3)}_k}{(2k-1)^3}
\end{eqnarray}
and by use of the following identity:
\begin{eqnarray}
\sum^{\infty}_{k=1}\frac{h^{(3)}_k h^{(3)}_k}{(2k-1)^2} &=& \frac{7}{8}\zeta(3)\sum^{\infty}_{k=1}\frac{h^{(2)}_k}{(2k-1)^3} + 
\frac{7}{8}\zeta(3)\sum^{\infty}_{k=1}\frac{h^{(3)}_k}{(2k-1)^2} + \sum^{\infty}_{k=1}\frac{h^{(3)}_k}{(2k-1)^5}  \nonumber \\ &-& 
\sum^{\infty}_{k=1}\frac{h^{(5)}_k}{(2k-1)^3} +\sum^{\infty}_{k=1}\frac{h^{(2)}_k}{(2k-1)^6} -2 \sum^{\infty}_{k=1}\frac{h^{(2)}_k h^{(3)}_k}{(2k-1)^3}~,
\end{eqnarray}
with
\begin{eqnarray}
\sum^{\infty}_{k=1}\frac{h^{(2)}_k h^{(3)}_k}{(2k-1)^3} &=& -\frac{1505}{4096}\zeta(8) + \frac{503}{512}\zeta(3)\zeta(5) + \frac{99}{512}\zeta(2)\zeta(3)^2 - 
\frac{255}{1024} \sum^{\infty}_{k=1}\frac{H^{(2)}_k}{k^6}~.
\end{eqnarray}
Last but not least it follows:
\begin{eqnarray}
\sum^{\infty}_{k=1}\frac{h_k h^{(3)}_k}{(2k-1)^4} &=& -\frac{2659}{12288}\zeta(8) + \frac{127}{256}ln(2)\zeta(7) + \frac{369}{512}\zeta(3)\zeta(5) 
- \frac{15}{64}ln(2)\zeta(2)\zeta(5) + \frac{45}{64}ln(2)\zeta(3)\zeta(4) \nonumber \\ &-& \frac{15}{256}\zeta(2)\zeta(3)^2 -
\frac{1}{256}\sum^{\infty}_{k=1} \frac{h_k}{k^3} - \frac{221}{1024} \sum^{\infty}_{k=1}\frac{H^{(2)}_k}{k^6}
\end{eqnarray}
and 
\begin{eqnarray}
\sum^{\infty}_{k=1}\frac{h_k h^{(2)}_k}{(2k-1)^5} &=& \frac{141}{512}\zeta(8) + \frac{127}{256}ln(2)\zeta(7) + \frac{29}{256}\zeta(3)\zeta(5) 
- \frac{15}{128}ln(2)\zeta(2)\zeta(5) + \frac{15}{64}ln(2)\zeta(3)\zeta(4) \nonumber \\ &-& \frac{15}{512}\zeta(2)\zeta(3)^2 -
\frac{1}{256}\sum^{\infty}_{k=1} \frac{h_k}{k^7} - \frac{17}{512} \sum^{\infty}_{k=1}\frac{H^{(2)}_k}{k^6}~.
\end{eqnarray}

This way all higher orders of this type of nonlinear Euler sums can be calculated explictly in terms of special linear Euler sums and corresponding zeta values. 

\section{Fifth family}

The following type of nonlinear Euler sums can be expressed in terms of linear Euler sums and zeta values:
\begin{eqnarray}
\sum^{\infty}_{k=1} \frac{H^{(a)}_k h^{(b)}_k}{(2k-1)^c}
\end{eqnarray}

with $a,b,c \in \mathbb{N}$, $a+b+c = $even.

The following identity holds:
\subsection{Lemma 5a}
\begin{eqnarray}
\sum^{\infty}_{k=1} \frac{H_k h_k}{(2k-1)^c} &=& ln(2)\sum^{\infty}_{k=1}\frac{H_k}{(2k-1)^c} + \sum^{c-2}_{n=1}\frac{(-)^{n+1}}{2^{n+1}}\zeta(n+1)
\sum^{\infty}_{k=1}\frac{H_k}{(2k-1)^{c-n}} \nonumber \\ &+& \frac{1}{2^{c-1}} 
\left( \sum^{\infty}_{k=1} \frac{2ln(2)}{k^{c-1}(2k-1)} - ln(2)\frac{h_k}{k^c} - \frac{h_k}{k^c(2k-1)} + \frac{h^{(2)}_k}{2k^c} + \frac{h^2_k}{2k^c} \right)
\end{eqnarray}
with $c = 2n \in \mathbb{N}$.

\subsection{proof}
We start with the following double valued help function: 
\begin{eqnarray}
\sum^{\infty}_{k=1} \frac{H_k}{(2k-1)(2i+2k-1)} = \frac{2ln(2)}{2i-1} - ln(2)\frac{h_i}{i} - \frac{h_i}{i(2i-1)} + \frac{h^{(2)}_i}{2i} + \frac{h^2_i}{2i}~. 
\end{eqnarray}
This identity represents the corresponding solution of the inhomogeneous difference equation of first order, which is obtained from the definition:
\begin{eqnarray}
F(i) = \sum^{\infty}_{k=1} \frac{H_k}{(2k-1)(2i+2k-1)}~.
\end{eqnarray}
From this we get
\begin{eqnarray}
F(i+1) &=& \sum^{\infty}_{k=1} \frac{H_k}{(2k-1)(2i+2k+1)} = \frac{1}{2i+2}\sum^{\infty}_{k=1} \frac{H_k}{2k-1} - \frac{1}{2i+2}\sum^{\infty}_{k=1} \frac{H_k}{2i+2k+1}
\nonumber \\ &=& \frac{1}{2i+2}\sum^{\infty}_{k=1} \frac{H_k}{2k-1} - \frac{1}{2i+2}\sum^{\infty}_{k=1} \frac{H_{k-1}}{2i+2k-1} \nonumber \\ &=&
\frac{1}{2i+2}\sum^{\infty}_{k=1} \frac{H_k}{2k-1} - \frac{1}{2i+2}\sum^{\infty}_{k=1} \frac{H_{k}}{2i+2k-1} + \frac{1}{2i+2}\sum^{\infty}_{k=1} \frac{1}{k(2i+2k-1)}
\nonumber \\ &=&  \frac{i}{i+1} F(i) + \frac{1}{i+1} \left( \frac{h_i}{2i-1} - \frac{2ln(2)}{2i-1} \right)~, 
\end{eqnarray}
or
\begin{eqnarray}
F(i+1) - \frac{i}{i+1} F(i) = \frac{1}{i+1} \left( \frac{h_i}{2i-1} - \frac{2ln(2)}{2i-1} \right)~. 
\end{eqnarray}
The solution is simply obtained by a summation over the inhomogeneity. Dividing now Eq.~(201) on both sides by $i^{c-1}$ und summing
up over i from 1 to infinity we get:
\begin{eqnarray}
\sum^{\infty}_{i=1} \frac{1}{i^{c-1}} \left( \sum^{\infty}_{k=1} \frac{H_k}{(2k-1)(2i+2k+1)} \right) = \sum^{\infty}_{k=1} \frac{H_k}{2k-1}
\left( \sum^{\infty}_{i=1} \frac{1}{i^{c-1}(2i+2k+1)} \right)~. 
\end{eqnarray}
Inserting Eq.~(201) on the left side of Eq.~(205), applying a partial fraction decomposition of the denominator on the right side of Eq.~(205) and performing some
standard algebraic manipulations Eq.~(200) results. Thus lemma 5a is proved. 

To be able to calculate corresponding Euler sums of this family where b$>1$ or $a>1$ lemma4 and an additional expression is needed. It follows:
\begin{eqnarray}
\sum^{\infty}_{k=1} \frac{h^{(b)}_k}{(2k-1)^c} \left( \sum^{\infty}_{i=1} \frac{1}{(i+k)^a} \right) = \sum^{\infty}_{i=1}
\left( \sum^{\infty}_{k=1} \frac{h^{(b)}_k}{(i+k)^a (2k-1)^c} \right) ~.
\end{eqnarray}
We define
\begin{eqnarray}
F(i) = \sum^{\infty}_{k=1} \frac{h^{(b)}_k}{(i+k)^a}~.
\end{eqnarray}
From this we get
\begin{eqnarray}
F(i+1) = \sum^{\infty}_{k=1} \frac{h^{(b)}_k}{(i+k+1)^a} = \sum^{\infty}_{k=1} \frac{h^{(b)}_{k-1}}{(i+k)^a} = \sum^{\infty}_{k=1} \frac{h^{(b)}_k}{(i+k)^a} -
\sum^{\infty}_{k=1} \frac{1}{(2k-1)^b(i+k)^a}~, 
\end{eqnarray}
or
\begin{eqnarray}
F(i+1) - F(i) = - \sum^{\infty}_{k=1} \frac{1}{(2k-1)^b(i+k)^a}~.
\end{eqnarray}
This way we result in an inhomogeneous difference equation of first order with a constant coefficient. The solution is simply obtained by a summation over the
inhomogeneity. Therefore we obtain:  
\begin{eqnarray}
F(i) = \sum^{\infty}_{k=1} \frac{h^{(b)}_k}{(k+1)^a} - \sum^{i-1}_{k=1} \left( \sum^{\infty}_{m=1}\frac{1}{(2m-1)^b(k+m)^a} \right)~.
\end{eqnarray}

As an example it holds for $a=2, b=1$:
\subsection{Lemma 5c}
\begin{eqnarray}
\sum^{\infty}_{k=1} \frac{h_k}{(i+k)^2} = \frac{7}{4}\zeta(3) + \zeta(2)h_i - 4ln(2)h^{(2)}_i - 2\sum^{i}_{k=1} \frac{H_{k-1}}{(2k-1)^2} - 
\sum^{i}_{k=1} \frac{H^{(2)}_{k-1}}{(2k-1)}~.
\end{eqnarray}

\subsection{proof}
It follows first:
\begin{eqnarray}
\sum^{\infty}_{k=1} \frac{h_k}{(i+k)^2} &=& \sum^{\infty}_{k=1} \frac{h_k}{(i+k)^2} - \sum^{i-1}_{k=1} \left( \sum^{\infty}_{m=1}\frac{1}{(2m-1)(k+m)^2} \right)
\nonumber \\ &=& \frac{7}{4}\zeta(3) - 4ln(2) + \zeta(2) - \sum^{i-1}_{k=1}\frac{1}{(2k+1)} \left( \sum^{\infty}_{m=1}\frac{2}{(2m-1)(k+m)} -
\sum^{\infty}_{m=1}\frac{2}{(k+m)^2} \right) \nonumber \\ &=&  \frac{7}{4}\zeta(3) - 4ln(2) + \zeta(2) - \sum^{i-1}_{k=1}\frac{2H_k}{(2k+1)^2} - 
\sum^{i-1}_{k=1}\frac{4ln(2)}{(2k+1)^2} + \sum^{i-1}_{k=1}\frac{\zeta(2)-H^{(2)}_k}{(2k+1)} \nonumber \\ &=& 
\frac{7}{4}\zeta(3) - 4ln(2) + \zeta(2) + \zeta(2) \left(h_i -1 \right) - 4ln(2)  \left(h^{(2)}_i -1 \right) -  2\sum^{i}_{k=1} \frac{H_{k-1}}{(2k-1)^2}
-\sum^{i}_{k=1} \frac{H^{(2)}_{k-1}}{(2k-1)} \nonumber \\ &=&
\frac{7}{4}\zeta(3) + \zeta(2)h_i - 4ln(2)h^{(2)}_i - 2\sum^{i}_{k=1} \frac{H_{k-1}}{(2k-1)^2} - \sum^{i}_{k=1} \frac{H^{(2)}_{k-1}}{(2k-1)}~.
\end{eqnarray}
Thus lemma 5c is proved.

\subsection{examples}
The order 4 case:
\begin{eqnarray}
\sum^{\infty}_{k=1} \frac{H_k h_k}{(2k-1)^2} &=& \frac{45}{32}\zeta(4) + \frac{7}{8}ln(2)\zeta(2)\zeta(3) + \frac{7}{8}\zeta(3) - \frac{3}{2}\left(ln(2)\right)^2\zeta(2)
+ \frac{3}{2}ln(2)\zeta(2) - \zeta(2) - \frac{1}{4} \sum^{\infty}_{k=1} \frac{h_k}{k^3}~.
\end{eqnarray}

The order 6 case:
\begin{eqnarray}
\sum^{\infty}_{k=1} \frac{H_k h_k}{(2k-1)^4} &=& \frac{405}{256}\zeta(6) + \frac{93}{32}ln(2)\zeta(5) + \frac{31}{32}\zeta(5) - \frac{15}{8}\left(ln(2)\right)^2\zeta(4)
+ \frac{15}{8}ln(2)\zeta(4) - \frac{15}{16}\zeta(4) - \frac{91}{128}\zeta(3)^2 \nonumber \\ &-& \frac{9}{8}ln(2)\zeta(2)\zeta(3) - \frac{3}{16}\zeta(2)\zeta(3) -
\frac{7}{4}ln(2)\zeta(3) + \frac{7}{8}\zeta(3) + \frac{3}{2}ln(2)\zeta(2) - \zeta(2) + \frac{1}{4} \sum^{\infty}_{k=1} \frac{h_k}{k^3} \nonumber \\ &-&
\frac{1}{8} \sum^{\infty}_{k=1} \frac{h_k}{k^5}~.
\end{eqnarray}

For example in the order 6 case with $a=b=c=2$ we get:
\begin{eqnarray}
\sum^{\infty}_{k=1} \frac{H^{(2)}_k h^{(2)}_k}{(2k-1)^2} = \frac{3}{4}\zeta(2)\sum^{\infty}_{k=1} \frac{H^{(2)}_k}{(2k-1)^2} - \sum^{\infty}_{k=1} \frac{H^{(2)}_k}{(2k-1)^2}
\left( \sum^{\infty}_{i=1} \frac{1}{(2i+2k-1)^2} \right) ~.
\end{eqnarray}

Calculating explicitly this expression it follows:
\begin{eqnarray}
\sum^{\infty}_{k=1} \frac{H^{(2)}_k h^{(2)}_k}{(2k-1)^2} &=& -\frac{705}{128}\zeta(6) + \frac{15}{2}\zeta(4) - \frac{161}{32}\zeta(3)^2 - \frac{21}{2}\zeta(3) + 
6ln(2)\zeta(2) \nonumber \\ &-& \sum^{\infty}_{k=1} \frac{h_k}{k^3} - \frac{3}{4}\zeta(2) \sum^{\infty}_{k=1} \frac{h_k}{k^3} + \frac{1}{4}
\sum^{\infty}_{k=1} \frac{h_k}{k^5}~.
\end{eqnarray}

As further examples we present:
\begin{eqnarray}
\sum^{\infty}_{k=1} \frac{H_k h^{(3)}_k}{(2k-1)^2} &=& \frac{1035}{512}\zeta(6) - \frac{31}{32}ln(2)\zeta(5) + \frac{31}{32}\zeta(5) - \frac{75}{32}\zeta(4) +
\frac{7}{128}\zeta(3)^2 - \frac{3}{4}ln(2)\zeta(2)\zeta(3) + \frac{3}{4}\zeta(2)\zeta(3) \nonumber \\ &+& \frac{7}{4}ln(2)\zeta(3) -\frac{1}{4} \sum^{\infty}_{k=1}
\frac{h_k}{k^3}~,
\end{eqnarray}
and
\begin{eqnarray}
\sum^{\infty}_{k=1} \frac{H_k h^{(2)}_k}{(2k-1)^3} &=& \frac{585}{256}\zeta(6) - \frac{31}{32}ln(2)\zeta(5) + \frac{31}{32}\zeta(5) - \frac{75}{32}\zeta(4) -
\frac{105}{128}\zeta(3)^2 - \frac{9}{16}ln(2)\zeta(2)\zeta(3) + \frac{9}{16}\zeta(2)\zeta(3) \nonumber \\ &+& \frac{21}{8}\zeta(3) - \frac{3}{2}ln(2)\zeta(2) 
+\frac{3}{16}\zeta(2) \sum^{\infty}_{k=1} \frac{h_k}{k^3}~.
\end{eqnarray}

Next we calculate the corresponding sum of this family for $a=2, b=1, c=3$:
\begin{eqnarray}
\sum^{\infty}_{k=1} \frac{h_k}{(2k-1)^3} \left( \sum^{\infty}_{i=1} \frac{1}{(i+k)^2} \right) = \zeta(2) \sum^{\infty}_{k=1} \frac{h_k}{(2k-1)^3} - 
\sum^{\infty}_{k=1} \frac{H^{(2)}_k h_k}{(2k-1)^3}~.
\end{eqnarray}
From this we get:
\begin{eqnarray}
\sum^{\infty}_{k=1} \frac{H^{(2)}_k h_k}{(2k-1)^3} = \zeta(2) \sum^{\infty}_{k=1} \frac{h_k}{(2k-1)^3} - \sum^{\infty}_{i=1} \left( \sum^{\infty}_{k=1} 
\frac{h_k}{(2k-1)^3(i+k)^2} \right)~.
\end{eqnarray}
From the corresponding partial fraction decomposition and use lemma 5c we get: 
\begin{eqnarray}
\sum^{\infty}_{k=1} \frac{H^{(2)}_k h_k}{(2k-1)^3} &=& - \frac{165}{32}\zeta(6) - \frac{93}{8}ln(2)\zeta(5) + \frac{15}{8}\zeta(4) + \frac{217}{64}\zeta(3)^2
+ \frac{49}{8}ln(2)\zeta(2)\zeta(3) + \frac{7}{2}ln(2)\zeta(3) \nonumber \\ &-& \frac{21}{4}\zeta(3) - 6ln(2)\zeta(2) + 6\zeta(2) - \frac{1}{2}
\sum^{\infty}_{k=1} \frac{h_k}{k^3} - \frac{1}{8}\zeta(2) \sum^{\infty}_{k=1} \frac{h_k}{k^3} + \frac{5}{8} \sum^{\infty}_{k=1} \frac{h_k}{k^5}~. 
\end{eqnarray}
The corresponding sum with $a=3, b=1, c=2$ finally follows by use of Eq.~(166). 

This way all sums of the fifth family to the order of six are known. This
calculational scheme of course works for higher order sums as well. This will be shown by the last two examples which are of order eight. We start with the
following expression:
\begin{eqnarray}
\sum^{\infty}_{k=1} \frac{H_k}{(2k-1)^5} \left( \sum^{\infty}_{i=1} \frac{1}{(2i+2k-1)^2} \right) =  \frac{3}{4}\zeta(2) \sum^{\infty}_{k=1}
\frac{H_k}{(2k-1)^5} - \sum^{\infty}_{k=1} \frac{H_k h^{(2)}_k}{(2k-1)^5}~.
\end{eqnarray}
From this we get:
\begin{eqnarray}
\sum^{\infty}_{k=1} \frac{H_k h^{(2)}_k}{(2k-1)^5} = \frac{3}{4}\zeta(2)\sum^{\infty}_{k=1}
\frac{H_k}{(2k-1)^5} - \sum^{\infty}_{i=1} \left( \sum^{\infty}_{k=1} \frac{H_k}{(2k-1)^5(2i+2k-1)^2} \right)~. 
\end{eqnarray}
By partial fraction decomposition, by use of Eq.~(205) and by use of Eq.~(167) for $a=1,b=2$ the corresponding Euler sum can be calculated explicitly. First
it follows:
\begin{eqnarray}
\sum^{\infty}_{k=1} \frac{H_k}{(2i+2k-1)^2} &=& \sum^{\infty}_{k=1} \frac{H_k}{(2i+1)^2} - \frac{2ln(2)}{(2i-1)^2} - \frac{3}{2}\zeta(2)\frac{1}{(2i-1)} +
\frac{3}{2}\zeta(2)h_i + 2ln(2)h^{(2)}_i -2h^{(3)}_i - 2h_i h^{(2)}_i \nonumber \\ &+& \frac{2h_i}{(2i-1)^2} +  \frac{2h^{(2)}_i}{(2i-1)}~. 
\end{eqnarray}
The proof is analogous to that of lemma 4b or lemma 4c, where the cases $a=2,b=2$ and $a=3,b=2$ have been proved. This way we result after standard but tedious
algebraic manipulations at the following expression:
\begin{eqnarray}
\sum^{\infty}_{k=1} \frac{H_k h^{(2)}_k}{(2k-1)^5} &=& \frac{3}{16}\zeta(8) - \frac{127}{128}ln(2)\zeta(7) + \frac{127}{128}\zeta(7) - \frac{27}{16}\zeta(6) +
\frac{643}{256}\zeta(3)\zeta(5) - \frac{15}{64}ln(2)\zeta(2)\zeta(5) \nonumber \\ &+& \frac{15}{64}\zeta(2)\zeta(5) + \frac{31}{32}\zeta(5) -
\frac{15}{32}ln(2)\zeta(3)\zeta(4) + \frac{15}{32}\zeta(3)\zeta(4) - \frac{75}{32}\zeta(4) - \frac{21}{64}\zeta(2)\zeta(3)^2 \nonumber \\ &-& \frac{7}{32}\zeta(3)^2
+ \frac{9}{16}\zeta(2)\zeta(3) + \frac{21}{8}\zeta(3) - \frac{3}{2}ln(2)\zeta(2) + \frac{3}{64}\zeta(2) \sum^{\infty}_{k=1}\frac{h_k}{k^5} -
\frac{153}{128}\sum^{\infty}_{k=1}\frac{H^{(2)}_k}{k^6}~,
\end{eqnarray}
and as a last example for this family:
\begin{eqnarray}
\sum^{\infty}_{k=1} \frac{H_k h^{(5)}_k}{(2k-1)^2} &=& \frac{5853}{1024}\zeta(8) - \frac{127}{128}ln(2)\zeta(7) + \frac{127}{128}\zeta(7) - \frac{441}{128}\zeta(6) -
\frac{1077}{256}\zeta(3)\zeta(5) - \frac{39}{32}ln(2)\zeta(2)\zeta(5) \nonumber \\ &+& \frac{39}{32}\zeta(2)\zeta(5) + \frac{31}{16}ln(2)\zeta(5) +
\frac{15}{32}ln(2)\zeta(3)\zeta(4) - \frac{15}{32}\zeta(3)\zeta(4) + \frac{21}{128}\zeta(2)\zeta(3)^2 \nonumber \\ &+& \frac{35}{64}\zeta(3)^2
- \frac{1}{16}\sum^{\infty}_{k=1}\frac{h_k}{k^5} + \frac{153}{128}\sum^{\infty}_{k=1}\frac{H^{(2)}_k}{k^6}~.
\end{eqnarray}

\section{Sixth family}

The following type of nonlinear Euler sums can be expressed in terms of linear Euler sums and zeta values:
\begin{eqnarray}
\sum^{\infty}_{k=1} \frac{h^{(a)}_k h^{(b)}_k}{k(2k-1)}
\end{eqnarray}

with $a,b,c \in \mathbb{N}$, $a+b+1 = $even.

The following identity holds:
\subsection{Lemma 6}
\begin{eqnarray}
\sum^{\infty}_{k=1} \frac{h^{(a)}_k h_k}{k(2k-1)} &=& \sum^{a}_{n=1} (-)^{n+1} \frac{2^{a+2-n}-1}{2^{a+1}}\zeta(a+2-n)\sum^{\infty}_{i=1}\frac{1}{i^n (2i-1)}
\nonumber \\ &+& \sum^{a-1}_{n=1}(-)^{n+1}\frac{2^{a+1-n}-1}{2^{a+1}}\zeta(a+1-n)\sum^{\infty}_{i=1}\frac{H_i}{i^n (2i-1)} 
+ln(2)\frac{(-)^{a+1}}{2^{a-1}}\sum^{\infty}_{i=1}\frac{h_i}{i^a} \nonumber \\ &+& \frac{(-)^{a+1}}{2^a}\sum^{\infty}_{i=1}\frac{h_i}{i^{a+1}(2i-1)} -
\frac{(-)^{a+1}}{2^{a-1}}\sum^{\infty}_{i=1}\frac{h^2_i}{i^a} + \frac{(-)^{a+1}}{2^a}\sum^{\infty}_{i=1}\frac{H_i h_i}{i^a} ~. 
\end{eqnarray}

\subsection{proof}
It follows by use of lemma 4a:
\begin{eqnarray}
\sum^{\infty}_{i=1} \frac{1}{2i-1} \left( \sum^{\infty}_{k=1} \frac{h^{(a)}_k}{(2k-1)(2i+2k-1)} \right) = \sum^{\infty}_{k=1} \frac{h^{(a)}_k}{2k-1} \left( \sum^{\infty}_{i=1}
\frac{1}{(2i-1)(2i+2k-1)} \right) = \frac{1}{2}\sum^{\infty}_{k=1} \frac{h^{(a)}_k h_k}{k(2k-1)}~.  
\end{eqnarray}
From this we get:
\begin{eqnarray}
\sum^{\infty}_{k=1} \frac{h^{(a)}_k h_k}{k(2k-1)} &=& \sum^{a}_{n=1} (-)^{n+1} \frac{2^{a+2-n}-1}{2^{a+1}}\zeta(a+2-n)\sum^{\infty}_{i=1}\frac{1}{i^n (2i-1)}
\nonumber \\ &+& \sum^{a-1}_{n=1}(-)^{n+1}\frac{2^{a+1-n}-1}{2^{a+1}}\zeta(a+1-n)\sum^{\infty}_{i=1}\frac{H^{(n)}_i}{i(2i-1)} 
+ \frac{(-)^{a+1}}{2^{a}} \sum^{\infty}_{i=1} \frac{1}{i(2i-1)} \sum^{i-1}_{k=1}\frac{h_k}{k^a}~. 
\end{eqnarray}
Finally, performing some standard algebraic manipulations on the last term of right side Eq.~(228) follows. Thus lemma 6 is proved.

In order to calculate corresponding Euler sums where b$>1$ additional identities of the following form are needed,
\begin{eqnarray}
\sum^{\infty}_{k=1} \frac{1}{k(2k-1)} \left( \sum^{k}_{i=1} \frac{h^{(a)}_k}{(2k-1)^b} \right) = \sum^{\infty}_{k=1} 
\frac{h^{(a)}_k h^{(b)}_k}{k(2k-1)} + \sum^{\infty}_{k=1} \frac{h^{(a+b)}_k}{k(2k-1)} -
\sum^{\infty}_{k=1} \frac{1}{k(2k-1)} \left( \sum^{k}_{i=1} \frac{h^{(b)}_k}{(2k-1)^a} \right)~.
\end{eqnarray}
These identities can reformulated by rearranging the corresponding double summations. It follows then:
\begin{eqnarray}
&& 2ln(2)\sum^{\infty}_{k=1} \frac{h^{(a)}_k}{(2k-1)^b} + \sum^{\infty}_{k=1} \frac{h^{(a)}_k}{k(2k-1)^{b+1}} - \sum^{\infty}_{k=1} \frac{h^{(a)}_k}{k(2k-1)^b}
\left(2h_k - H_k \right) = \sum^{\infty}_{k=1} \frac{h^{(a)}_k h^{(b)}_k}{k(2k-1)} + \sum^{\infty}_{k=1} \frac{h^{(a+b)}_k}{k(2k-1)} \nonumber \\ &-&
2ln(2)\sum^{\infty}_{k=1} \frac{h^{(b)}_k}{(2k-1)^a} - \sum^{\infty}_{k=1} \frac{h^{(b)}_k}{k(2k-1)^{a+1}} + \sum^{\infty}_{k=1} \frac{h^{(b)}_k}{k(2k-1)^a}
\left(2h_k - H_k \right)~.
\end{eqnarray}
This way, we get:
\begin{eqnarray}
\sum^{\infty}_{k=1} \frac{h^{(a)}_k h^{(b)}_k}{k(2k-1)} &=&  2ln(2) \left(\sum^{\infty}_{k=1} \frac{h^{(a)}_k}{(2k-1)^b} + 
\sum^{\infty}_{k=1} \frac{h^{(b)}_k}{(2k-1)^a} \right) + \sum^{\infty}_{k=1} \frac{h^{(a)}_k}{k(2k-1)^{b+1}} + \sum^{\infty}_{k=1} \frac{h^{(b)}_k}{k(2k-1)^{a+1}} 
\nonumber \\ &+& \sum^{\infty}_{k=1} \frac{H_k h^{(a)}_k}{(2k-1)^b} + \sum^{\infty}_{k=1} \frac{H_k h^{(b)}_k}{(2k-1)^a} 
-2\sum^{\infty}_{k=1} \frac{h_k h^{(a)}_k}{(2k-1)^b} - 2\sum^{\infty}_{k=1} \frac{h_k h^{(b)}_k}{(2k-1)^a} - \sum^{\infty}_{k=1} \frac{h^{(a+b)}_k}{k(2k-1)}~.
\end{eqnarray}
As one can see, on the right side of this identity only known Euler sums appear together with Euler sums belonging to the fourth and fifth family. This shows that
the corresponding calculational scheme for the family six Euler sums is also highly recursive as family five members with $a=1$ are needed together with lemma 4.

\subsection{examples}
The order 4 case with $b=1$:
\begin{eqnarray}
\sum^{\infty}_{k=1} \frac{h_k h^{(2)}_k}{k(2k-1)} = \frac{105}{64}\zeta(4) + \frac{7}{8}ln(2)\zeta(3) - \frac{3}{4}\left(ln(2)\right)^2 \zeta(2) - 
\frac{1}{8}\sum^{\infty}_{k=1} \frac{h_k}{k^3}~.
\end{eqnarray}

The order 6 case with $b=1$:
\begin{eqnarray}
\sum^{\infty}_{k=1} \frac{h_k h^{(4)}_k}{k(2k-1)} &=& \frac{975}{512}\zeta(6) + \frac{31}{32}ln(2)\zeta(5) - \frac{15}{16}\left(ln(2)\right)^2 \zeta(4) - 
\frac{7}{32}\zeta(3)^2 + \frac{3}{16}ln(2)\zeta(2)\zeta(3) \nonumber \\ &-& \frac{3}{16}\zeta(2)\sum^{\infty}_{k=1}\frac{h_k}{k^3} - \frac{1}{32}\zeta(2)
\sum^{\infty}_{k=1} \frac{h_k}{k^5}~.
\end{eqnarray}

The order 6 case:
\begin{eqnarray}
\sum^{\infty}_{k=1} \frac{h^{(2)}_k h^{(3)}_k}{k(2k-1)} =  \frac{1083}{512}\zeta(6) + \frac{63}{128}\zeta(3)^2 - \frac{21}{16}ln(2)\zeta(2)\zeta(3) + 
\frac{3}{16}\zeta(2)\sum^{\infty}_{k=1} \frac{h_k}{k^3}~.
\end{eqnarray}

The order 8 case with $b=1$:
\begin{eqnarray}
\sum^{\infty}_{k=1} \frac{h_k h^{(6)}_k}{k(2k-1)} &=& \frac{2241}{1024}\zeta(8) + \frac{127}{128}ln(2)\zeta(7) -\frac{63}{64}\left(ln(2)\right)^2 \zeta(6) - 
\frac{67}{128}\zeta(3)\zeta(5) + \frac{3}{64}ln(2)\zeta(2)\zeta(5) + \frac{15}{64}ln(2)\zeta(3)\zeta(4) \nonumber \\ &-& \frac{3}{256}\zeta(2)\zeta(3)^2 -
\frac{15}{64}\zeta(4) \sum^{\infty}_{k=1} \frac{h_k}{k^3} - \frac{3}{64}\zeta(2) \sum^{\infty}_{k=1} \frac{h_k}{k^5} - \frac{1}{128}\sum^{\infty}_{k=1} \frac{h_k}{k^7} +
\frac{17}{256}\sum^{\infty}_{k=1} \frac{H^{(2)}_k}{k^6}~.
\end{eqnarray}

Finally we present, as an example for the order 8 case, the Euler sum for $a=2, b=5$. Using Eq.~(233) it follows:
\begin{eqnarray}
\sum^{\infty}_{k=1} \frac{h^{(2)}_k h^{(5)}_k}{k(2k-1)} &=& 2ln(2) \left( \sum^{\infty}_{k=1}\frac{h^{(2)}_k}{(2k-1)^5} + \sum^{\infty}_{k=1}\frac{h^{(5)}_k}{(2k-1)^2}
\right) + \sum^{\infty}_{k=1}\frac{h^{(2)}_k}{k(2k-1)^6} + \sum^{\infty}_{k=1}\frac{h^{(5)}_k}{k(2k-1)^3} \nonumber \\ &+& \sum^{\infty}_{k=1}\frac{H_k h^{(2)}_k}{(2k-1)^5}
+ \sum^{\infty}_{k=1}\frac{H_k h^{(5)}_k}{(2k-1)^2} - 2\sum^{\infty}_{k=1}\frac{h_k h^{(2)}_k}{(2k-1)^5} -2\sum^{\infty}_{k=1}\frac{h_k h^{(5)}_k}{(2k-1)^2} -
\sum^{\infty}_{k=1}\frac{h^{(7)}_k}{k(2k-1)}~.
\end{eqnarray}
As all terms on the right side of Eq.~(238) are known explicitly, it follows:
\begin{eqnarray}
\sum^{\infty}_{k=1} \frac{h^{(2)}_k h^{(5)}_k}{k(2k-1)} &=& \frac{21263}{6144}\zeta(8) - \frac{47}{128}\zeta(3)\zeta(5) - \frac{93}{64}ln(2)\zeta(2)\zeta(5)
-\frac{33}{256}\zeta(2)\zeta(3)^2 + \frac{3}{64}\zeta(2)\sum^{\infty}_{k=1} \frac{h_k}{k^5} \nonumber \\ &+& \frac{187}{512}\sum^{\infty}_{k=1} \frac{H^{(2)}_k}{k^6}~.
\end{eqnarray}

\section{Seventh family}

The following type of nonlinear Euler sums can be expressed in terms of linear Euler sums and zeta values:
\begin{eqnarray}
\sum^{\infty}_{k=1} \frac{H^{(b)}_k h^{(a)}_k}{k(2k-1)}
\end{eqnarray}

with $a,b \in \mathbb{N}$, $a+b$ = odd.

The following identity holds for $a=1$:
\subsection{Lemma 7}

\begin{eqnarray}
\sum^{\infty}_{k=1} \frac{H^{(2a)}_k h_k}{k(2k-1)} &=& 2ln(2)\sum^{\infty}_{k=1}\frac{h_k}{k^{2a}} + \zeta(2)\zeta(2a) - \sum^{2a-2}_{n=1}(-)^{n+1}\zeta(2a-n)
\sum^{\infty}_{k=1}\frac{h_k}{k^{n+1}} + \sum^{\infty}_{k=1}\frac{h_k}{k^{2a+1}(2k-1)} \nonumber \\ &-& \sum^{\infty}_{k=1}\frac{h^{(2)}_k}{k^{2a}} -
\sum^{\infty}_{k=1}\frac{h^2_k}{k^{2a}}~.
\end{eqnarray}

\subsection{proof}
We start with the identity introduced in \cite{bra20}
\begin{eqnarray}
\sum^{\infty}_{k=1} \frac{h_k}{k(i+k)} &=& 2ln(2)\frac{h_i}{i} + \frac{H_i h_i}{i} - \frac{1}{i}\sum^{i}_{k=1}\frac{h_k}{k} =  2ln(2)\frac{h_i}{i} + \frac{H_i h_i}{i} +
\frac{1}{i}\sum^{i}_{k=1}\frac{h_k}{k(2k-1)} - \frac{2}{i} \sum^{i}_{k=1}\frac{h_k}{2k-1} \nonumber \\ &=&
2ln(2)\frac{h_i}{i} + \frac{H_i h_i}{i} - \frac{h^{(2)}_i}{i} - \frac{h^2_i}{i} + \frac{1}{i}\sum^{i}_{k=1}\frac{h_k}{k(2k-1)}~.
\end{eqnarray} 
Dividing now both sides by $i^b$ and summing up over i from 1 to $\infty$ we get:
\begin{eqnarray}
\sum^{\infty}_{i=1}\frac{1}{i^b} \left(\sum^{\infty}_{k=1} \frac{h_k}{k(i+k)} \right) &=& \sum^{\infty}_{k=1} \frac{h_k}{k} \left( \sum^{\infty}_{i=1} \frac{1}{i^b(k+i)} \right)
\nonumber \\ &=& 2ln(2) \sum^{\infty}_{i=1}\frac{h_i}{i^{b+1}} + \sum^{\infty}_{i=1}\frac{H_i h_i}{i^{b+1}} -  \frac{h^{(2)}_i}{i^{b+1}} - \frac{h^2_i}{i^{b+1}} +
\sum^{\infty}_{i=1}\frac{1}{i^{b+1}} \left( \sum^{i}_{k=1}\frac{h_k}{k(2k-1)} \right)~.
\end{eqnarray}
With
\begin{eqnarray}
\sum^{\infty}_{i=1} \frac{1}{i^b(k+i)} = (-)^{b+1} \frac{H_k}{k^b} + \sum^{b-1}_{n=1} \frac{(-)^{n+1}}{k^n}\zeta(b-n+1) 
\end{eqnarray}
it follows:
\begin{eqnarray}
\sum^{\infty}_{k=1} (-)^{b+1} \frac{H_k h_k}{k^{b+1}} + \sum^{b-1}_{n=1} (-)^{n+1}\zeta(b-n+1)\sum^{\infty}_{k=1} \frac{h_k}{k^{n+1}} &=&
2ln(2) \sum^{\infty}_{k=1}\frac{h_k}{k^{b+1}} + \sum^{\infty}_{k=1}\frac{H_k h_k}{k^{b+1}} -  \frac{h^{(2)}_k}{k^{b+1}} - \frac{h^2_k}{k^{b+1}} \nonumber \\ &+&
\zeta(2)\zeta(b+1) + \sum^{\infty}_{k=1} \frac{h_k}{k^{b+2}(2k-1)} - \sum^{\infty}_{k=1} \frac{H^{(b+1)}_k h_k}{k(2k-1)}~.
\end{eqnarray}
With $b =2a-1$ we get after some standard algebraic manipulations Eq.~(241). Thus the lemma is proved.

In order to calculate corresponding Euler sums where a,b$>1$ an additional identity is needed. With Adegoke \cite{ade16} it follows: 
\begin{eqnarray}
H^{(a)}_k h^{(b)}_k = \sum^{k}_{i=1} \frac{h^{(b)}_i}{i^a} + \sum^{k}_{i=1} \frac{H^{(a)}_{i-1}}{(2i-1)^b}~,
\end{eqnarray}
and from this we get: 
\begin{eqnarray}
\sum^{\infty}_{k=1} \frac{H^{(a)}_k h^{(b)}_k}{k(2k-1)} = \sum^{\infty}_{k=1} \frac{1}{k(2k-1)} \left( \sum^{k}_{i=1} \frac{h^{(b)}_i}{i^a} \right)
+ \sum^{\infty}_{k=1} \frac{1}{k(2k-1)} \left( \sum^{k}_{i=1} \frac{H^{(a)}_{i-1}}{(2i-1)^b} \right)~.
\end{eqnarray}
Rearranging the double summations on the right side of this expression one arrives at:
\begin{eqnarray}
\sum^{\infty}_{k=1} \frac{H^{(a)}_k h^{(b)}_k}{k(2k-1)} &=& 2ln(2) \left( \sum^{\infty}_{k=1} \frac{h^{(b)}_k}{k^a} + \sum^{\infty}_{k=1} \frac{H^{(a)}_k}{(2k-1)^b}
\right) - \sum^{\infty}_{k=1} \frac{1}{k(2k-1)} \left( \sum^{k}_{i=1} \frac{1}{i^a(2i-1)^b} \right) \nonumber \\ &+& 
\sum^{\infty}_{k=1} \frac{h^{(b)}_k}{k^{a+1}(2k-1)} + \sum^{\infty}_{k=1} \frac{H^{(a)}_k}{k(2k-1)^{b+1}} - \sum^{\infty}_{k=1} \frac{h^{(b)}_k}{k^a}
\left(2h_k - H_k \right) - \sum^{\infty}_{k=1} \frac{H^{(a)}_k}{(2k-1)^b} \left(2h_k - H_k \right)~, 
\end{eqnarray}
or
\begin{eqnarray}
\sum^{\infty}_{k=1} \frac{H^{(a)}_k h^{(b)}_k}{k(2k-1)} &=& 2ln(2) \left( \sum^{\infty}_{k=1} \frac{h^{(b)}_k}{k^a} + \sum^{\infty}_{k=1} \frac{H^{(a)}_k}{(2k-1)^b}
\right) - \sum^{\infty}_{k=1} \frac{1}{k(2k-1)} \left( \sum^{k}_{i=1} \frac{1}{i^a(2i-1)^b} \right) + 
\sum^{\infty}_{k=1} \frac{h^{(b)}_k}{k^{a+1}(2k-1)} \nonumber \\ &+& \sum^{\infty}_{k=1} \frac{H^{(a)}_k}{k(2k-1)^{b+1}} 
- 2\sum^{\infty}_{k=1} \frac{h_k h^{(b)}_k}{k^a} + \sum^{\infty}_{k=1} \frac{H_k h^{(b)}_k}{k^a} 
- 2\sum^{\infty}_{k=1} \frac{H^{(a)}_k h_k}{(2k-1)^b} + \sum^{\infty}_{k=1} \frac{H^{(a)}_k H_k}{(2k-1)^b}~. 
\end{eqnarray}
Again, on the right side of this identity only known Euler sums appear together with Euler sums belonging to the first, second and fifth family. Furthermore,
Euler sums appear on the right side of Eq.~(249) that have been discussed in \cite{bra20}.  This way, the corresponding calculational scheme for the seventh family
of Euler sums is highly recursive.

\subsection{examples}
The order 4 case:
\begin{eqnarray}
\sum^{\infty}_{k=1} \frac{H^{(2)}_k h_k}{k(2k-1)} &=& -\frac{25}{4}\zeta(4) + \frac{7}{2}ln(2)\zeta(3) - \frac{7}{2}\zeta(3) + 4\zeta(2)~.
\end{eqnarray}

\begin{eqnarray}
\sum^{\infty}_{k=1} \frac{H_k h^{(2)}_k}{k(2k-1)} &=& \frac{15}{8}\zeta(4) - \frac{35}{4}ln(2)\zeta(3) + \frac{21}{4}\zeta(3) + 3\left(ln(2)\right)^2\zeta(2) -
3ln(2)\zeta(2) + \frac{3}{2}\sum^{\infty}_{k=1} \frac{h_k}{k^3}~.
\end{eqnarray}

The order 6 case:
\begin{eqnarray}
\sum^{\infty}_{k=1} \frac{H^{(4)}_k h_k}{k(2k-1)} &=& -\frac{439}{32}\zeta(6) + \frac{31}{2}ln(2)\zeta(5) - \frac{31}{2}\zeta(5) + \frac{49}{8}\zeta(3)^2 -
7ln(2)\zeta(2)\zeta(3) + 7\zeta(2)\zeta(3) - 14\zeta(3) - 16\zeta(2) \nonumber \\ &-& 4\sum^{\infty}_{k=1} \frac{h_k}{k^3} + \zeta(2)\sum^{\infty}_{k=1}
\frac{h_k}{k^3} + \sum^{\infty}_{k=1} \frac{h_k}{k^5}~.
\end{eqnarray}

\begin{eqnarray}
\sum^{\infty}_{k=1} \frac{H_k h^{(4)}_k}{k(2k-1)} &=& \frac{465}{256}\zeta(6) - \frac{279}{16}ln(2)\zeta(5) + \frac{155}{16}\zeta(5) + 
\frac{15}{4}\left(ln(2)\right)^2\zeta(4) - \frac{15}{4}ln(2)\zeta(4) + \frac{49}{64}\zeta(3)^2 \nonumber \\ &+& \frac{39}{8}ln(2)\zeta(2)\zeta(3) - 
\frac{9}{4}\zeta(2)\zeta(3) -\frac{3}{8}\zeta(2) \sum^{\infty}_{k=1} \frac{h_k}{k^3} + \frac{3}{4} \sum^{\infty}_{k=1} \frac{h_k}{k^5}~.
\end{eqnarray}
The two remaining sums of order 6 with $a=2, b=3$ and $a=3, b=2$ can be calculated using Eq.~(249). It follows:
\begin{eqnarray}
\sum^{\infty}_{k=1} \frac{H^{(2)}_k h^{(3)}_k}{k(2k-1)} &=& -\frac{735}{128}\zeta(6) + \frac{93}{2}ln(2)\zeta(5) - \frac{93}{4}\zeta(5) + \frac{75}{8}\zeta(4)
-\frac{49}{32}\zeta(3)^2 - \frac{91}{4}ln(2)\zeta(2)\zeta(3) + \frac{21}{2}\zeta(2)\zeta(3) \nonumber \\ &-& 7ln(2)\zeta(3) + \sum^{\infty}_{k=1} \frac{h_k}{k^3}
+\frac{13}{4}\zeta(2) \sum^{\infty}_{k=1} \frac{h_k}{k^3} - \frac{7}{4} \sum^{\infty}_{k=1} \frac{h_k}{k^5}
\end{eqnarray}
and
\begin{eqnarray}
\sum^{\infty}_{k=1} \frac{H^{(2)}_k h^{(5)}_k}{k(2k-1)} &=& -\frac{29079}{512}\zeta(8) + \frac{1905}{16}ln(2)\zeta(7) -\frac{1905}{32}\zeta(7) + \frac{441}{32}\zeta(5)
+\frac{3763}{64}\zeta(3)\zeta(5) - \frac{775}{16}ln(2)\zeta(2)\zeta(5) \nonumber \\ &+&  \frac{93}{4}\zeta(2)\zeta(5) - \frac{31}{4}ln(2)\zeta(5) -
\frac{105}{4}ln(2)\zeta(3)\zeta(4) + \frac{105}{4}\zeta(3)\zeta(4) - \frac{119}{64}\zeta(2)\zeta(3)^2 \nonumber \\ &-& \frac{35}{16}\zeta(3)^2 + \frac{15}{4}\zeta(4)
\sum^{\infty}_{k=1} \frac{h_k}{k^3} + \frac{1}{4}\sum^{\infty}_{k=1} \frac{h_k}{k^5} + \frac{25}{16}\zeta(2)\sum^{\infty}_{k=1} \frac{h_k}{k^5}
- \frac{39}{32}\sum^{\infty}_{k=1} \frac{h_k}{k^7} - \frac{2703}{128}\sum^{\infty}_{k=1} \frac{H^{(2)}_k}{k^6}~.
\end{eqnarray}

\section{Eighth family}

The following type of nonlinear Euler sums can be expressed in terms of linear Euler sums and zeta values:
\begin{eqnarray}
\sum^{\infty}_{k=1} \frac{h^{(a)}_k h^{(b)}_k}{k^c}
\end{eqnarray}
with $a,b \in \mathbb{N}$, $a+b+c = $even.

The calculational scheme for this family starts, very similar to that for the sixth family, with the identity \cite{ade16}:
\begin{eqnarray}
\sum^{\infty}_{k=1} \frac{1}{k^c} \left( \sum^{k}_{i=1} \frac{h^{(a)}_k}{(2k-1)^b} \right) = \sum^{\infty}_{k=1} 
\frac{h^{(a)}_k h^{(b)}_k}{k^c} + \sum^{\infty}_{k=1} \frac{h^{(a+b)}_k}{k^c} -
\sum^{\infty}_{k=1} \frac{1}{k^c} \left( \sum^{k}_{i=1} \frac{h^{(b)}_k}{(2k-1)^a} \right)~.
\end{eqnarray}
From this it follows by a proper rearrangement of the corresponding double summations:
\begin{eqnarray}
\sum^{\infty}_{k=1}\frac{h^{(a)}_k h^{(b)}_k}{k^c} &=& \sum^{\infty}_{k=1} \frac{1}{k^c} \left( \sum^{\infty}_{i=1} \frac{h^{(b)}_i}{(2i-1)^a} \right) +
\sum^{\infty}_{k=1} \frac{1}{k^c} \left( \sum^{\infty}_{i=1} \frac{h^{(a)}_i}{(2i-1)^b} \right) + 
\sum^{\infty}_{k=1} \frac{h^{(a)}_k}{k^c(2k-1)^b} + \sum^{\infty}_{k=1} \frac{h^{(b)}_k}{k^c(2k-1)^a} \nonumber \\ &-&
\sum^{\infty}_{k=1} \frac{h^{(a+b)}_k}{k^c} - \sum^{\infty}_{k=1} \frac{H^{(c)}_k h^{(a)}_k}{(2k-1)^b} - \sum^{\infty}_{k=1} \frac{H^{(c)}_k h^{(b)}_k}{(2k-1)^a}~.
\end{eqnarray}
Again, on the right side of this identity only known Euler sums appear together with Euler sums belonging to the fifth family. 
This way, all members of the eight family can be calculated explicitly in terms of zeta values and linear Euler sums.

\subsection{examples}
The order 4 case:
\begin{eqnarray}
\sum^{\infty}_{k=1} \frac{h_k h_k}{k^2} = \frac{45}{16}\zeta(4)~.
\end{eqnarray}

The order 6 case:
\begin{eqnarray}
\sum^{\infty}_{k=1} \frac{h_k h^{(2)}_k}{k^3} = \frac{49}{8} \zeta(3)^2 - \frac{945}{128}\zeta(6)~,
\end{eqnarray}

\begin{eqnarray}
\sum^{\infty}_{k=1} \frac{h_k h^{(3)}_k}{k^2} = \frac{945}{256}\zeta(6) -  \frac{3}{4}\zeta(2)\sum^{\infty}_{k=1} \frac{h_k}{k^3}~,
\end{eqnarray}

\begin{eqnarray}
\sum^{\infty}_{k=1} \frac{h^{(2)}_k h^{(2)}_k}{k^2} = \frac{945}{128}\zeta(6) - \frac{49}{8}\zeta(3)^2 + \frac{3}{2}\zeta(2)\sum^{\infty}_{k=1} \frac{h_k}{k^3}~.
\end{eqnarray}
This way all members of the eight family to the order 6 have been explicitly calculated.

The order 8 case:
\begin{eqnarray}
\sum^{\infty}_{k=1} \frac{h^{(3)}_k h^{(3)}_k}{k^2} = \frac{2025}{256}\zeta(8) - \frac{21}{8}\zeta(2)\zeta(3)^2~,
\end{eqnarray}

\begin{eqnarray}
\sum^{\infty}_{k=1} \frac{h_k h^{(2)}_k}{k^5} = \frac{651}{8} \zeta(3)\zeta(5) - \frac{343}{16}\zeta(2)\zeta(3)^2 - 
\frac{1575}{32} \zeta(8)~,
\end{eqnarray}

\begin{eqnarray}
\sum^{\infty}_{k=1} \frac{h_k h^{(5)}_k}{k^2} = \frac{2475}{512}\zeta(8) - \frac{21}{64}\zeta(2)\zeta(3)^2 - \frac{15}{16}\zeta(4)\sum^{\infty}_{k=1} \frac{h_k}{k^3} -
\frac{3}{8}\zeta(2)\sum^{\infty}_{k=1} \frac{h_k}{k^5}~.
\end{eqnarray}
The remaining members of this family for the order 8 case can be found in the appendix.

\section{Summary}
We have introduced a special summation method that allows to calculate explicitly eight families of quadratic Euler sums in terms of zeta values and linear Euler sums
of types t(1,2n-1) and s(2,2n) which can be seen as basic numbers like for example $ln(2), \zeta(3)$ or $Li_4(1/2)$. The special linear Euler sums $t$ and $s$ occur pairwise
in the explicit calculation of the corresponding even order nonlinear Euler sums. This means that, for example, in the computation of the order 4 nonlinear Euler sums
in addition to various zeta values the series $t(1,3)$ and $s(2,2)$ are needed for an explicit representation of the corresponding members of the different families.
In the oder 6 case in addition to $t(1,3)$ and $s(2,2)$ the series $t(1,5)$ and $s(2,4)$ appear in the explicit computation of corresponding nonlinear Euler sums.
This calculational scheme continues for higher orders with the additional appearance of $t(1,7)$ and $s(2,6)$ and so on. The quadratic Euler sums themselves consist of
products of even-even, even-odd and odd-odd linear Euler sums where all three possible combinations of even and odd fundamental denominators are considered. This way our
theoretical approach allows for an explicit calculation of a variety of quadratic Euler sums of even order. Even more this approach should be able to open the way for
an explicit computation of the corresponding nonlinear Euler sums of odd order.

\section{Appendix A}
Order 8 case:
\begin{eqnarray}
\sum^{\infty}_{k=1} \frac{H_k h^{(4)}_k}{k^3} &=& 2\sum^{\infty}_{k=1} \frac{h^{(4)}_k}{k^4} + 4\sum^{\infty}_{k=1} \frac{h^{(5)}_k}{k^3}
-\frac{15}{8}\zeta(4)\sum^{\infty}_{k=1} \frac{h_k}{k^3} - \frac{3}{2}\zeta(2)\sum^{\infty}_{k=1} \frac{h^{(3)}_k}{k^3} \nonumber \\ &=&
\frac{1363}{64}\zeta(8) - \frac{39}{8}\zeta(3)\zeta(5) - \frac{273}{32}\zeta(2)\zeta(3)^2 - \frac{15}{8}\zeta(4)\sum^{\infty}_{k=1} \frac{h_k}{k^3} - 
\frac{9}{4}\zeta(2)\sum^{\infty}_{k=1} \frac{h_k}{k^5} - \frac{5}{4}\sum^{\infty}_{k=1} \frac{h_k}{k^7}  \nonumber \\ &+& 
\frac{221}{16}\sum^{\infty}_{k=1} \frac{H^{(2)}_k}{k^6}~.
\end{eqnarray}

\begin{eqnarray}
\sum^{\infty}_{k=1} \frac{H^{(2)}_k h^{(3)}_k}{k^3} &=& -\frac{225}{4}\zeta(8) + \frac{279}{8}\zeta(3)\zeta(5) + \frac{259}{16}\zeta(2)\zeta(3)^2 +
\frac{3}{2}\zeta(2)\sum^{\infty}_{k=1} \frac{h_k}{k^5} - \frac{5}{2}\sum^{\infty}_{k=1} \frac{h^{(3)}_k}{k^5} \nonumber \\ &-& 9\sum^{\infty}_{k=1} \frac{h^{(4)}_k}{k^4}
- 12\sum^{\infty}_{k=1} \frac{h^{(5)}_k}{k^3} \nonumber \\ &=&
\frac{3323}{128}\zeta(8) - \frac{515}{8}\zeta(3)\zeta(5) + \frac{259}{16}\zeta(2)\zeta(3)^2 + \frac{3}{2}\zeta(2)\sum^{\infty}_{k=1}\frac{h_k}{k^5} +
\frac{15}{8}\sum^{\infty}_{k=1} \frac{h_k}{k^7} + \frac{391}{32}\sum^{\infty}_{k=1} \frac{H^{(2)}_k}{k^6}~.
\end{eqnarray}

\begin{eqnarray}
\sum^{\infty}_{k=1} \frac{H^{(3)}_k h^{(2)}_k}{k^3} &=& -\frac{4797}{32}\zeta(8) - \frac{887}{4}\zeta(3)\zeta(5) - \frac{217}{8}\zeta(2)\zeta(3)^2 -
\frac{3}{2}\sum^{\infty}_{k=1} \frac{h_k}{k^7} - \frac{459}{8}\sum^{\infty}_{k=1} \frac{H^{(2)}_k}{k^6}~.
\end{eqnarray}

\begin{eqnarray}
\sum^{\infty}_{k=1} \frac{H^{(4)}_k h_k}{k^3} &=& \frac{2113}{32}\zeta(8) - 93\zeta(3)\zeta(5) + \frac{35}{4}\zeta(2)\zeta(3)^2 +
\zeta(4)\sum^{\infty}_{k=1} \frac{h_k}{k^3} + \frac{1}{2}\sum^{\infty}_{k=1} \frac{h_k}{k^7} + \frac{221}{8}\sum^{\infty}_{k=1} \frac{H^{(2)}_k}{k^6}~. 
\end{eqnarray}

\begin{eqnarray}
\sum^{\infty}_{k=1} \frac{H_k h^{(5)}_k}{k^2} &=& \frac{14335}{512}\zeta(8) - \frac{2859}{64}\zeta(3)\zeta(5) + \frac{189}{32}\zeta(2)\zeta(3)^2 +
\frac{15}{8}\zeta(4)\sum^{\infty}_{k=1} \frac{h_k}{k^3} + \frac{3}{4}\zeta(2)\sum^{\infty}_{k=1} \frac{h_k}{k^5} + \frac{15}{32}\sum^{\infty}_{k=1}
\frac{h_k}{k^7}  \nonumber \\ &+& \frac{1445}{128}\sum^{\infty}_{k=1} \frac{H^{(2)}_k}{k^6}~.
\end{eqnarray}

\begin{eqnarray}
\sum^{\infty}_{k=1} \frac{H^{(2)}_k h^{(4)}_k}{k^2} &=& \frac{12675}{128}\zeta(8) - \frac{357}{16}\zeta(2)\zeta(3)^2 -
\frac{1}{2}\zeta(2)\sum^{\infty}_{k=1} \frac{h_k}{k^5} - \sum^{\infty}_{k=1} \frac{h^{(4)}_k}{k^4} - 8\sum^{\infty}_{k=1} \frac{h^{(5)}_k}{k^3} -
20\sum^{\infty}_{k=1} \frac{h^{(6)}_k}{k^2} \nonumber \\ &=&
-\frac{14723}{128}\zeta(8) + \frac{1449}{8}\zeta(3)\zeta(5) - \frac{357}{16}\zeta(2)\zeta(3)^2 - \frac{1}{2}\zeta(2)\sum^{\infty}_{k=1}\frac{h_k}{k^5} -
\frac{5}{4}\sum^{\infty}_{k=1} \frac{h_k}{k^7} - \frac{833}{16}\sum^{\infty}_{k=1} \frac{H^{(2)}_k}{k^6}~.
\end{eqnarray}

\begin{eqnarray}
\sum^{\infty}_{k=1} \frac{H^{(3)}_k h^{(3)}_k}{k^2} &=& \frac{24447}{128}\zeta(8) - \frac{2251}{8}\zeta(3)\zeta(5) + \frac{133}{4}\zeta(2)\zeta(3)^2 +
\frac{15}{8}\sum^{\infty}_{k=1} \frac{h_k}{k^7} +  \frac{2499}{32}\sum^{\infty}_{k=1} \frac{H^{(6)}_k}{k^2}~.
\end{eqnarray}

\begin{eqnarray}
\sum^{\infty}_{k=1} \frac{H^{(4)}_k h^{(2)}_k}{k^2} &=& \frac{4605}{32}\zeta(8) + 31\zeta(3)\zeta(5) - \frac{77}{4}\zeta(2)\zeta(3)^2 -
\zeta(4)\sum^{\infty}_{k=1} \frac{h_k}{k^3} - 2\sum^{\infty}_{k=1} \frac{h^{(2)}_k}{k^6} - 8\sum^{\infty}_{k=1} \frac{h^{(3)}_k}{k^5} 
\nonumber \\ &-& 18\sum^{\infty}_{k=1} \frac{h^{(4)}_k}{k^4} - 32\sum^{\infty}_{k=1} \frac{h^{(5)}_k}{k^3} - 40\sum^{\infty}_{k=1} \frac{h^{(6)}_k}{k^2}
\nonumber \\ &=&
-\frac{199}{2}\zeta(8) + \frac{615}{4}\zeta(3)\zeta(5) - \frac{77}{4}\zeta(2)\zeta(3)^2 - \zeta(4)\sum^{\infty}_{k=1}\frac{h_k}{k^3} -
\frac{3}{2}\sum^{\infty}_{k=1} \frac{h_k}{k^7} - \frac{323}{8}\sum^{\infty}_{k=1} \frac{H^{(2)}_k}{k^6}~.
\end{eqnarray}

\begin{eqnarray}
\sum^{\infty}_{k=1} \frac{H^{(5)}_k h_k}{k^2} &=& -\frac{855}{16}\zeta(8) -\frac{45}{2}\zeta(3)\zeta(5) + \frac{7}{2}\zeta(2)\zeta(3)^2 +
\frac{7}{2}\sum^{\infty}_{k=1} \frac{h_k}{k^7} + 5\sum^{\infty}_{k=1} \frac{h^{(2)}_k}{k^6} + 8\sum^{\infty}_{k=1} \frac{h^{(3)}_k}{k^5}
\nonumber \\ &+& 12\sum^{\infty}_{k=1} \frac{h^{(4)}_k}{k^4} + 16\sum^{\infty}_{k=1} \frac{h^{(5)}_k}{k^3} + 16\sum^{\infty}_{k=1} \frac{h^{(6)}_k}{k^2}
\nonumber \\ &=&
\frac{301}{16}\zeta(8) - 24\zeta(3)\zeta(5) + \frac{7}{2}\zeta(2)\zeta(3)^2 + \frac{1}{2}\sum^{\infty}_{k=1}\frac{h_k}{k^7} +
\frac{17}{4}\sum^{\infty}_{k=1} \frac{H^{(2)}_k}{k^6}~.
\end{eqnarray}

\section{Appendix B}
\begin{eqnarray}
s(2,8) = \sum^{\infty}_{k=1}\frac{h^{(2)}_k}{k^8} = -\frac{515}{40}\zeta(10) + 16\zeta(3)\zeta(7) + 11\zeta(5)^2 -4\sum^{\infty}_{k=1}\frac{h_k}{k^9} -
\frac{53}{4}\sum^{\infty}_{k=1}\frac{H^{(2)}_k}{k^8}~,
\end{eqnarray}
\begin{eqnarray}
s(3,7) = \sum^{\infty}_{k=1}\frac{h^{(3)}_{k}}{k^7} = \frac{719}{4}\zeta(10) - 122\zeta(3)\zeta(7) - \frac{831}{8}\zeta(5)^2 + 7\sum^{\infty}_{k=1}\frac{h_k}{k^9} +
\frac{293}{4}\sum^{\infty}_{k=1}\frac{H^{(2)}_{k}}{k^8}~,
\end{eqnarray}
\begin{eqnarray}
s(4,6) = \sum^{\infty}_{k=1}\frac{h^{(4)}_{k}}{k^6} = -\frac{28051}{64}\zeta(10) + \frac{4365}{16}\zeta(3)\zeta(7) + \frac{3951}{16}\zeta(5)^2 -
7\sum^{\infty}_{k=1}\frac{h_k}{k^9} - \frac{4757}{32}\sum^{\infty}_{k=1}\frac{H^{(2)}_{k}}{k^8}~,
\end{eqnarray}
\begin{eqnarray}
s(5,5) = \sum^{\infty}_{k=1}\frac{h^{(5)}_{k}}{k^5} = \frac{61185}{128}\zeta(10) - \frac{9425}{32}\zeta(3)\zeta(7) - \frac{16699}{64}\zeta(5)^2 +
\frac{35}{8}\sum^{\infty}_{k=1}\frac{h_k}{k^9} + \frac{9915}{64}\sum^{\infty}_{k=1}\frac{H^{(2)}_{k}}{k^8}~,
\end{eqnarray}
\begin{eqnarray}
s(6,4) = \sum^{\infty}_{k=1}\frac{h^{(6)}_{k}}{k^4} = -\frac{70037}{256}\zeta(10) + \frac{11477}{64}\zeta(3)\zeta(7) + \frac{9079}{64}\zeta(5)^2 -
\frac{7}{4}\sum^{\infty}_{k=1}\frac{h_k}{k^9} - \frac{11869}{128}\sum^{\infty}_{k=1}\frac{H^{(2)}_{k}}{k^8}
\end{eqnarray}
and
\begin{eqnarray}
s(7,3) = \sum^{\infty}_{k=1}\frac{h^{(7)}_{k}}{k^3} = \frac{37247}{512}\zeta(10) - \frac{3409}{64}\zeta(3)\zeta(7) - \frac{2293}{64}\zeta(5)^2 +
\frac{7}{16}\sum^{\infty}_{k=1}\frac{h_k}{k^9} + \frac{7903}{256}\sum^{\infty}_{k=1}\frac{H^{(2)}_{k}}{k^8}~.
\end{eqnarray}

\section{Appendix C}
Here we present more nonlinear Euler sums belonging to different families and orders.

Family 2, order eight case:
\begin{eqnarray}
\sum^{\infty}_{k=1} \frac{H_k^2}{(2k-1)^6} &=& \frac{4847}{256}\zeta(8) - \frac{381}{16}ln(2)\zeta(7) + \frac{381}{16}\zeta(7) + \frac{63}{16}\left(ln(2)\right)^2\zeta(6)
-\frac{63}{8}ln(2)\zeta(6) - \frac{189}{32}\zeta(6) - \frac{351}{16}\zeta(3)\zeta(5) \nonumber \\ &+& \frac{93}{16}ln(2)\zeta(2)\zeta(5) - \frac{93}{16}\zeta(2)\zeta(5)
+ \frac{31}{4}ln(2)\zeta(5) + \frac{31}{4}\zeta(5) + \frac{105}{16}ln(2)\zeta(3)\zeta(4) - \frac{105}{16}\zeta(3)\zeta(4) \nonumber \\ &-& \frac{15}{2}ln(2)\zeta(4) +
\frac{45}{8}\zeta(4) + \frac{147}{64}\zeta(2)\zeta(3)^2 + \frac{49}{16}\zeta(3)^2 - \frac{21}{4}\zeta(2)\zeta(3) + 7ln(2)\zeta(3) - 7\zeta(3) - 6ln(2)\zeta(2) 
\nonumber \\ &+& 14\zeta(2) - 20ln(2) + 4\left(ln(2)\right)^2 + \frac{3}{16}\sum^{\infty}_{k=1}\frac{h_k}{k^7} + \frac{289}{64}\sum^{\infty}_{k=1}\frac{H^{(2)}_{k}}{k^6}~.
\end{eqnarray}

\begin{eqnarray}
\sum^{\infty}_{k=1} \frac{H_k H^{(2)}_k}{(2k-1)^5} &=& -\frac{25285}{512}\zeta(8) + \frac{1905}{32}ln(2)\zeta(7) - \frac{1905}{32}\zeta(7) + \frac{705}{32}\zeta(6) + 
\frac{2189}{32}\zeta(3)\zeta(5) - \frac{403}{16}ln(2)\zeta(2)\zeta(5) \nonumber \\ &+& \frac{403}{16}\zeta(2)\zeta(5) -  \frac{31}{4}ln(2)\zeta(5) - \frac{93}{2}\zeta(5) -
\frac{105}{8}ln(2)\zeta(3)\zeta(4) + \frac{105}{8}\zeta(3)\zeta(4) + 15ln(2)\zeta(4) - \frac{15}{4}\zeta(4)  \nonumber \\ &-& \frac{637}{64}\zeta(2)\zeta(3)^2 -
\frac{175}{16}\zeta(3)^2 + \frac{91}{4}\zeta(2)\zeta(3) - 21ln(2)\zeta(3) + \frac{21}{2}\zeta(3) + 26ln(2)\zeta(2) - 60\zeta(2) \nonumber \\ &+& 120ln(2) -
20\left(ln(2)\right)^2 - 2\sum^{\infty}_{k=1} \frac{h_k}{k^3} - \sum^{\infty}_{k=1} \frac{h_k}{k^5} - \frac{15}{32}\sum^{\infty}_{k=1} \frac{h_k}{k^7}
- \frac{1445}{128}\sum^{\infty}_{k=1} \frac{H^{(2)}_k}{k^6}~.
\end{eqnarray}

\begin{eqnarray}
\sum^{\infty}_{k=1} \frac{H_k H^{(3)}_k}{(2k-1)^4} &=& -\frac{203}{16}\zeta(8) - \frac{635}{4}ln(2)\zeta(7) + \frac{635}{4}\zeta(7) - \frac{135}{8}\zeta(6) - 
\frac{147}{8}\zeta(3)\zeta(5) + \frac{341}{4}ln(2)\zeta(2)\zeta(5) \nonumber \\ &-& \frac{341}{4}\zeta(2)\zeta(5) + 62\zeta(5) +
\frac{45}{4}ln(2)\zeta(3)\zeta(4) - \frac{45}{4}\zeta(3)\zeta(4) - 15ln(2)\zeta(4) - \frac{75}{4}\zeta(4)  \nonumber \\ &+& \frac{161}{16}\zeta(2)\zeta(3)^2 +
\frac{77}{8}\zeta(3)^2 - \frac{53}{2}\zeta(2)\zeta(3) + 40ln(2)\zeta(3) + 44\zeta(3) - 72ln(2)\zeta(2) + 120\zeta(2) \nonumber \\ &-& 320ln(2) +
80\left(ln(2)\right)^2 + 2\sum^{\infty}_{k=1} \frac{h_k}{k^3} + 3\sum^{\infty}_{k=1} \frac{h_k}{k^5} - \frac{1}{2}\zeta(2)\sum^{\infty}_{k=1} \frac{h_k}{k^5} 
+ \frac{5}{4}\sum^{\infty}_{k=1} \frac{h_k}{k^7} - \frac{221}{16}\sum^{\infty}_{k=1} \frac{H^{(2)}_k}{k^6}~.
\end{eqnarray}
In order to calculate explicitly Eq.~(283) the double valued help function
\begin{eqnarray}
\sum^{\infty}_{k=1} \frac{H^{(3)}_k}{(i+k)(2k-1)} &=& \frac{1}{2i+1} \sum^{\infty}_{k=1}\frac{H^{(3)}_k}{k(2k-1)} + \frac{\zeta(4)}{2i+1} + \zeta(3)\frac{H_{i-1}}{2i+1} 
-\zeta(2)\frac{H^{(2)}_{i-1}}{2i+1} + \frac{1}{2i+1}\sum^{i-1}_{k=1}\frac{H_k}{k^3}
\end{eqnarray}
is needed. This identity represents the generalization of Eq.~(118) to the order three case. Analogously it follows:
\begin{eqnarray}
\sum^{\infty}_{k=1} \frac{H_k H^{(4)}_k}{(2k-1)^3} &=& -\frac{9011}{128}\zeta(8) + \frac{1905}{8}ln(2)\zeta(7) - \frac{1905}{8}\zeta(7) + \frac{69}{8}\zeta(6) - 
\frac{1605}{16}\zeta(3)\zeta(5) - \frac{279}{2}ln(2)\zeta(2)\zeta(5) \nonumber \\ &+& \frac{279}{2}\zeta(2)\zeta(5) - \frac{31}{2}\zeta(5) -
\frac{7}{4}ln(2)\zeta(3)\zeta(4) + \frac{7}{4}\zeta(3)\zeta(4) + 2ln(2)\zeta(4) + 15\zeta(4)  \nonumber \\ &-& \frac{7}{16}\zeta(2)\zeta(3)^2 -
\frac{7}{2}\zeta(3)^2 + 7\zeta(2)\zeta(3) - 28ln(2)\zeta(3) - 132\zeta(3) + 96ln(2)\zeta(2) - 120\zeta(2) \nonumber \\ &+& 480ln(2) -
160\left(ln(2)\right)^2 - 4\sum^{\infty}_{k=1} \frac{h_k}{k^5} + \frac{3}{2}\zeta(2)\sum^{\infty}_{k=1} \frac{h_k}{k^5} 
- \frac{15}{8}\sum^{\infty}_{k=1} \frac{h_k}{k^7} + \frac{1717}{32}\sum^{\infty}_{k=1} \frac{H^{(2)}_k}{k^6}~.
\end{eqnarray}

Family 3, order eight case:
\begin{eqnarray}
\sum^{\infty}_{k=1} \frac{H^{(2)}_k H^{(5)}_k}{k(2k-1)} &=& \frac{12931}{144}\zeta(8) + \frac{25}{3}\zeta(6) - 118\zeta(3)\zeta(5) + 2ln(2)\zeta(2)\zeta(5) -
2\zeta(2)\zeta(5) + 8ln(2)\zeta(5) \nonumber \\ &+& 36\zeta(5) - 28\zeta(4) + \frac{11}{2}\zeta(2)\zeta(3)^2 - 6\zeta(3)^2 - 24\zeta(2)\zeta(3) - 80\zeta(3) +
64ln(2)\zeta(2) \nonumber \\ &-& 64\zeta(2) + 256ln(2) - 8\sum^{\infty}_{k=1} \frac{h_k}{k^5} - 2\zeta(2)\sum^{\infty}_{k=1} \frac{h_k}{k^5} + \frac{173}{4}
\sum^{\infty}_{k=1} \frac{H^{(2)}_k}{k^6}~.
\end{eqnarray}

Family 5, order eight case:
\begin{eqnarray}
\sum^{\infty}_{k=1} \frac{H^{(2)}_k h_k}{(2k-1)^5} &=& -\frac{363}{1024}\zeta(8) - \frac{1905}{64}ln(2)\zeta(7) + \frac{63}{32}\zeta(6) - \frac{101}{64}\zeta(3)\zeta(5)
+ \frac{403}{32}ln(2)\zeta(2)\zeta(5) + \frac{31}{8}ln(2)\zeta(5) \nonumber \\ &-& \frac{31}{5}\zeta(5) + \frac{105}{16}ln(2)\zeta(3)\zeta(4) -
\frac{15}{2}ln(2)\zeta(4) + \frac{45}{8}\zeta(4) - \frac{49}{64}\zeta(2)\zeta(3)^2 - \frac{7}{16}\zeta(3)^2 + \frac{3}{4}\zeta(2)\zeta(3) \nonumber \\ &+&
\frac{21}{2}ln(2)\zeta(3) - \frac{35}{4}\zeta(3) - 12ln(2)\zeta(2) + 10\zeta(2) - \frac{3}{2}\sum^{\infty}_{k=1}\frac{h_k}{k^3} - \frac{1}{8}\sum^{\infty}_{k=1}\frac{h_k}{k^5}
- \frac{1}{32}\zeta(2)\sum^{\infty}_{k=1}\frac{h_k}{k^5} \nonumber \\ &+& \frac{21}{64}\sum^{\infty}_{k=1}\frac{h_k}{k^7} +
\frac{969}{256}\sum^{\infty}_{k=1}\frac{H^{(2)}_k}{k^6}~.
\end{eqnarray}

\begin{eqnarray}
\sum^{\infty}_{k=1} \frac{H_k h_k}{(2k-1)^6} &=& \frac{1889}{1024}\zeta(8)+\frac{635}{128}ln(2)\zeta(7)+\frac{127}{128}\zeta(7)-\frac{63}{32}\left(ln(2)\right)^2\zeta(6) 
+ \frac{63}{32}ln (2)\zeta(6) -\frac{63}{64}\zeta(6) - \frac{349}{256}\zeta(3)\zeta(5) \nonumber \\ &-& \frac{45}{32}ln(2)\zeta(2)\zeta(5) - \frac{3}{64}\zeta(2)\zeta(5)
- \frac{31}{16}ln(2)\zeta(5) +\frac{31}{32}\zeta(5) - \frac{45}{32}ln(2)\zeta(3)\zeta(4) - \frac{15}{64}\zeta(3)\zeta(4)  \nonumber \\ &+& \frac{15}{8}ln(2)\zeta(4) -
\frac{15}{16}\zeta(4) + \frac{21}{128}\zeta(2)\zeta(3)^2 + \frac{7}{32}\zeta(3)^2 - \frac{3}{16}\zeta(2)\zeta(3) - \frac{7}{4}ln(2)\zeta(3) + \frac{7}{8}\zeta(3)
\nonumber \\ &+& \frac{3}{2}ln(2)\zeta(2) - \zeta(2) + \frac{1}{4} \sum^{\infty}_{k=1} \frac{h_k}{k^3} + \frac{1}{16} \sum^{\infty}_{k=1} \frac{h_k}{k^5} -
\frac{3}{64} \sum^{\infty}_{k=1}\frac{h_k}{k^7} - \frac{17}{256}\sum^{\infty}_{k=1}\frac{H^{(2)}_k}{k^6}~.
\end{eqnarray}

\begin{eqnarray}
\sum^{\infty}_{k=1} \frac{H^{(2)}_k h^{(3)}_k}{(2k-1)^3} &=& -\frac{2743}{512}\zeta(8) + \frac{63}{32}\zeta(6) - \frac{51}{64}\zeta(3)\zeta(5)
- \frac{31}{2}\zeta(5) + \frac{225}{16}\zeta(4) + \frac{217}{128}\zeta(2)\zeta(3)^2 \nonumber \\ &+& \frac{49}{32}\zeta(3)^2 + \frac{9}{4}\zeta(2)\zeta(3) -
\frac{21}{2}ln(2)\zeta(3) + \frac{3}{2}\sum^{\infty}_{k=1}\frac{h_k}{k^3} + \frac{3}{32}\sum^{\infty}_{k=1}\frac{h_k}{k^7} + 
\frac{289}{128}\sum^{\infty}_{k=1}\frac{H^{(2)}_k}{k^6}~.
\end{eqnarray}

\begin{eqnarray}
\sum^{\infty}_{k=1} \frac{H^{(5)}_k h_k}{(2k-1)^2} &=& -\frac{2307}{128}\zeta(8)+\frac{381}{4}ln(2)\zeta(7) + \frac{395}{16}\zeta(3)\zeta(5) - \frac{107}{2}\zeta(2)\zeta(5)
+ 31\zeta(5) - \frac{7}{2}ln(2)\zeta(3)\zeta(4) \nonumber \\ &+& \frac{63}{8}\zeta(2)\zeta(3)^2 - 14\zeta(2)\zeta(3) + 70\zeta(3)  
+ 24ln(2)\zeta(2) - 80\zeta(2) + 12\sum^{\infty}_{k=1} \frac{h_k}{k^3} \nonumber \\ &+& \frac{1}{2}\zeta(4)\sum^{\infty}_{k=1} \frac{h_k}{k^3} + \sum^{\infty}_{k=1} \frac{h_k}{k^5}
+  \zeta(2)\sum^{\infty}_{k=1} \frac{h_k}{k^5} - \frac{21}{8} \sum^{\infty}_{k=1}\frac{h_k}{k^7} - \frac{969}{32}\sum^{\infty}_{k=1}\frac{H^{(2)}_k}{k^6}~.
\end{eqnarray}

\begin{eqnarray}
\sum^{\infty}_{k=1} \frac{H^{(4)}_k h^{(2)}_k}{(2k-1)^2} &=& -\frac{89}{128}\zeta(8) + \frac{135}{24}\zeta(6) + \frac{163}{8}\zeta(3)\zeta(5) - 62\zeta(5)
+ \frac{105}{2}\zeta(4) - \frac{287}{16}\zeta(2)\zeta(3)^2 - \frac{7}{4}\zeta(3)^2  \nonumber \\ &+& 35\zeta(2)\zeta(3) - 84\zeta(3)  
+ 48ln(2)\zeta(2) - 12\sum^{\infty}_{k=1} \frac{h_k}{k^3} - 2\sum^{\infty}_{k=1} \frac{h_k}{k^5} - \frac{3}{2}\zeta(2)\sum^{\infty}_{k=1} \frac{h_k}{k^5}
+ \frac{5}{4} \sum^{\infty}_{k=1}\frac{h_k}{k^7} \nonumber \\ &+& \frac{153}{8}\sum^{\infty}_{k=1}\frac{H^{(2)}_k}{k^6}~.
\end{eqnarray}

\begin{eqnarray}
\sum^{\infty}_{k=1} \frac{H_k h^{(3)}_k}{(2k-1)^4} &=& -\frac{883}{512}\zeta(8) - \frac{127}{128}ln(2)\zeta(7) + \frac{127}{128}\zeta(7) - \frac{63}{64}\zeta(6) +
\frac{1411}{256}\zeta(3)\zeta(5) + \frac{15}{32}ln(2)\zeta(2)\zeta(5) \nonumber \\ &-& \frac{31}{32}\zeta(2)\zeta(5) + \frac{31}{32}\zeta(5)- \frac{45}{32}ln(2)\zeta(3)\zeta(4) +
\frac{45}{32}\zeta(3)\zeta(4) - \frac{75}{32}\zeta(4) - \frac{105}{128}\zeta(2)\zeta(3)^2  \nonumber \\ &-& \frac{49}{64}\zeta(3)^2 + \frac{3}{4}\zeta(2)\zeta(3)+
\frac{7}{4}ln(2)\zeta(3) - \frac{1}{4}\sum^{\infty}_{k=1} \frac{h_k}{k^3} - \frac{221}{128}\sum^{\infty}_{k=1}\frac{H^{(2)}_k}{k^6}~.
\end{eqnarray}
In order to calculate explicitly Eq.~(292) the double valued help function
\begin{eqnarray}
\sum^{\infty}_{k=1} \frac{H_k}{(2i+2k-1)^3} &=& \sum^{\infty}_{k=1}\frac{H_k}{(2k+1)^3} + \frac{7}{4}\zeta(3) h_{i-1} + \frac{3}{2}\zeta(2)h^{(2)}_{i-1} +
2ln(2)h^{(3)}_{i-1} - 2h_{i-1}h^{(3)}_{i-1} -2h^{(4)}_{i-1} \nonumber \\ &-& 2\sum^{i-1}_{k=1}\frac{h^{(2)}_k}{(2k-1)^2}
\end{eqnarray}
is needed. Analogously it follows:
\begin{eqnarray}
\sum^{\infty}_{k=1} \frac{H_k h^{(4)}_k}{(2k-1)^3} &=& \frac{3793}{512}\zeta(8) - \frac{127}{128}ln(2)\zeta(7) + \frac{127}{128}\zeta(7) - \frac{351}{128}\zeta(6) -
\frac{1845}{256}\zeta(3)\zeta(5) - \frac{15}{32}ln(2)\zeta(2)\zeta(5) \nonumber \\ &+& \frac{15}{32}\zeta(2)\zeta(5) + \frac{155}{32}\zeta(5) - \frac{15}{64}ln(2)\zeta(3)\zeta(4) +
\frac{15}{64}\zeta(3)\zeta(4) - \frac{15}{8}ln(2)\zeta(4) \nonumber \\ &+& \frac{105}{256}\zeta(2)\zeta(3)^2 + \frac{7}{32}\zeta(3)^2 -
\frac{9}{8}\zeta(2)\zeta(3) + \frac{15}{64}\zeta(4)\sum^{\infty}_{k=1}\frac{h_k}{k^3} + \frac{221}{128}\sum^{\infty}_{k=1}\frac{H^{(2)}_k}{k^6}~.
\end{eqnarray}
 
\begin{eqnarray}
\sum^{\infty}_{k=1} \frac{H^{(2)}_k h^{(2)}_k}{(2k-1)^4} &=& \frac{6319}{512}\zeta(8) + \frac{27}{8}\zeta(6) - \frac{1291}{64}\zeta(3)\zeta(5) - \frac{31}{8}\zeta(5)
+ \frac{135}{8}\zeta(4) + \frac{35}{32}\zeta(2)\zeta(3)^2 + \frac{7}{16}\zeta(3)^2 \nonumber \\ &-& \frac{9}{4}\zeta(2)\zeta(3) -21\zeta(3) + 12ln(2)\zeta(2) 
- \sum^{\infty}_{k=1}\frac{h_k}{k^3} - \frac{3}{8}\zeta(2)\sum^{\infty}_{k=1}\frac{h_k}{k^5} + \frac{3}{32}\sum^{\infty}_{k=1}\frac{h_k}{k^7} 
\nonumber \\ &+& \frac{1343}{128}\sum^{\infty}_{k=1}\frac{H^{(2)}_k}{k^6}~.
\end{eqnarray}
In order to calculate explicitly Eq.~(295) the double valued help function
\begin{eqnarray}
\sum^{\infty}_{k=1} \frac{h_k}{(i+k)^3} &=& \sum^{\infty}_{k=1}\frac{h_k}{(k+1)^3} + 8ln(2) - 2\zeta(2) - \zeta(3) + \zeta(3)h_i + 2\zeta(2)h^{(2)}_i - 8ln(2)h^{(3)}_i
- 4H_i h^{(3)}_i \nonumber \\ &-& 2 H^{(2)}_i h^{(2)}_i - H^{(3)}_i h_i + 4\sum^i_{k=1}\frac{h^{(3)}_k}{k} + 2\sum^i_{k=1}\frac{h^{(2)}_k}{k^2} + 
\sum^i_{k=1}\frac{h_k}{k^3}
\end{eqnarray}
is needed. Analogously it follows:
\begin{eqnarray}
\sum^{\infty}_{k=1} \frac{H^{(2)}_k h^{(4)}_k}{(2k-1)^2} &=& -\frac{6615}{256}\zeta(8) + \frac{423}{32}\zeta(6) + \frac{2057}{64}\zeta(3)\zeta(5) - \frac{155}{8}\zeta(5)
+ \frac{15}{2}ln(2)\zeta(4) - \frac{175}{64}\zeta(2)\zeta(3)^2 - \frac{35}{8}\zeta(3)^2 \nonumber \\ &+& \frac{9}{2}\zeta(2)\zeta(3) - 
\frac{15}{16}\zeta(4)\sum^{\infty}_{k=1}\frac{h_k}{k^3} - \frac{1}{2}\sum^{\infty}_{k=1}\frac{h_k}{k^5} + \frac{3}{32}\sum^{\infty}_{k=1}\frac{h_k}{k^7} 
- \frac{765}{128}\sum^{\infty}_{k=1}\frac{H^{(2)}_k}{k^6}~.
\end{eqnarray}

\begin{eqnarray}
\sum^{\infty}_{k=1} \frac{H^{(4)}_k h_k}{(2k-1)^3} &=& \frac{12597}{256}\zeta(8) - \frac{1905}{16}ln(2)\zeta(7) - \frac{2207}{32}\zeta(3)\zeta(5) 
+ \frac{279}{4}ln(2)\zeta(2)\zeta(5) - \frac{31}{4}\zeta(5) + \frac{7}{8}ln(2)\zeta(3)\zeta(4) \nonumber \\ &+& \frac{15}{2}\zeta(4) - \frac{91}{32}\zeta(2)\zeta(3)^2
+ \frac{7}{2}\zeta(2)\zeta(3) + 14ln(2)\zeta(3) - 70\zeta(3) -  48ln(2)\zeta(2) + 80\zeta(2) - 8\sum^{\infty}_{k=1}\frac{h_k}{k^3} \nonumber \\ &-& 
\frac{1}{8}\zeta(4)\sum^{\infty}_{k=1}\frac{h_k}{k^3} - \frac{3}{4}\zeta(4)\sum^{\infty}_{k=1}\frac{h_k}{k^5} + \frac{35}{16}\sum^{\infty}_{k=1}\frac{h_k}{k^7} 
+ \frac{2669}{64}\sum^{\infty}_{k=1}\frac{H^{(2)}_k}{k^6}~.
\end{eqnarray}
In order to calculate explicitly Eq.~(298) the double valued help function defined by Eq.~(171) is needed.
\begin{eqnarray}
\sum^{\infty}_{k=1} \frac{H^{(3)}_k h_k}{(2k-1)^4} &=& -\frac{11257}{512}\zeta(8) + \frac{635}{8}ln(2)\zeta(7) + \frac{1057}{32}\zeta(3)\zeta(5) 
- \frac{341}{8}ln(2)\zeta(2)\zeta(5) + \frac{31}{8}\zeta(5) - \frac{45}{8}ln(2)\zeta(3)\zeta(4) \nonumber \\ &+& \frac{15}{2}ln(2)\zeta(4) - 
\frac{45}{4}\zeta(4) + \frac{39}{32}\zeta(2)\zeta(3)^2 - \frac{3}{4}\zeta(2)\zeta(3) - 21ln(2)\zeta(3) + 35\zeta(3) + 36ln(2)\zeta(2) \nonumber \\ &-& 40\zeta(2) + 
4\sum^{\infty}_{k=1}\frac{h_k}{k^3} + \frac{1}{4}\zeta(4)\sum^{\infty}_{k=1}\frac{h_k}{k^5} - \frac{35}{32}\sum^{\infty}_{k=1}\frac{h_k}{k^7} 
- \frac{2669}{128}\sum^{\infty}_{k=1}\frac{H^{(2)}_k}{k^6}~.
\end{eqnarray}
In order to calculate explicitly Eq.~(299) the double valued help function
\begin{eqnarray}
\sum^{\infty}_{k=1} \frac{H^{(3)}_k}{(2k-1)(2i+2k-1)} &=& \sum^{\infty}_{k=1}\frac{1}{k^3(2k-1)(2i+2k-1)} + \frac{1}{2}\zeta(3)\frac{h_i}{i} -
\zeta(2)\frac{h^{(2)}_i}{i} - 4ln(2)\frac{h^{(3)}_i}{i} \nonumber \\ &+& \frac{4}{i}\sum^i_{k=1}\frac{h_k}{(2k-1)^3}
\end{eqnarray}
is needed.
\begin{eqnarray}
\sum^{\infty}_{k=1} \frac{H^{(3)}_k h^{(2)}_k}{(2k-1)^3} &=& -\frac{17909}{512}\zeta(8) + \frac{1383}{32}\zeta(3)\zeta(5) + \frac{155}{8}\zeta(5) - 45\zeta(4) + 
\frac{221}{64}\zeta(2)\zeta(3)^2 - \frac{13}{2}\zeta(2)\zeta(3) \nonumber \\ &+& 63\zeta(3) - 36ln(2)\zeta(2) + 6\sum^{\infty}_{k=1}\frac{h_k}{k^3} + 
\frac{9}{8}\zeta(2)\sum^{\infty}_{k=1}\frac{h_k}{k^5} - \frac{15}{32}\sum^{\infty}_{k=1}\frac{h_k}{k^7} - \frac{3553}{128}\sum^{\infty}_{k=1}\frac{H^{(2)}_k}{k^6}~.
\end{eqnarray}
In order to calculate explicitly Eq.~(301) the double valued help function
\begin{eqnarray}
\sum^{\infty}_{k=1} \frac{H^{(2)}_k}{(2i+2k-1)^3} &=& \sum^{\infty}_{k=1}\frac{H^{(2)}_k}{(2k+1)^3} - \frac{7}{2}\zeta(3)h^{(2)}_{i-1} - 7\zeta(2)h^{(3)}_{i-1} -
12ln(2)h^{(4)}_{i-1} + 12\sum^{i-1}_{k=1}\frac{h_k}{(2k-1)^4}  \nonumber \\ &+& 8\sum^{i-1}_{k=1}\frac{h^{(2)}_k}{(2k-1)^3} + 4\sum^{i-1}_{k=1}\frac{h^{(3)}_k}{(2k-1)^2}
\end{eqnarray}
is needed. Analogously it follows:

\begin{eqnarray}
\sum^{\infty}_{k=1} \frac{H^{(3)}_k h^{(3)}_k}{(2k-1)^2} &=& \frac{20739}{512}\zeta(8) - \frac{135}{16}\zeta(6) - \frac{1351}{32}\zeta(3)\zeta(5) + \frac{403}{8}\zeta(5)
- \frac{225}{8}\zeta(4) + \frac{97}{16}\zeta(2)\zeta(3)^2 + \frac{91}{16}\zeta(3)^2 \nonumber \\ &-& 18\zeta(2)\zeta(3) + 21ln(2)\zeta(3) - 3\sum^{\infty}_{k=1}\frac{h_k}{k^3} + 
+ \frac{3}{2}\sum^{\infty}_{k=1}\frac{h_k}{k^5} - \frac{15}{32}\sum^{\infty}_{k=1}\frac{h_k}{k^7} + \frac{663}{128}\sum^{\infty}_{k=1}\frac{H^{(2)}_k}{k^6}~.
\end{eqnarray}

Family 6, order eight case:
\begin{eqnarray}
\sum^{\infty}_{k=1} \frac{h^{(3)}_k h^{(4)}_k}{k(2k-1)} &=& \frac{10321}{6144}\zeta(8) + \frac{407}{256}\zeta(3)\zeta(5) - \frac{105}{64}ln(2)\zeta(3)\zeta(4)
- \frac{9}{32}\zeta(2)\zeta(3)^2 + \frac{15}{64}\zeta(4)\sum^{\infty}_{k=1} \frac{h_k}{k^3} \nonumber \\ &-& \frac{221}{512}\sum^{\infty}_{k=1} \frac{H^{(2)}_k}{k^6}~.
\end{eqnarray}

Family 7, order eight case:
\begin{eqnarray}
\sum^{\infty}_{k=1} \frac{H^{(3)}_k h^{(4)}_k}{k(2k-1)} &=& \frac{28481}{256}\zeta(8) - \frac{635}{2}ln(2)\zeta(7) + \frac{635}{4}\zeta(7) - \frac{495}{16}\zeta(6)
- \frac{1793}{16}\zeta(3)\zeta(5) + \frac{341}{2}ln(2)\zeta(2)\zeta(5) \nonumber \\ &-&  \frac{341}{4}\zeta(2)\zeta(5) + \frac{155}{4}\zeta(5) +
\frac{195}{8}ln(2)\zeta(3)\zeta(4) - \frac{105}{8}\zeta(3)\zeta(4) - 15ln(2)\zeta(4) - \frac{1}{4}\zeta(2)\zeta(3)^2 \nonumber \\ &+& \frac{63}{4}\zeta(3)^2 - 9\zeta(2)\zeta(3) 
- \frac{15}{8}\zeta(4)\sum^{\infty}_{k=1} \frac{h_k}{k^3} + 2\sum^{\infty}_{k=1} \frac{h_k}{k^5} - \frac{11}{2}\zeta(2)\sum^{\infty}_{k=1} \frac{h_k}{k^5}
+ \frac{45}{16}\sum^{\infty}_{k=1} \frac{h_k}{k^7} \nonumber \\ &+& \frac{2227}{64}\sum^{\infty}_{k=1} \frac{H^{(2)}_k}{k^6}~.
\end{eqnarray}

\begin{eqnarray}
\sum^{\infty}_{k=1} \frac{H^{(4)}_k h^{(3)}_k}{k(2k-1)} &=& -\frac{17109}{128}\zeta(8) + \frac{1905}{4}ln(2)\zeta(7) - \frac{1905}{8}\zeta(7) + \frac{135}{4}\zeta(6)
- \frac{243}{2}\zeta(3)\zeta(5) - 279ln(2)\zeta(2)\zeta(5) \nonumber \\ &+&  \frac{279}{2}\zeta(2)\zeta(5) - 93\zeta(5) -
\frac{7}{4}ln(2)\zeta(3)\zeta(4) + \frac{75}{2}\zeta(4) - \frac{91}{4}\zeta(3)^2 + 42\zeta(2)\zeta(3) - 28ln(2)\zeta(3) \nonumber \\ &+&
4\sum^{\infty}_{k=1} \frac{h_k}{k^3} + \frac{1}{4}\zeta(4)\sum^{\infty}_{k=1}\frac{h_k}{k^3} - 6\sum^{\infty}_{k=1} \frac{h_k}{k^5} + 9\zeta(2)\sum^{\infty}_{k=1} \frac{h_k}{k^5}
- \frac{25}{8}\sum^{\infty}_{k=1} \frac{h_k}{k^7} - \frac{1003}{32}\sum^{\infty}_{k=1} \frac{H^{(2)}_k}{k^6}~.
\end{eqnarray}

Family 8, order eight case:
\begin{eqnarray}
\sum^{\infty}_{k=1} \frac{h^{(2)}_k h^{(3)}_k}{k^3} = -\frac{3375}{128}\zeta(8) + \frac{651}{16}\zeta(3)\zeta(5) - \frac{567}{64}\zeta(2)\zeta(3)^2 -
\frac{9}{8}\zeta(2)\sum^{\infty}_{k=1}\frac{h_k}{k^5}~.
\end{eqnarray}

\begin{eqnarray}
\sum^{\infty}_{k=1} \frac{h^{(2)}_k h^{(2)}_k}{k^4} = \frac{2925}{64}\zeta(8) - \frac{217}{2}\zeta(3)\zeta(5) + \frac{287}{8}\zeta(2)\zeta(3)^2 +
3\zeta(2)\sum^{\infty}_{k=1}\frac{h_k}{k^5}~.
\end{eqnarray}

\begin{eqnarray}
\sum^{\infty}_{k=1} \frac{h_k h^{(3)}_k}{k^4} = \frac{2475}{64}\zeta(8) - \frac{651}{16}\zeta(3)\zeta(5) + \frac{105}{16}\zeta(2)\zeta(3)^2 -
\frac{3}{2}\zeta(2)\sum^{\infty}_{k=1}\frac{h_k}{k^5}~.
\end{eqnarray}

\begin{eqnarray}
\sum^{\infty}_{k=1} \frac{h_k h^{(4)}_k}{k^3} = -\frac{4275}{256}\zeta(8) + \frac{217}{16}\zeta(3)\zeta(5) - \frac{21}{64}\zeta(2)\zeta(3)^2 +
\frac{9}{8}\zeta(2)\sum^{\infty}_{k=1}\frac{h_k}{k^5}~.
\end{eqnarray}

\section{Appendix D}
Here we present the numerical computation of the 'corrector functions $O2n_1$ and $O2n_2$, n$\in \mathbb{N}$, as well as the numerical computation of selected nonlinear Euler
sums of different families. It follows: 
\begin{eqnarray}
\sum^{\infty}_{k=1} \frac{h_k}{k^3} = 1.29817551577186712475~,
\end{eqnarray}
\begin{eqnarray}
\sum^{\infty}_{k=1} \frac{h_k}{k^5} = 1.05070828839842708728~,
\end{eqnarray}
\begin{eqnarray}
\sum^{\infty}_{k=1} \frac{h_k}{k^7} = 1.01125393384735909982~,
\end{eqnarray}
\begin{eqnarray}
\sum^{\infty}_{k=1} \frac{h_k}{k^9} = 1.00268963654979441867~,
\end{eqnarray}
\begin{eqnarray}
\sum^{\infty}_{k=1} \frac{H^{(2)}_k}{k^6} = 1.02189709661478032774~,
\end{eqnarray}
\begin{eqnarray}
\sum^{\infty}_{k=1} \frac{H^{(2)}_k}{k^8} = 1.00511704480621791756~,
\end{eqnarray}

\begin{eqnarray}
\sum^{\infty}_{k=1} \frac{H_k h^{(3)}_k }{k^4} = 1.13880957043316224083~,
\end{eqnarray}
\begin{eqnarray}
\sum^{\infty}_{k=1} \frac{H^{(2)}_k H^{(3)}_k }{(2k-1)^3} = 1.07653048082931893543~,
\end{eqnarray}
\begin{eqnarray}
\sum^{\infty}_{k=1} \frac{h^{(2)}_k h^{(3)}_k }{(2k-1)^3} =  1.06071267806165873612~,
\end{eqnarray}
\begin{eqnarray}
\sum^{\infty}_{k=1} \frac{H_k h^{(2)}_k }{(2k-1)^5} = 1.00776268435455472628~,
\end{eqnarray}
\begin{eqnarray}
\sum^{\infty}_{k=1} \frac{h_k h^{(3)}_k }{k^4} = 1.12021909420689660746~,
\end{eqnarray}
where the numerical accuracy has always been chosen to $\delta = 10^{-20}$.
\begin{eqnarray}
\sum^{\infty}_{k=1} \frac{H^{(3)}_k h^{(4)}_k }{k(2k-1)} = 1.479405080355585 
\end{eqnarray}
For the last Euler series we were able to achieve a numerical accuracy of $\delta = 10^{-15}$ only as the summation procedure converges very slowly
for this series. In fact we had to consider $10^{13}$ terms in the numerical calculation.

\end{document}